\documentclass[pdflatex,sn-mathphys-num]{sn-jnl}% Math and Physical Sciences Numbered Reference Style
\usepackage{amssymb}

\usepackage{mathtools}

\usepackage{color,xcolor,ucs}
\usepackage{mathtools}   \usepackage{tikz} 
\usepackage{ amssymb }
\usepackage{extarrows} 
\usepackage{pgf,tikz}
\usepackage{float}
\usetikzlibrary{positioning}
\usetikzlibrary{shapes.geometric}
\usetikzlibrary{shapes.misc}
\usetikzlibrary{arrows}
\usepackage{caption}
\usepackage{mathrsfs}
\usetikzlibrary{arrows,shapes,automata,backgrounds,petri,positioning}
\usetikzlibrary{decorations.pathmorphing}
\usetikzlibrary{decorations.shapes}
\usetikzlibrary{decorations.text}
\usetikzlibrary{decorations.fractals}
\usetikzlibrary{decorations.footprints}
\usetikzlibrary{shadows}
\usetikzlibrary{calc}
\usetikzlibrary{spy}
\usepackage{amsmath}
\usepackage{array}
\usepackage{ amssymb }
\usepackage{braket}
\usepackage{qcircuit}
\usepackage{soul}
\usepackage{braket} 
\usepackage{relsize}

\usepackage{amsmath}
\usepackage{ amssymb }
\usepackage{braket}
\usepackage{qcircuit}
\usepackage{soul}
\usepackage{braket}
\usepackage{amsmath}
\usepackage{graphicx}%
\usepackage{multirow}%
\usepackage{amsmath,amssymb,amsfonts}%
\usepackage{amsthm}%
\usepackage{mathrsfs}%
\usepackage[title]{appendix}%
\usepackage{xcolor}%
\usepackage{textcomp}%
\usepackage{manyfoot}%
\usepackage{booktabs}%
\usepackage{algorithm}%
\usepackage{algorithmicx}%
\usepackage{algpseudocode}%
\usepackage{listings}%
\theoremstyle{thmstyleone}%
\theoremstyle{thmstyletwo}%

\theoremstyle{thmstylethree}%

\begin{document}

\title[Article Title]{Renormalization of crossing probabilities in the dilute Potts model}


\author[1]{\fnm{Pete} \sur{Rigas}}\email{pbr43@cornell.edu}

\affil[1]{\city{Newport Beach}, \postcode{92625}, \state{CA}, \country{United States}}


\abstract{A recent paper due to Duminil-Copin and Tassion from $2019$ introduces a novel argument for obtaining estimates on horizontal crossing probabilities of the Random-Cluster model, in which a range of four possible behaviors, through a quadrichotomy, is established. Novel renormalization arguments for crossing probabilities that the authors propose are studied in other models of interest that are not self-dual, specifically for the dilute Potts model. The probability measure of this model, through a suitably defined spin representation, is obtained from the high-temperature expansion of the loop $O(n)$ measure. The dilute Potts model was originally introduced in $1991$ by Nienhuis and is another model whose possible range of behaviors can be analyzed through a quadrichotomy; the range of four possible behaviors of the model can be respectively characterized at subcriticality, or supercriticality, in addition to either being critical discontinuous, or critical continuous. The exponential factor that is inserted into the Loop $O(n)$ model to quantify properties of the high-temperature phase  is proportional to the summation over all spins, and the number of monochromatically colored triangles over a finite volume, which is in exact correspondence with the parameter of a Boltzmann weight introduced in Nienhuis' paper detailing extensions of the $q$-state Potts model.}




\maketitle

\section{Introduction}

\subsection{Overview}

Russo-Seymour-Welsh (RSW) theory provides estimates regarding the crossing probabilities across rectangles of specified aspect ratios, and was studied by Russo, and then by Seymour and Welsch on the square lattice, with results specifying the finite mean size of percolation clusters {[23]}, in addition to a relationship that critical probabilities satisfy through a formalization of the sponge problem {[24]}. With such results, other models in statistical physics have been examined, particularly ones exhibiting sharp threshold phenomena {[1,7]} and continuous phase transitions {[13]}, with RSW type estimates obtained for Voronoi percolation {[27]}, critical site percolation on the square lattice {[28]}, the Kostlan ensemble {[2]}, and the FK-Ising model {[9]}, to name a few. To further develop RSW-type results for two dimensional systems, we demonstrate how a quadrichotomy of possible behaviors for the Random-Cluster model over $\textbf{Z}^2$ can be leverage for expressing that a similar quadrichotomy of possible behaviors holds for the spin representation of the loop $O(n)$ model over $\textbf{H}$. Irrespective of the fact that RSW-type results have been obtained for several models, it continues to remain of interest to obtain such results \textit{without} the use of duality.

RSW arguments which have classically relied upon model self-duality are captured through the following mathematical characteristics. The \textit{self-duality} property is said to be  enjoyed by a model iff the probability measure can be expressed with a dual measure. By self-dualiy in the case of $\textbf{Z}^2$ for the planar random-cluster model RSW results have provided a range of four possible behaviors. However the \textit{self-duality} property is enjoyed by neither the random cluster nor the dilute Potts models. Through an adaptation of the 2019 renormalization of crossing probabilities argument due to Duminil-Copin and Tassion {[14]}, crucial modifications for renormalizing crossing probabilities in the dilute Potts model arise not only from analogues of the Spatial Markov Property (SMP) and Comparison between Boundary Conditions (CBC) properties which intrinscially capture the model's dependence with the hexagonal lattice, but also in arguments for proving that the (PushPrimal), (PushDual), (PushPrimal Strip) and (PushDual Strip) conditions hold. The constants provided in the lower bounds of PushPrimal and PushDual conditions additionally impact arguments throughout Section \textit{7} and Section \textit{8} surrounding strip and renormalization inequalities, hence allowing for a classification of four phases of behavior of the dilute Potts model (the dichotomy of four possible behaviors is introduced in Section \textit{3.3}). Although classical RSW arguments have proven to be successful for analyzing self-duality, previous arguments have not been applicable to the dilute Potts model. The dilute Potts, as a model in Statistical Physics which is not self-dual can be defined in correspondence with the high-temperature loop $O(n)$ model in presence of two external fields.

Incorporating the two external fields, which are respectively denoted with $h,h^{\prime}>0$ is provided in the definition of the probability measure in Section \textit{3}. With respect to Nienhuis' critical point of the loop $O(n)$ model, which holds when the number of loops $n$ of a configuration is restricted between $0 \leq  n < 2$, we exhibit how the (SMP), and (CBC), properties holding for the planar random-cluster model also imply that RSW-type results hold, also with the Fortuin-Kestelyn-Ginibre (FKG) lattice conditions. Previous studies of the loop $O(n)$ model have been initiated by Nienhuis {[15,19,20]} who not only conjectured that the critical point of the model should be $1/\sqrt{2+\sqrt{2-n}}$ for $0 \leq n < 2$, but also by Duminil-Copin and Smirnov, {[22]}, which proved that the connective constant of the honeycomb lattice is $\sqrt{2+\sqrt{2}}$ {[12]}. It is also known that the loop $O(n)$ model, as a model for random collections of loop configurations on the hexagonal lattice, is conjectured to exhibit a phase transition with critical parameter $1/\sqrt{2+\sqrt{2-n}}$, in which, for \textit{subcritical} behavior, the probability of obtaining a macroscopic loop configuration of length $k$ decays exponentially fast in $k$, while \textit{at criticality}, the probability of obtaining infinitely many macroscopic loop configurations, taken to be of the same length $k$ and centered about the origin is bound below by $c$ and above by $1-c$, for $c \in (0,1)$, irrespective of boundary conditions {[8]}. While behavior below, and above, Nienhuis' critical point has not been formally established, it is widely believed to be true. The existence of macroscopic loops in the loop $O(n)$ model has also been proved in {[3]} with the XOR trick.

Given the above discussion of the loop model and its accompanying parameters, we return to the planar random-cluster model, and its quadrichotomy of possible behaviors that is central to the RSW-type result. In {[14]}, Duminil-Copin \textit{\&} Tassion proposed alternative arguments to obtain RSW estimates for models that are not self-dual at criticality. The novel quantities of interest in the argument involve renormalization inequalities, which in the case of Bernoulli percolation can be viewed as a coarse graining argument, as well as the introduction of strip densities which are quantities defined as a limit supremum over a real parameter $\alpha$. Ultimately, the paper proves RSW estimates for measures with free or wired boundary conditions in \textit{subcritical}, \textit{supercritical}, \textit{critical discontinuous} $\&$ \textit{critical continuous} cases, with applications of the two theorems relating to the mixing times of the random cluster measure, for systems undergoing discontinuous phase transitions {[14,18]}. Near the end of the introduction, the authors mention that potential generaliations of their novel renormalization argument can be realized in the dilute Potts model studied by Nienhuis which is equivalent to the loop $O(n)$ model, a model conjectured to exist in the same universality class as the spin $O(n)$ model. 

Obtaining RSW-type results for several models in Statistical Physics, whether the model depends upon self-duality or not, is continually of great interest to pursue. Relatedly, on the discrete level being able to argue that RSW-type results hold for a broad class of models would also provide information pertaining to universality. With regards to the loop $O(n)$ model, previous arguments have demonstrated that the connective constant, along with other information that is crucial to critical behavior, has been obtained with the \textit{parafermionic observable}. While such an observable was originally introduced to study conformal invariance of different models in several celebrated works {[11,25,26]}, determining the classes of boundary conditions for which the observable is sufficiently well behaved has historically been very difficult. 

As a holomorphic function, in certain cases representations for the discrete contour integral of the observable vanish under specific assumptions on a prefactor to the winding term in the power of the exponential. Straightforwardly, the wining number of a loop configuration over the hexagonal lattice for the loop $O(n)$ model is supported over subgraphs of even degree and can hence be determined from the number of self-intersection points. Under such assumptions on $\sigma$, when $x \equiv x_c$ at criticality Duminil-Copin $\&$ coauthors make use of several arguments with the Cauchy Riemann equations, which is related to strong holomorphicity which has been identified for several classes of related observables under a narrow class of boundary conditions {[8]}. For several observables which satisfy discrete holomorphicity, scaling relations, critical points, and several related objects of interest can be deduced. In the past disorder operators share connections with the parafermionic observable and have been studied to prove the existence of phase transitions through examination of the behavior of expectations of random variables exactly at the critical point {[11,16]}, while other novel uses of the parafermionic observable have been introduced in {[10]}. To this end, it is of interest to formulate RSW-type arguments for the six-vertex model which is in preparation in another paper for next year.

\subsection{Organization of results}

We define the models of interest to introduce Spin and Loop configurations, from which modifications to the (SMP) and (CBC) properties (defined in Section \textit{3.2}) yields bounds for crossings across symmetric domains (such domains are obtained through the high-temperature expansion for the loop $O(n)$ model, one configuration of which is provided in \textbf{Figure} \textit{1}). The parameters that we introduce to define such models, the length and number of loops, are related to the dimensions across which crossings occur in the RSW-type result. Moreover, while deterministic exploration algorithms which have been introduced for the plana random-cluster model also hold for the loop $O(n)$ model, stochastic exploration algorithms are also formulated. In comparison to deterministic exploration algorithms which define the top, bottom, left, and right boundaries of a symmetric domain, the stochastic counterpart determines each boundary of the domain depending upon crossing probabilities. For the spin representation of the loop $O(n)$ model, we describe how RSW behavior can be captured through deterministic and stochastic exploration.

In \textit{3} we define the Loop $O(n)$ measure, from which the dilute Potts measure in the presence of two external fields $h$ and $h^{\prime}$ can be defined, where $h, h^{\prime}$ are taken to be real parameters. This high-temperature expansion of the loop $O(n)$ model is used for arguments appearing in Section \textit{4} $\&$ Section \textit{5} with the proof of \textbf{Theorem 1} and \textbf{Lemma $9$}. In Section \textit{6}, we apply the $\mu$ homomorphism to lower bound vertical crossings with horizontal crossings, from which (PushPrimal) $\&$ (PushDual) conditions are introduced in Section \textit{7} to prove horizontal and vertical strip density formulas. In Section \textit{8} $\&$ Section \textit{9}, we characterize two behaviors of the quadrichotomy, finalizing our characterization of the  discontinuous-continuous phases of the quadrichotomy behavior by making use of the parafermionic observable which has already been manipulated to characterize properties of the phase transition for the random cluster model {[5,6,13]}. We conclude Section \textit{9} by obtaining classical results for the dilute Potts measure in each of the four regimes of behavior.

\section{Background}

To execute steps of the renormalization argument in the hexagonal case, we introduce quantities to avoid making use of self duality arguments. In all forthcoming arguments, the graphs that will be considered for the loop $O(n)$ model are the hexagonal lattice, $\textbf{H}$. First, we will describe how loops can be sampled uniformly at random from the sample space of the loop $O(n)$ model. Second, to obtain the high-temperature expansion, we incorporate labels on vertices of $\textbf{H}$, while assigning values of the spins for the faces of $\textbf{H}$. For $G=(V,E)$, $n \geq 1$ and the strip $\textbf{R} \times [-N,2N] \equiv S_N \subset G$, let $\phi_{S_N}^{\xi}$, for $\xi \in \{0,1, 0/ 1 \}$, corresponding to the boundary conditions,

{\small \[
\left\{\!\begin{array}{ll@{}>{{}}l}
  0 \equiv \textit{Boundary conditions for the Random-Cluster measure where edges of the boundary} \\ \textit{ of a finite volume of $\textbf{Z}^2$  are left open, ie unoccupied}  , \\ \\  1 \equiv \textit{Boundary conditions for the Random-Cluster measure where edges of the boundary } \\ \textit{of a finite  volume of $\textbf{Z}^2$ are left closed, ie occupied}   , \\ \\ 0 \backslash 1 \equiv \textit{Boundary conditions for the Random-Cluster measure where a partition of } \\ \textit{boundary of a finite volume of $\textbf{Z}^2$ are left closed, ie occupied}   , 
\end{array}\right.
\] }

\noindent respectively denote the measures with free, wired and Dobrushin boundary conditions in which all vertices at the bottom of the strip are wired. From such measures on the square lattice, several planar crossing events are defined in order to obtain RSW estimates for all four parameter regimes ({\textit{subcritical, supercritical, discontinuous $\&$ continuous critical}}), including analyses of the intersection of crossing probabilities across a family of non disjoint rectangles $\mathcal{R}$, each of aspect ratio $[0, \rho n] \times [0,n]$ for $\rho >0 $, to obtain crossings across long rectangles from the FKG inequality, three arm events which establish lower bounds of the crossing probabilities across $\mathcal{R}$ under translation and reflection invariance of $\phi$, in addition to horizontal rectangular crossings which are used to prove renormalization inequalities through use of (PushPrimal) $\&$ (PushDual) relations. To begin, we define the horizontal and vertical crossing strip densities. Denote,

\begin{align*}
  \phi^{\xi}_{\Lambda} \big[ \omega  \big] =  \frac{p^{o( \omega )} \big( 1 - p \big)^{c(\omega) } q^{k(\omega)}}{Z^{\xi}_{\Lambda} \big( \omega \big) }      ,
\end{align*}

\noindent corresponding to the probability measure of the random-cluster model under boundary conditions $\xi$, where,

\begin{align*}
  p^{o( \omega )} \equiv \textit{The probability $p$ raised to the number of open edges of $\omega$}  , \\ \\ 1 -   p^{c( \omega )} \equiv \textit{The probability $1-p$ raised to the number of closed edges of $\omega$}  , \\ \\  q^{k(\omega)} \equiv \textit{the number of clusters contained within $\omega$} , 
\end{align*}

\noindent and normalizing constant $Z^{\chi}_{\Lambda} \big( \omega \big)$, the partition function, so that the probabilities distributed under the law $\phi^{\xi}_{\Lambda} \big[ \cdot \big]$ sum to $1$.

\begin{figure}[H]
\begin{center}
\includegraphics[width=0.85\columnwidth]{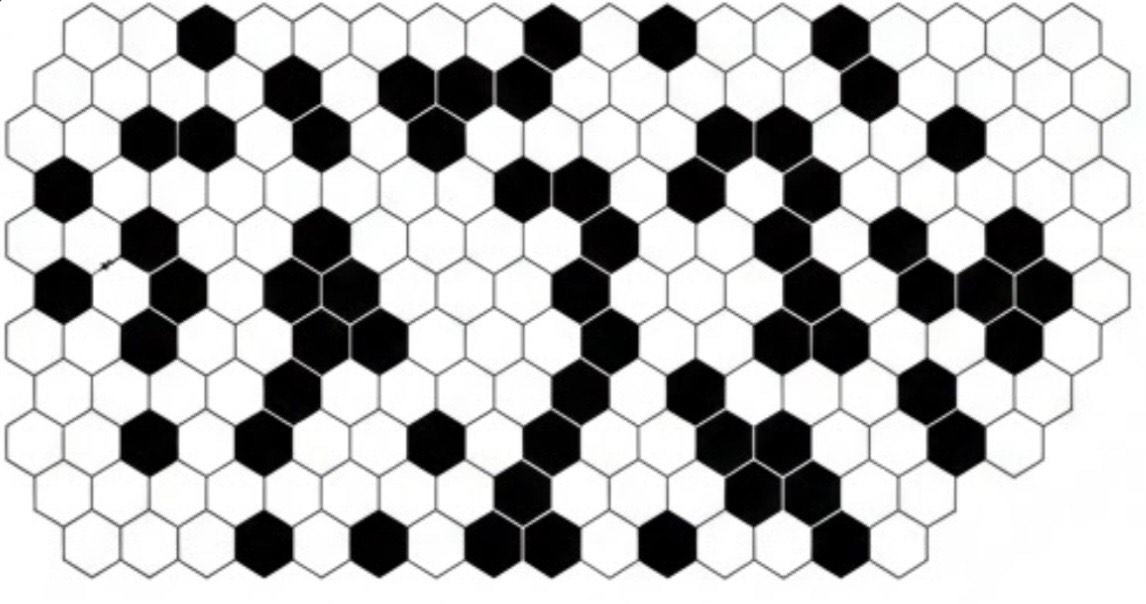}
\end{center}
\caption{A depiction of percolation over faces of the hexagonal lattice. With probability $p$, a face of the hexagonal lattice is colored black and otherwise white with probability $1-p$. To quantify connectivity properties of the underlying lattice, one is interested in constructing paths of neighboring faces which are all colored black.}
\end{figure}

\noindent \textbf{Definition 1} (\textit{the vertical strip density for the planar random-cluster model}, {[14]}): The \textit{strip density} corresponding to the measure across a rectangle $\mathcal{R}$ of aspect ratio $[0,\alpha N ] \times [- N , 2 N]$ with free boundary conditions is of the form, 
\begin{align*}
   p_N = \mathrm{lim \text{ } sup}_{\alpha \rightarrow \infty} \big\{ \phi^0_{[0,\alpha N]\times [-N,2N]}\big[ \mathcal{H}_{[0,\alpha N] \times [0,N]}\big] \big\}^{\frac{1}{\alpha}}   \text{ , }
\end{align*}

\noindent where $\mathcal{H}$ denotes the event that $\mathcal{R}$ is crossed horizontally, whereas for the measure supported over $\mathcal{R}$ with wired boundary conditions, the crossing density is of the form,

\begin{align*}
    q_N = \mathrm{lim \text{ } sup}_{\alpha \rightarrow \infty} \big\{ \phi^1_{[0,\alpha N]\times [-N,2N]}\text{ } \big[  \mathcal{V}_{[0,\alpha N] \times [0,N]}^c \big]  \big\}^{\frac{1}{\alpha}}\text{ , }
\end{align*}

\noindent where $\mathcal{V}^c$ denotes the complement of a vertical crossing across $\mathcal{R}$.

\bigskip

The strip and renormalization inequalities provided in this section are dependent on different quantities for $\pm$ spin configurations rather than the corresponding inequalities for the random cluster model which only depend on the cluster weight $q$. Besides the definition of the strip densities $p_n$ and $q_n$, another key step in the argument involves inequalities relating $p_n$ and $q_n$. The behavior of such inequalities is related too the manner in which connected components coalesce over the hexagonal lattice, a few examples of which are depicted in \textbf{Figure} \textit{2}, and in \textbf{Figure} \textit{3}. The statement of the Lemma below holds under the assumption that the planar random cluster model is neither in the subcritical nor supercritical phase.

\begin{figure}[H]
\begin{center}
\includegraphics[width=0.92\columnwidth]{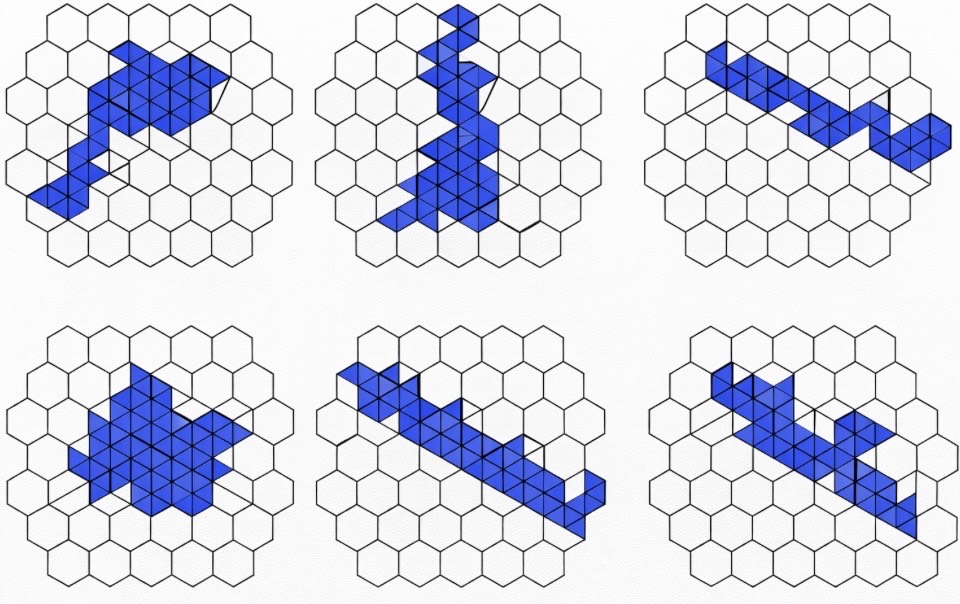}
\end{center}
\caption{A depiction of six percolation configurations corresponding to percolation over faces of the triangular lattice. Macroscopic crossings are obtained either from one side of the domain to another, or between several other possible collections of faces. Faces across which percolation occurs are colored in blue.}
\end{figure}

\noindent\textbf{Lemma 1} (\textbf{Lemma 12}, {[14]}) There exists a constant $C > 0$ such that for every integer $\lambda \geq 2$, and for every $n \in 3\textbf{N}$,

\begin{align*}
   p_{3n} \geq \frac{1}{\lambda^C} q_n^{3 + \frac{3}{\lambda}} \text{ , }
\end{align*}

\noindent while a similar inequality holds between horizontal and the complement of vertical crossing probabilities of the complement $\mathcal{V}^c$ across $\mathcal{R}$, which takes the form,

\begin{align*}
    q_{3n} \geq \frac{1}{\lambda^C} p_n^{3 + \frac{3}{\lambda}} \text{ . }
\end{align*}

Finally, we introduce the renormalization inequalities.

\bigskip

\noindent\textbf{Lemma 2} (\textbf{Lemma 15}, {[14]}) There exists $C > 0$ such that for every integer $\lambda \geq 2$ and for every $n \in 3 \textbf{N}$, 

\begin{align*}
      p_{3n} \leq \lambda^C p_n^{3 - \frac{9}{\lambda}} \text{         } \text{        \&        } \text{     } q_{3n} \leq \lambda^C q_n^{3 - \frac{9}{\lambda}} \text{ . }
\end{align*}

To readily generalize the renormalization argument to the dilute Potts model, we proceed in the spirit of {[14]} by introducing hexagonal analogues of the crossing events discussed at the beginning of the section.

\section{Towards hexagonal analogues of crossing events from the random cluster renormalization argument}

\subsection{Loop $O(n)$ measure, hexagonal lattice crossing events}

\noindent Denote $\Lambda$ as some finite volume over the hexagonal lattice $\textbf{H}$, and boundary conditions $\xi$. An instance of a boundary condition can be obtained by coloring, with probability $p$ and $1-p$ for some $p>0$, each face of $\textbf{H}$ as white, and black, respectively. Before passing to the high-temperature expansion, the Gibbs measure on a random configuration $\sigma$ in the loop $O(n)$ model is of the form,

\begin{align*}
    {\textbf{P}}_{\Lambda,x,n}^{\xi}(\sigma) = \frac{x^{e(\sigma)} n^{l(\sigma)}}{Z_{\Lambda, x,n}^{\xi}} \tag{Loop measure} \text{ , }
\end{align*}

\noindent where $e (\sigma)$ denotes the number of edges, $l( \sigma ) $ the number of loops, given by,

\begin{align*}
   l \big( \sigma \big) =    \# \big\{ \textit{loops: loops belong to $\sigma$} \big\}              , 
\end{align*}

\noindent $\Lambda \subset \textbf{H}$, $\xi \in \{0,1, 0 / 1 \}$ corresponding to the collection of possible boundary conditions,

{\small \[
\left\{\!\begin{array}{ll@{}>{{}}l} 
 0 \equiv \textit{Boundary conditions for the high-temperature spin measure where the face of a } \\ \textit{hexagon incident  to the boundary is not colored black}   , \\ \\  1 \equiv \textit{Boundary conditions for the high-temperature spin measure where the face of a } \\ \textit{hexagon incident to the boundary is  colored black}   \\ \\ 0 \backslash 1 \equiv  \textit{Boundary conditions for the high-temperature spin measure where a fixed } \\ \textit{proportion of faces of hexagons are not colored black, while the remaining proportion } \\ \textit{of faces of hexagons are colored black}           , 
\end{array}\right.
\] }

\noindent and $Z_{\Lambda, x,n}^{\xi}$ is the partition function which normalizes ${\textbf{P}}_{\Lambda, x,n}^{\xi} \big[ \cdot \big] $ so that it is a probability measure. Equivalently, the loop configurations distributed according to the above law can also be thought of as \textit{loop} configurations, namely subgraphs of even degree. To this end, we denote $G$ as a finite subgraph of $\textbf{H}$, and would like to study the behavior of the Spin Measure below as $G \longrightarrow \textbf{H}$. The notion by which $G$ will be said to converge to that of $\textbf{H}$ is obtained through the weak limit $\longrightarrow$. In particular, we restrict the parameter regime of $x$ to that of {[8]}, in which the loop $O(n)$ model satisfies the strong FKG lattice condition and monotonicity through a spin representation measure albeit ${\textbf{P}}_{\Lambda,x,n}^{\xi}$ not being monotonic. By construction at criticality when $x = x_c \big( n \big)$ ${\textbf{P}}_{\Lambda,x,n}^{\xi}$ is invariant under $\frac{2 \pi}{3}$ rotations. For the high-temperature expansion below, we take the number of loops $n$ to satisfy $n \geq 2$. This serves the role of considering the supercritical behavior of the loop $O(n)$ model. Through this extension for $n \geq 2$ of the spin representation of $\textbf{P}_{\Lambda,e(\sigma),l(\sigma)}^{\xi} \big[ \cdot \big]$, the measure on spin configurations $\sigma^{\prime} \in \Sigma(G, \tau)$, for configurations

\begin{align*}
    \sigma^{\prime} =  \big\{  \textit{Spin configuration of the high-temperature expansion obtained from $\sigma \sim \textbf{P}^{\xi}_{\Lambda,x,n} \big[ \cdot \big]$}  \big\}  \\ \in \big\{ \pm 1 \big\}^{F ( \textbf{H}^* ) }  , 
\end{align*}

\noindent with,

{\small \begin{align*}
    \big\{ \pm 1 \big\}^{F ( \textbf{H}^* ) } \equiv \big\{ F \big( \textbf{H}^{*} \big) :  \big\{  F \big( \textbf{H}^{*} \big) =  +1 \Longleftrightarrow \textit{three vertices of the face of $\textbf{H}^{*}$ are colored black, or white} \big\}  , \\ \big\{ F \big( \textbf{H}^{*} \big) =  - 1  \Longleftrightarrow \textit{three vertices of the face of $\textbf{H}^{*}$ are not all colored black, or white} \big\}  \big\} , 
\end{align*}     }

\noindent whose length, and number of connected components, are defined as,

\begin{align*}
 l \big( \sigma^{\prime} \big) =  \textit{length of $\sigma^{\prime}$}  , \\ \\ k \big( \sigma^{\prime} \big) =   \textit{number of connected components of $\sigma^{\prime}$}  =      \# \big\{ \textit{faces} \in F \big( \textbf{H}^{*} \big) :  \textit{faces} = \big\{ + 1 \big\}^{F(\textbf{H}^{*} )}   \big\}                   , 
\end{align*}

\noindent supported over

\begin{align*}
  G =   \big( V \big( \textbf{H} \big) , E \big( \textbf{H} \big) \big)      , 
\end{align*}

\noindent is of the form,

\begin{align*}
     \mu_{G, x , n}^{\tau}(\sigma^{\prime}) = \frac{ {\LARGE n^{k(\sigma^{\prime})} x^{e(\sigma^{\prime})} \text{ } \mathrm{exp} \big[   h r(\sigma^{\prime}) + \frac{h^{\prime}}{2} r^{\prime}(\sigma^{\prime})  \big]  } }{Z_{G,x,n}^{\tau}} \tag{Spin Measure} \text{ , }
\end{align*}

\noindent where $\tau \in \{-1,+1 \}^{\textbf{T}}$, $\Sigma(G,\tau)$ is the set of spin configurations coinciding with $\tau$ outside of $G$, $r(\sigma^{\prime}) = \sum_{u \in G} \sigma^{\prime}_u$ is the summation of spins inside $G$, $r^{\prime}(\sigma^{\prime}) = \sum_{\{u,v,w\} \in G} \sigma^{\prime}_u \textbf{1}_{\{ \sigma^{\prime}_u = \sigma^{\prime}_v = \sigma^{\prime}_w \} }$ is the difference between the spins of monochromatic triangles, and $Z_{G,n,x}^{\tau}$ is the partition function which makes $\mu_{G , x, n }^{\tau} \big[ \cdot \big]$ a probability measure. The extension enjoys translation invariance, the Spatial Markov property that will be mentioned in \textit{5.1}, comparison between boundary conditions that is mentioned in \textit{3.2}, $\&$ (FKG) for $n \geq 1$ and $nx^2 \leq 1$. The dual measure of $\mu_{G,x,n}^{1}$ is $\mu_{G^{*},x,n}^{0}$. Simply put, the superscripts above $\mu$ indicate whether the pushforward of a horizontal or vertical crossing event under the measure is under free, wired, or mixed boundary conditions. 

Additionally, the model can be placed into correspondence with the dilute Potts model, originally characterized by occupied, and vacant, faces of $\textbf{H}$. The exponential factor introduced to characterize high temperature behavior of the Loop $O(n)$ model in the Spin Measure equality is in direct correspondence with the dilute Potts model, with Boltzmann weight,

\begin{align*}
  \mathscr{W}_{ij} \text{ } \equiv \text{ }      \prod_{i \sim j} \big(    1 - t_i t_j + t_i t_j \delta_{s_j,s_k}     \big)   \text{ } \mathrm{exp} \bigg[      K_1 \sum_i t_i + K_2 \sum_{i \sim j}   t_i t_j + K_3 \sum_{i \sim j \sim k} t_i t_j t_k     \bigg]   \text{, } 
\end{align*}

\noindent where the quantities in the nearest-neighbor product above include the occupation number $t_j$, which is either equal to $0$ or $1$, and the spins $s_j$ and $s_k$ indicated in the nearest neighbor product of $i$ and $j$ can take values between $1,\cdots , q$ corresponding to $q$-state Potts model spins [22]. In the power of the exponential term, we consider occupancy numbers across sites, edges, and faces of triangles, respectively, with $K_1,K_2,K_3$ real constants.

\begin{figure}[H]
\begin{center}
\includegraphics[width=1\columnwidth]{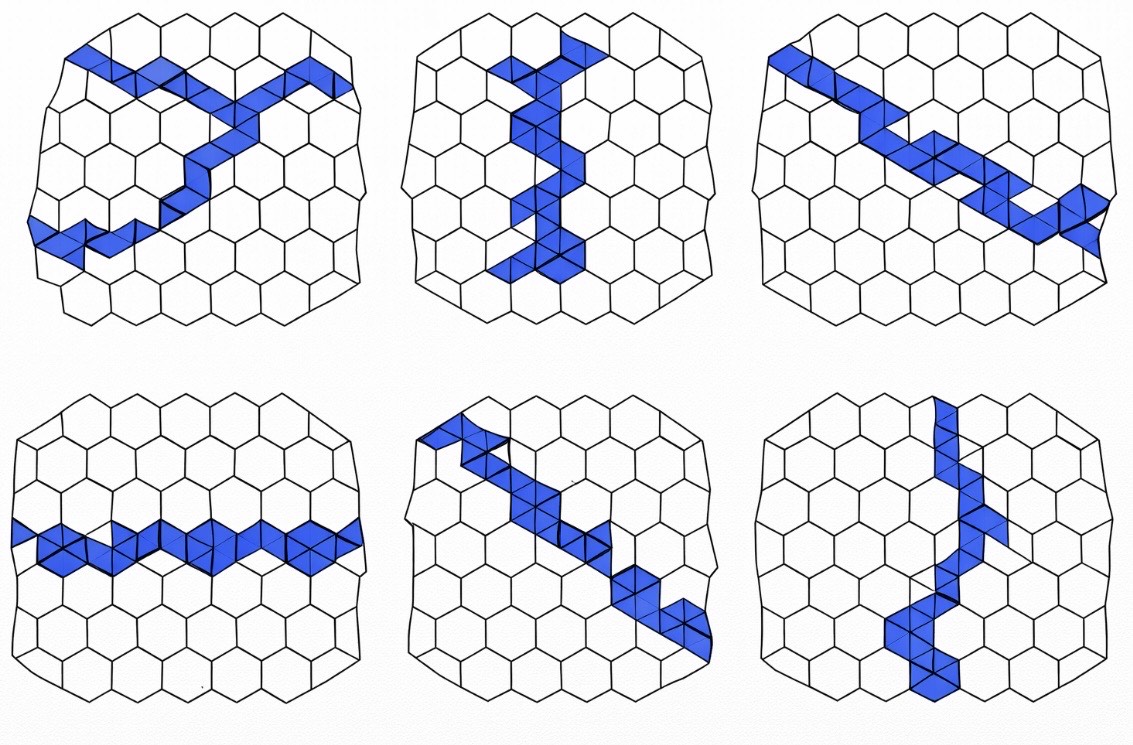}
\end{center}
\caption{Crossing probabilities across the triangular lattice with blue connected components can be constructed in several possible configuations within hexagons, as shown above.}
\end{figure}

\noindent To obtain boundary dependent RSW results on $\textbf{H}$ in all $4$ cases, we identify crossing events in the planar renormalization argument in addition to difficulties associated with applying the planar argument to the push forward of similarly defined horizontal and vertical crossing events under $\mu_{G,x,n}^{\tau}$ on $(\textbf{T})^{*} = \textbf{H}$. In what follows, we describe all planar crossing events in the argument.

\bigskip

\noindent First, planar crossing events across translates of horizontal crossings across short rectangles of equal aspect ratio are combined to obtain horizontal crossings across long rectangles, through the introduction of a lower bound to the probability of the intersection that all short rectangles are simultaneously crossed horizontally with (FKG). On $\textbf{H}$, the probability of the intersection of horizontal crossing events can be readily generalized to produce longer horizontal crossings from the intersection of shorter ones, through an adaptation of  \textbf{Lemma 9} in {[14]}. The forthcoming arguments for the loop $O(n)$ model are related to conditioning on the existence of distinct connected components, a few configurations of which are depicted in \textbf{Figure} \textit{4}, \textbf{Figure} \textit{5} and \textbf{Figure} \textit{7}. Possible choices of symmetric domains that can be constructed from distinct connected components are depicted in \textbf{Figure} \textit{6}.

\noindent Explicitly, we define crossing events for the Spin Measure supported over $\textbf{H}$ as collections of spins satisfying the following conditions,

{\small \begin{align*}
  \textit{Crossings over $\textbf{H}$ for spin configurations $\sim \mu^{\tau}_{G,x,n} \big[ \cdot \big]$} \equiv \big\{\forall  \textit{path} \subsetneq \textbf{H} , \exists \sigma \neq  \sigma^{\prime} \in \big\{ \pm 1\big\}^{F(\textbf{H}^{*})}\\  : \textit{path connects $\sigma$ and $\sigma^{\prime}$ with only $+1$ faces} \big\}  . 
\end{align*} }

\begin{figure}[H]
\begin{center}
\includegraphics[width=0.92\columnwidth]{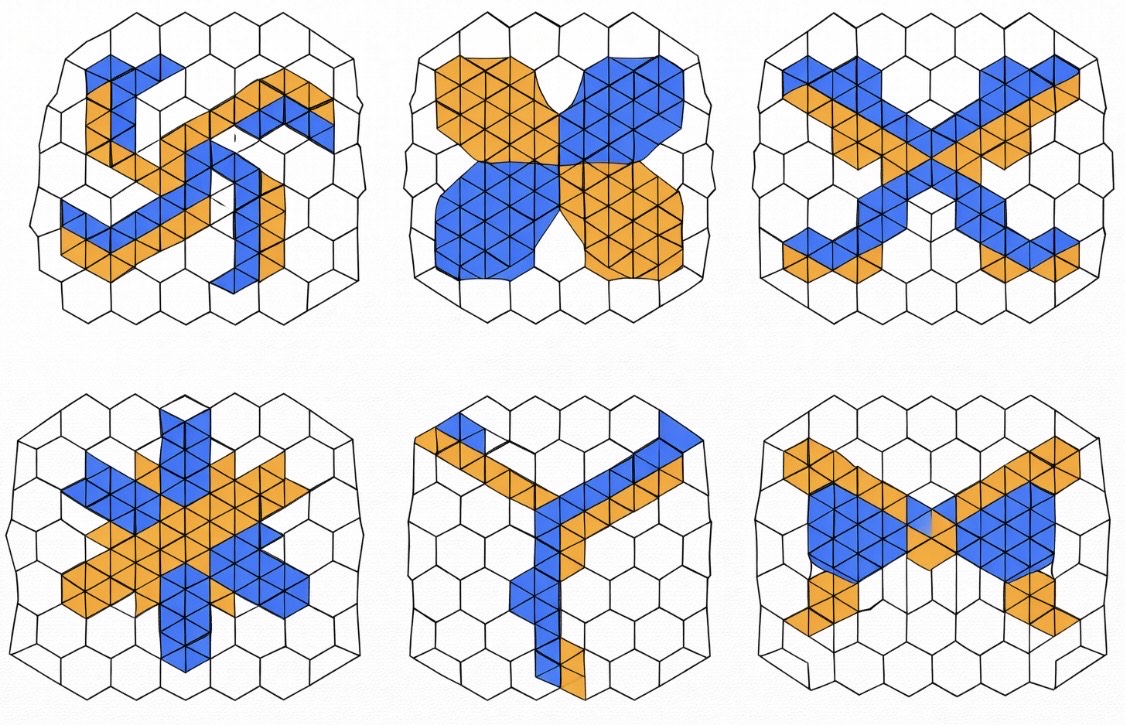}
\end{center}
\caption{Disjoint connected components of the high-temperature expansion of the Loop $O(n)$ model, as depicted above, can  be used to construct symmetric domains.}
\end{figure}

Second, three arm events which determine whether two horizontal crossings to the top of a rectangle of aspect ratio $[0,N] \times [0,\rho N]$ intersect. Planar crossings used to create symmetric domains over which the conditional probability of horizontal crossings in the symmetric domain can be determined, which for the renormalization argument rely on comparison between random cluster measures with free and wired boundary conditions. For random cluster configurations, comparison between boundary conditions is established in how the number of clusters in a configuration is counted. Comparison between boundary conditions applies to $\mu_{G,x,n}^{\tau} \big[ \cdot \big]$ from {[8]}, with hexagonal symmetric domains enjoying $\frac{2 \pi}{3}$ rotational symmetry.

Third, crossing events with wired boundary conditions, induce wired boundary conditions within close proximity of vertical crossings in planar strips. Long horizontal crossings are guaranteed through applications of (FKG) across dyadic translates of horizontal crossings across shorter rectangles. For hexagonal domains, modifications to planar crossings of first type permit ready generalizations of third type planar crossings. 

Fourth, planar horizontal crossing events across rectangles establish relations between the strip densities $p_N$ $\&$ $q_N$ (provided in \textbf{Lemma 1}). Finally, planar crossing events satisfying (PushPrimal) $\&$ (PushDual) properties prove \textbf{Lemma 2}.

\begin{figure}[H]
\begin{center}
\includegraphics[width=1\columnwidth]{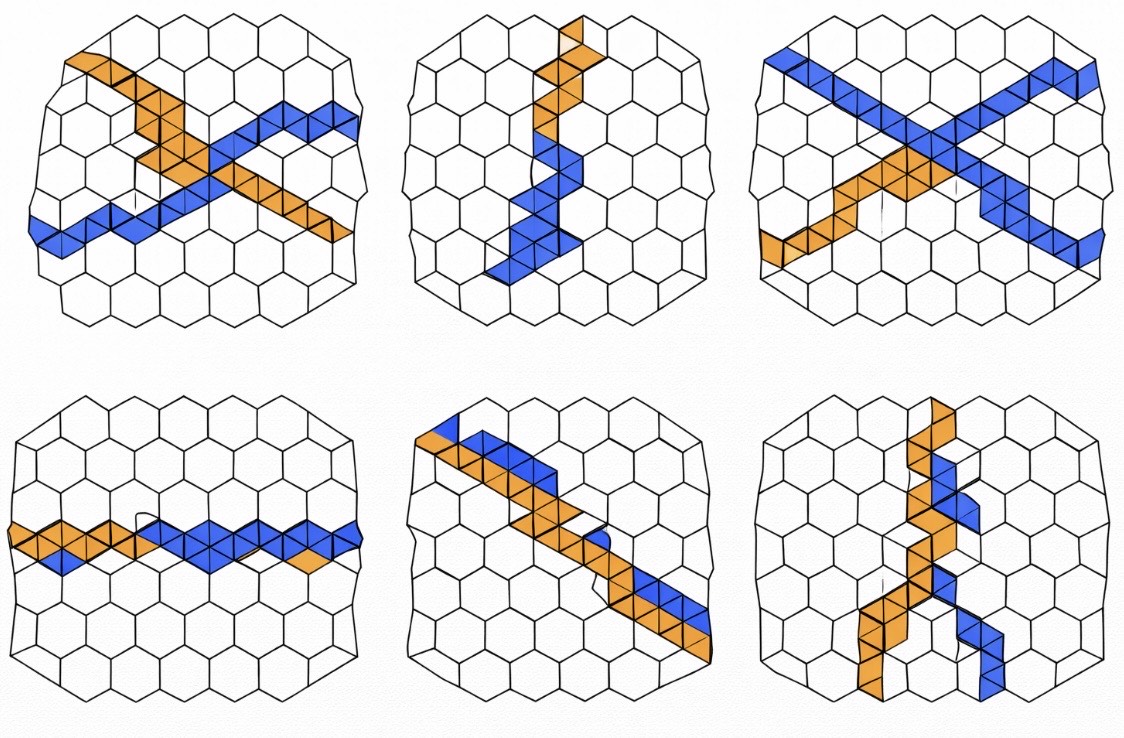}
\end{center}
\caption{Besides connected components displayed in Figure 4, arguments for demonstrating that crossings across symmetric domains occur with high probability are first shown to hold within the interior of a box over $\textbf{H}$.}
\end{figure}

\subsection{Comparison of boundary conditions $\&$ Spatial Markovianity for the $n \geq 2$ extension of the loop $O(n)$ measure}

For suitable comparison of boundary conditions in the presence of external fields $h,h^{\prime}$, the influence of boundary conditions from the fields on the spin representation amounts to enumerating configurations differently for wired and free boundary conditions than for the random cluster model in {[14]}. In particular, one makes use of straightforward adaptations of two properties corresponding to (CBC) and (SMP) for the loop $O(n)$ model which are defined below.

For the loop model, crossings across admissible symmetric domains, $\mathrm{Sym}$, inherit boundary conditions from partitions on the outermost layer of hexagons along loop configurations (see \textbf{Figure} \textit{8} corresponding to crossings over a finite volume of $\textbf{H}$ in Section \textit{4.2}). Through distinct partitions of the $+/-$ assignment on hexagons on the outermost layer to the boundary, appearing in arguments for symmetric domains appearing in Sections \textit{5.1 - 5.4}.

\bigskip

\noindent \textbf{Corollary}  (\textbf{Corollary 10}, {[8]}), \textit{comparison between boundary conditions for the Spin measure}): Consider $G \subsetneq  \textbf{T}$ finite and fix $(n,x,h,h^{\prime})$ such that $n \geq 1$ and $nx^2 \leq \mathrm{exp}(-|h^{\prime}|)$. For any increasing event $A$ and any $\tau \leq \tau^{\prime}$,

\begin{align*}
       \mu^{\tau}_{G,n,x,h,h^{\prime}}[A] \leq \mu^{\tau^{\prime}}_{G,n,x,h,h^{\prime}}[A] \text{, } \tag{$\mathcal{S}- \mathrm{CBC}$} 
\end{align*}

\noindent Altogether, modifications to comparison of boundary conditions and the spatial Markov property between measurable spin configurations for $\mu$ is also achieved. We recall the (CBC) inequality for the random cluster model, and for the loop model make use of an 'analogy' discussed in {[8]}, in which we associate wired boundary conditions to the $+$ spin, and free boundary conditions to the $-$ spin over the set of faces $F \big( \textbf{T} \big) = F \big( \textbf{H}^{*} \big)$. Specifically, for boundary conditions $\xi , \psi$ distributed under the random cluster measure $\phi$, the measure supported over $G$ satisfies

\begin{align*}
     \phi^{\xi}_{G}[\mathcal{A}]    \text{ }    \leq    q^{\mathrm{max}\{   k_{\xi}(\omega) - k_{\psi}(\omega) : \omega  \} - \mathrm{min}\{    k_{\xi}(\omega) - k_{\psi}(\omega)      \}} \text{ } \phi^{\psi}_G[\mathcal{A}]  \text{ . }
\end{align*}

\noindent For an arbitrary increasing event $\mathcal{A}^{\prime}$, a special case of $(\mathcal{S}-\mathrm{CBC})$ property above takes the form,

{\small \begin{align*}
          \mu^{\tau}_{\textbf{H}} [\mathcal{A}^{\prime}] \leq  n^{k_{\tau^{\prime}}(\sigma)  -  k_{\tau}(\sigma)}  \text{ } x^{e_{\tau^{\prime}}(\sigma) - e_{\tau}(\sigma)}  \text{ }  \mathrm{exp} \bigg[  h \big( r_{\tau^{\prime}}(\sigma) - r_{\tau}(\sigma)     \big)            + \frac{h^{\prime}}{2} \big(    r_{\tau^{\prime}}(\sigma^{\prime}) - r_{\tau}(\sigma^{\prime}) \big) \bigg] \text{ }  \mu^{\tau^{\prime}}_{\textbf{H}}[\mathcal{A}^{\prime}]          \tag{$\mathcal{S}- \mathrm{CBC}$}     \text{ , } 
\end{align*} }

\noindent Another property that the dilute Potts measure satisfies, for finite volumes $\mathcal{I} \subsetneq \mathcal{O}$,

\begin{align*}
    \mu^{\tau}_{\mathcal{I}} \big[   \cdot |_{\mathcal{I}} | \sigma^{\prime}_{\mathcal{I}} \equiv \sigma^{\prime}_{\mathcal{O}} \text{ for } \mathcal{I} \cap \mathcal{O}     \big] \equiv     \mu^{\tau\cup \tau^{\prime}}_{\mathcal{O}} \big[ \cdot \big]       \text{. }  \tag{$\mathcal{S}-\mathrm{SMP}$}
\end{align*}

\noindent where the exponential factor in front of the pushforward in the upper bound results from the difference between the number of monochromatically colored triangles in the configuration distributed under the Spin Measure, the edge weight associated with $x$, and the number of connected components $k(\sigma) + 1$, respectively with boundary conditions $\tau$ and $\tau^{\prime}$,

{\small \begin{align*}
  \tau \equiv   \textit{Boundary conditions for the Spin Measure associated with the finite vol-} \\ \textit{ume $\partial \mathcal{I}$}     , \\ \\ \tau^{\prime} \equiv   \textit{Boundary conditions for the Spin Measure associated with the finite volume} \\ \textit{ $\partial \big[ \mathcal{I} \cap \mathcal{O} \big]$}      . 
\end{align*} }

\noindent The multiplicative factor arises from comparisons between the Spin Measure and the Random-Cluster probability measure, particularly by associating the summation over all spins in Spin Measure configurations with the ratio of the number of open edges to the number of closed edges in an FK percolation configuration, the number of connected components in a spin configuration under the loop $O(n)$ model with the number of clusters in the Random cluster model, and also, the edge weights of $x$ of spin configurations under the loop $O(n)$ model with the cluster weights of $q$ in the Random-Cluster Model.

\begin{figure}[H]
\begin{center}
\includegraphics[width=1\columnwidth]{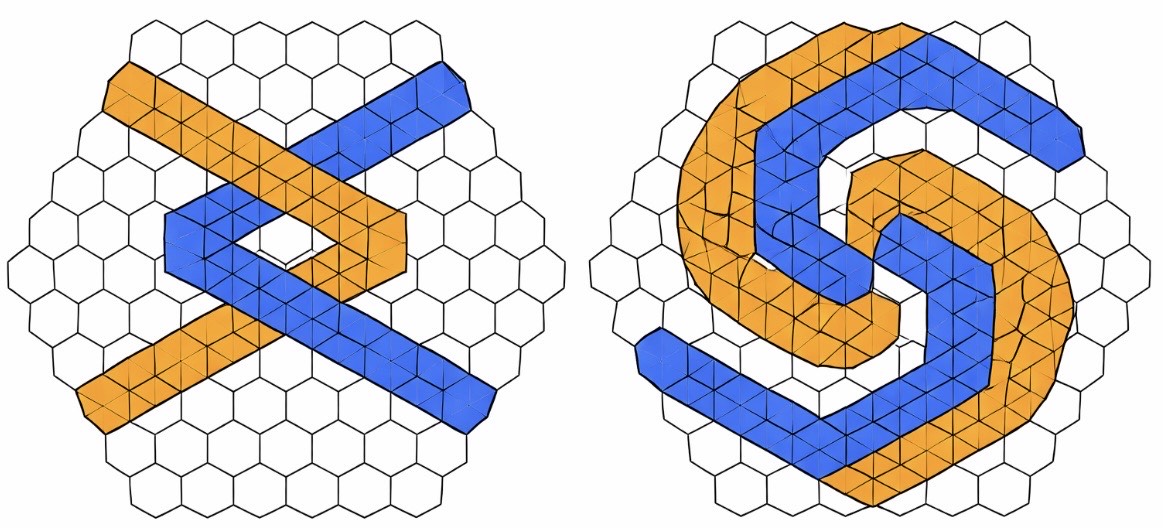}
\end{center}
\caption{Symmetric domains from connected components of faces over $\textbf{T}$ can be obtained by either taking the domain enclosed within the blue and orange paths, as displayed on the left; or by taking the region complementary to that spanned by the blue and orange connected components to the box over $\textbf{H}$, as displayed on the right.}
\end{figure}

A special case of the inequality above will be implemented several times throughout the renormalization argument, stating,

\begin{align*}
             \mu^{\tau}_{\textbf{H}} [\mathcal{A}^{\prime}] \leq  \text{ }  n^{k_{\tau^{\prime}}(\sigma) - k_{\tau}(\sigma)}    \text{ }  x \text{ }  \mathrm{exp}(h) \text{ }            \mu^{\tau^{\prime}}_{\textbf{H}}[\mathcal{A}^{\prime}]      \text{, }       
\end{align*}

\noindent which will be introduced when applying $(\mathcal{S} -\mathrm{SMP})$ in forthcoming arguments. The special case of the multiplicative factor above represents the difference in the number of clusters that are counted under boundary conditions $\tau, \tau^{\prime}$, in addition to the corresponding edge weights $x$ under each boundary condition. We denote the modified properties for spin representations $\mathcal{S}$ with $(\mathcal{S} -\mathrm{CBC})$ and $(\mathcal{S}-\mathrm{SMP})$.

\begin{figure}[H]
\begin{center}
\includegraphics[width=0.55\columnwidth]{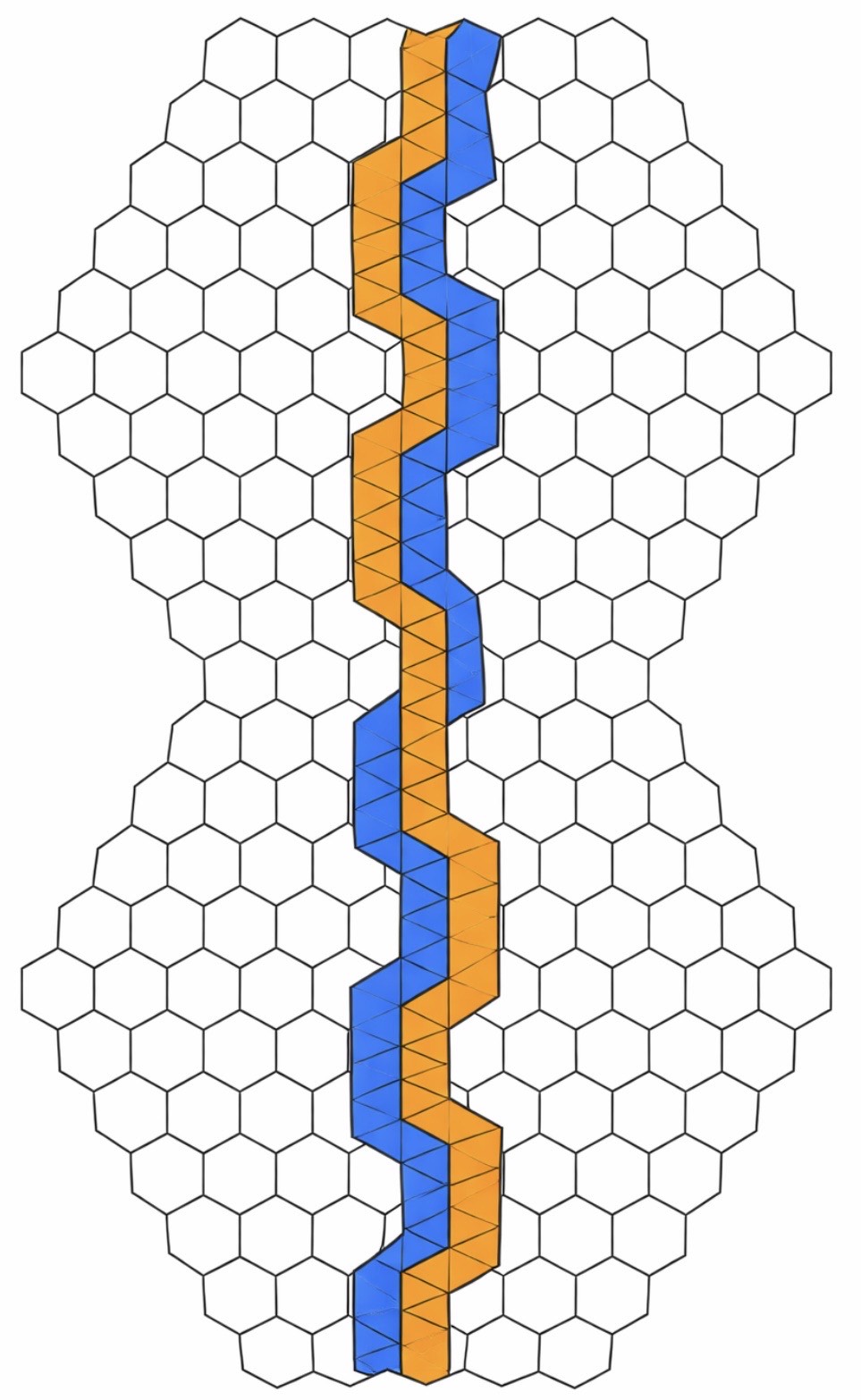}
\end{center}
\caption{Two connected components, as displayed above in blue and orange, traverse two boxes over $\textbf{H}$.}
\end{figure}

Besides such modifications, (MON) from {[14]} directly applies, and will be used repeatedly. To probabilistically capture the dependence of the high-temperature measure on the first external field, another special case of $(\mathcal{S}-\mathrm{SMP})$ is, 

{\small \begin{align*}
  \mathrm{exp} \bigg[ {   \#  \big\{ \textit{hexagonal boxes }  \mathcal{H}  :  1 \leq   i \sim j   \leq 6  ,  \sigma_v \in \{ \pm 1 \}^{V (\textbf{H} )}  , \textbf{1}_{\{  \sigma_{v_i} = \sigma_{v_j}   = 1\} }   \big\}   } \bigg]     \text{, } 
\end{align*} } 

where,

\begin{align*}
  \textit{hexagonal boxes} \equiv       \underset{F \in F ( \textbf{H} )}{\bigcup} \big\{\textit{F: F is spanned by the six vertices $\sigma_1, \cdots, \sigma_6 \in V \big( \textbf{H} \big)$}  \big\}         . 
\end{align*}

\noindent Finally, given two increasing events $A$ and $B$ (two events which have a higher probability of occurring if the number of edges in each respective event is increased), \textit{positive association} is an inequality of the form,

\begin{align*}
 \mu^{\tau} [  A \cap B  ]  \equiv  \mu^{\tau}_{\textbf{H}}[  A \cap B  ] \geq \mu^{\tau}_{\textbf{H}}[ A ] \mu^{\tau}_{\textbf{H}}[ B ] \equiv \mu^{\tau}[ A ] \mu^{\tau}[ B ]   \text{. }  \tag{$\mathrm{FKG}$}    
\end{align*}

\noindent We refer to free boundary conditions under the Spin measure which represents $-$ boundary conditions, and wired boundary conditions for $+$ boundary conditions. The $q^k$ analogue from the random cluster model for the high-temperature Spin measure will enter into the novel renormalization argument at several points, stating:

\begin{itemize}
     \item[$\bullet$]  \textbf{Lemma $9$}, in which $(\mathcal{S}-\mathrm{SMP})$ will be repeatedly used to compare boundary conditions between crossings across the second or third edge of a hexagon, and boundary conditions for crossings across symmetric regions $\mathrm{Sym}$,
     
     \item[$\bullet$]  \textbf{Corollary $11$}, in which $(\mathcal{S}-\mathrm{SMP})$ will be used to bound the pushfoward of horizontal crossings under wired boundary conditions, which in light of the homeomorphism $f$ in \textit{4.1}, yields a corresponding bound for the pushforward of a vertical crossing under free boundary conditions, 
     
     \item[$\bullet$]  \textbf{Lemma $1$}, in which an application of $(\mathcal{S}-\mathrm{SMP})$ and (MON) yield a lower bound for the probability of a horizontal crossing under free boundary conditions with a probability of a horizontal crossing under wired boundary conditions, 
     
     \item[$\bullet$] \textbf{Lemma $2$}, in which a modification to the lower bound obtained in the proof of \textbf{Lemma} $1$ is applied to obtain a lower bound for the probability of a horizontal crossing under wired boundary conditions with the probability of a vertical crossing under free boundary conditions, 
     
     \item[$\bullet$] \textit{Quadrichotomy proof}, in which the crossing events from previous results are compared to obtain the standard \textit{box crossing estimate} that the Gibbs measure on loop configurations satisfies, per conditions of the Continuous Critical phase of the dilute Potts model provided in \textbf{Theorem 1} {[8]}.

\end{itemize}

\subsection{Description of results}

The result presented for the loop $O(n)$ model mirrors the dichotomy of possible behaviors, in which the \textit{standard box crossing estimate} reflects the influence of boundary conditions on the nature of the phase transition, namely that the transition is discontinuous, from the \textit{{discontinuous critical}} case. To prove the \textit{subcritical} $\&$ \textit{supercritical} cases, the generalization to the dilute Potts model will make use of planar crossing events of first and second type, while crossing events of third and fourth type proves the remaining \textit{{continuous $\&$ discontinuous critical}} cases. We denote the vertical strip domain $\mathcal{S}_T$ with $T$ hexagons, $\mathcal{S}_{T,L}$ the finite domain of $\mathcal{S}_T$ of length $L > 0$, and any regular hexagon $H_j \subsetneq \mathcal{S}_T$ with side $j$, as extensively manipulated in {[12]}. The strip densities $p^{\mu}_n$ and $q^{\mu}_n$ are defined in Section \textit{7}.

\bigskip

\noindent\textbf{Theorem} $1$ (\textit{Spin Mesaure $\mu^{\tau}_{G,x,n} \big[ \cdot \big]$ homeomorphism}): For $L \in [0,1]$, there exits an increasing homeomorphism $f_{L}$ so that for every $n \geq 1$, where $\mathcal{H}_H \equiv \mathcal{H}$ and $\mathcal{V}_H \equiv \mathcal{V}$ denote the horizontal and vertical crossings across a regular hexagon $H$, $\mu^{\tau}_{G,x,n}(\mathcal{H}) \geq f(\mu^{\tau}_{G,x,n}(\mathcal{V}))$. In comparison to the aspect ratio of the Random-Cluster model which is taken as $\rho$, we take the aspect ratio length for crossings over $\textbf{H}$ as the sum of the lengths of a square inscribed in each face, in addition to .

\bigskip

\noindent\textbf{Theorem} $2$ (\textit{hexagonal crossing probabilities}): For $x \leq  \frac{1}{\sqrt{n}}$, aspect ratio $N$ of a regular hexagon (namely a hexagon with equal side lengths) $H \subsetneq \mathcal{S}_T$, $c > 0$, and horizontal crossing $\mathcal{H}$ across $H$, estimates on crossing probabilities with free, wired or mixed boundary conditions satisfy the following criterion provided in the following four possible behaviors.

\begin{itemize}
    \item[$\bullet$] \textit{{Subcritical}}: For every $N \geq 1$, under wired boundary conditions, $\mu_{G,x,n}^{1}[  \mathcal{H} ] \leq \mathrm{exp}(-c N)$,
    
    \item[$\bullet$] \textit{{Supercritical}}: For every $N \geq 1$, under free boundary conditions, $\mu_{G,x,n}^{0}[\mathcal{H}] \geq 1 - \mathrm{exp}(-c N)    $,
    
    \item[$\bullet$] \textit{{Continuous Critical}} (\textit{Russo-Seymour-Welsh property}): For every $N \geq 1$ \text{,} independent of boundary conditions $\tau$, under the admissibility stipulation,

{\small \begin{align*}
  \big\{  \tau \textit{ is admissible for the Russo-Seymour-Welsh property} \big\} \Longleftrightarrow \big\{ \tau \textit{ is placed sufficiently far away} \\ \textit{  from $\mathcal{H}$ such that }  \mu^{\tau}_{G,x,n}[\mathcal{H}] \neq 0 \big\} , 
\end{align*} }

    \noindent $c \leq  \mu^{\tau}_{G,x,n}[\mathcal{H}] \leq 1-c$,  

    \item[$\bullet$] \textit{{Discontinuous Critical}}: For every $N \geq 1$, $\mu_{G,x,n}^{1}[\mathcal{H}] \geq 1 -  \mathrm{exp}(-c N )$ for free boundary conditions, while $\mu_{G,x,n}^{0}[\mathcal{H}] \leq \mathrm{exp}(-c N )$ for wired boundary conditions. 
\end{itemize}

\bigskip

\noindent As in the proofs for each set of inequalities in Section \textit{7} and Section \textit{9.1}, we set $p^{\mu}_* \equiv p_*$ and $q^{\mu}_* \equiv q_*$ for simplicity. Each one of the estimates below before letting $\rho \longrightarrow +\infty$ is achieved by concluding the argument with the $q^k$ 'analogy' mentioned on the previous page. Making use of $(\mathcal{S}-\mathrm{CBC})$ leads to similar estimates for the Spin measure. In the statements below, the factor Stretch appearing in the strip density and renormalization inequalities denotes some nonzero factor that a regular hexagon is “stretched” by along in the vertical degree of freedom (which is denoted with Stretch in the renormalization inequalities for crossing probabilities in the horizontal and vertical directions). In comparison to the spacing that is required for renormalization inequalities of the Random-Cluster model, for the high-temperature expansion of the loop $O(n)$ model, Stretch factors below take the form,

\begin{align*}
  \mathrm{Stretch} \equiv          \rho +    \frac{ \rho}{\sqrt{3}}            .
\end{align*}

\noindent For any factors $\mathrm{Stretch}^{\prime}<\mathrm{Stretch}$, the following renormalization inequalities for the high-temperature expansion of the loop $O(n)$ model will not hold.

\begin{figure}[H]
\begin{center}
\begin{tikzpicture}[spy using overlays={size=15mm},connect spies]
\node[regular polygon, regular polygon sides=6, minimum width=6cm,draw=blue] (reg1) at (1.2,0){};
\node[regular polygon, regular polygon sides=6, minimum width=6cm,label=side 1:$4_{j+\delta^{\prime}}$, label=side 2:$5_{j+\delta^{\prime}}$, label=side 3:$6_{j+\delta^{\prime}}$,
    label=side 4:$1_{j+\delta^{\prime}}$, label=side 5:$2_{j+\delta^{\prime}}$, label=side 6:$3_{j+\delta^{\prime}}$,draw=red] (reg2) at (2.4,0){};
    \node[regular polygon, regular polygon sides=6, minimum width=6cm,label=side 1:$4_{j+2\delta^{\prime}}$, label=side 2:$5_{j+2\delta^{\prime}}$, label=side 3:$6_{j+2\delta^{\prime}}$,
    label=side 4:$1_{j+2\delta^{\prime}}$, label=side 5:$2_{j+2\delta^{\prime}}$, label=side 6:$3_{j+2\delta^{\prime}}$,draw=gray] (reg2) at (3.6,0){};\node[above right=10pt of {(3.9,1.7)}] {$x_1$};\node[above right=10pt of {(4.1,0.4)}] {$x_2$};\node[above right=10pt of {(-0.08,-0.3)}] {$\gamma_1$};\node[above right=10pt of {(0.8,-0.3)}] {$\gamma_2$};
\draw[thick] plot [smooth,tension=1.5] coordinates{(4.1,2.2) (2.5,0.9) (1.5,1.5) (0.9,0.2) (1.0,0.3) (0.9,-1.9) (1.1,-2.4) (1.5,-2.6)};
\draw[thick] plot [smooth,tension=1.5] coordinates{(5.8,1.4) (2.4,0.7) (2.1,0.1) (2.3,-0.2) (1.3,-0.8) (1.8,-1.6) (1.9,-2.2) (2.6,-2.6)};
\draw[thick] plot [smooth,tension=1.5] coordinates{(4.1,2.2) (2.5,0.9) (1.5,1.5) (0.9,0.2) (1.0,0.3) (0.9,-1.9) (1.1,-2.4) (1.5,-2.6)};
\draw[thick] plot [smooth,tension=1.5] coordinates{(5.8,1.4) (2.4,0.7) (2.1,0.1) (2.3,-0.2) (1.3,-0.8) (1.8,-1.6) (1.9,-2.2) (2.6,-2.6)};
\draw (-3,-2.7) -- (7,-2.7){};
\draw (-3,-2.6) -- (7,-2.6){};
\spy [size=3.5cm] on (3.7,1.8)
              in node [right] at (7,2);
\end{tikzpicture}
%\begin{align*}
   %\includegraphics[width=0.95\columnwidth]%{IMG_3576.jpeg}\\
%\end{align*}
\end{center}
\caption{$0 \leq n < 2$ construction of $\mathrm{Sym}$ from macroscopic $+ \backslash -$ crossings induced by $\mathscr{C}_{2_{j}}$ and $\mathscr{C}_{2_{j + 2 \delta^{\prime}}}$. Loop configurations with distribution $\textbf{P}$, with corresponding $+/-$ random coloring of faces in $\textbf{H}$ with distribution $\mu$, are depicted above with $\gamma_1$ and $\gamma_2$. Each configuration intersects $3_{j+\delta^{\prime}}$, at points $x_1$ and $x_2$, respectively, with crossing events occurring across the box $H_j$ and its translate $H_{j + 2 \delta^{\prime}}$. Under translation invariance of the spin representation, different classes of $\mathrm{Sym}$ domains are produced from the intersection of $\gamma_1$ and $\gamma_2$, with the third side $3_{j+\delta^{\prime}}$ of the hexagon depicted above in the center. In grey, a  magnification of the region over the hexagonal lattice is provided over which the symmetric domain will be constructed. As discussed in the construction of $\mathrm{Sym}$, such symmetric domains are dependent upon the connected components of the connected component of the path $\gamma_1$ which lies above $\gamma_2$. Across the side $3_{j+\delta^{\prime}}$, one half of $\mathrm{Sym}$ is rotated to obtain the other half of the symmetric domain about the crossed edge. From paths of the connected components of each configuration, $\mathrm{Sym}$ is determined by forming the region from the intersection of the connected components of $\gamma_1$ and $\gamma_2$. }
\end{figure}

\noindent\textbf{Lemma} $1$ (\textit{Section 7}, \textit{hexagonal strip density inequalities}): In the \textbf{Non}(\textit{{Subcritical}}) and \textbf{Non}(\textit{{Supcritical}}) regimes, for every integer $\lambda \geq 2$, and every $N \in \mathrm{Stretch} \text{ } \textbf{N}$, there exists a positive constant $C$ satisfying,

\begin{align*}
     p_{\text{ } \mathrm{Stretch}\text{ } N} \geq \text{ }   \frac{1}{\lambda^C}    \text{ }  \big( q_{\mathrm{Stretch} \text{ } N} \big)^{\mathrm{Stretch} + \frac{\mathrm{Stretch}}{\lambda}}   \text{, }
\end{align*}

\noindent while a similar upper bound for vertical crossings is of the form,

\begin{align*}
    q_{\mathrm{Stretch} \text{ } N} \geq \text{ }  \frac{1}{\lambda^C}     \text{ } \big( p_{\mathrm{Stretch} \text{ } N}  \big)^{\mathrm{Stretch} + \frac{\mathrm{Stretch}}{\lambda}} \text{. }
\end{align*}

\bigskip

\noindent With the strip densities for horizontal and vertical crossings, we state closely related renormalization inequalities.

\bigskip

\noindent\textbf{Lemma} $2$ (\textit{Section 9}, \textit{hexagonal renormalization inequalities}): In the \textbf{Non}(\textit{{Subcritical}}) and \textbf{Non}(\textit{{Supcritical}}) regimes, for every integer $\lambda \geq 2$, and every $N \in \mathrm{Stretch} \text{ } \textbf{N}$, there exists a positive constant $C$ satisfying,

\begin{align*}
      p_{\mathrm{Stretch}\text{ } N} \geq   \text{ }  \lambda^C    \text{ }  \big(p_{\mathrm{Stretch} \text{ } N} \big)^{\mathrm{Stretch} \text{ } - \frac{ N \text{ }  \mathrm{Stretch}}{\lambda} }       \text{         } \text{        \&        } \text{     }  q_{\mathrm{Stretch}\text{ } N} \geq   \text{ }  \lambda^C    \text{ }  \big(q_{\mathrm{Stretch} \text{ } N} \big)^{\mathrm{Stretch} - \frac{N \text{ } \mathrm{Stretch}}{\lambda} \text{ } }  \text{ . }
\end{align*}

\section{Proofs of Theorem $1$ $\&$ of Lemma $9$}

To prove \textbf{Theorem 1}, we introduce $6$-arm crossing events,

\begin{align*}
  \mu^{\tau}_{\textbf{H},x,n} \big[  \textit{six crossings exist from each vertex of a hexagonal face belonging to $F \big( \textbf{H} \big)$}        \big]  , 
\end{align*}

\noindent under boundary conditions $\tau$, from which symmetric domains will be crossed with good probability. The arguments hold for the $n \geq 2$ extension measure with free, wired or mixed boundary conditions. Previous use of such domains has been implemented to avoid using self duality throughout the renormalization argument {[1,13]}. Although more algebraic characterizations of fundamental domains on the hexagonal, and other, lattices exist {[4]}, we focus on defining crossing events, from which we compute the probability conditioned on a path $\Gamma$ crossing the symmetric region.

\subsection{Existence of the homeomorphism $\mu$}

The increasing homeomorphism is shown to exist with the following.

\bigskip

\noindent \textbf{Lemma $8$} (\textit{Existence of the Spin Measure homeomorphism}): For any $L > 0$, there exists $c_0 = c_0(L) > 0$ so that for $nL \geq 1$, $\mu^{\tau}_{G,x,n}[\mathcal{H}] \geq c_0 \mu^{\tau}_{G,x,n}[\mathcal{V}]^{\frac{1}{c_0}}$.

\noindent \textit{Proof of Theorem $1$}. With the statement of \textbf{Lemma} $8$, for $\mu = \mu^{\tau}$ on $\mathcal{S}_{T,L}$, $\mu^{*}$ is a measure supported on dual loop configurations, from which a correspondence between horizontal and vertical hexagonal crossings is well known. By making use of \textbf{Lemma} $8$, rearrangements across the following inequality demonstrate the existence of $f$ that is stated in \textbf{Theorem 1}, as,

{\small \begin{align*}
   \mu^{0}[\mathcal{H}] \geq c_0 \mu^{1}[\mathcal{V}]^{\frac{1}{c_0}} {\Longleftrightarrow} (\mathrm{Stretch} \big)^{-1}   - \mu^1[\mathcal{V}] \geq c_0 \big( (\mathrm{Stretch} \big)^{-1}   - \mu^0[\mathcal{H}] \big)^{\frac{1}{c_0}}  {\Longleftrightarrow}   \big( (\mathrm{Stretch} \big)^{-1} \\   - \mu^1[\mathcal{V}]  \big)^{c_0} \geq c_0^{c_0} \big( (\mathrm{Stretch} \big)^{-1}  - \mu^0[\mathcal{H}] \big) \text{, } 
   \end{align*} } 
   
   \noindent where the final inequality is equivalent to,

   \begin{align*}
   \mu^0[\mathcal{H}] \leq 1 - \frac{1}{c_0^{c_0}}\big( (\mathrm{Stretch} \big)^{-1}   - \mu^1[\mathcal{V}] \big) \text{, }
\end{align*}

\noindent from the fact that $\mu^0[\mathcal{H}]+ ( \mathrm{Stretch} ) \mu^{1}[\mathcal{V}] \leq  1$. The existence of a homeomorphism satisfying $\mu(\mathcal{H}) \geq f(\mu(\mathcal{V}))$ is equivalent to $ (\mathrm{Stretch} \big)^{-1}   - \mu(\mathcal{V}) \geq f(\mu(\mathcal{V}))$, implying from the upper bound,

{\small \begin{align*}
       1 - \frac{1}{c_0^{c_0}}\big( (\mathrm{Stretch} \big)^{-1}  - \mu^1[\mathcal{V}] \big) = \frac{c_0^{c_0} - (\mathrm{Stretch} \big)^{-1}  + \mu^1[\mathcal{V}] }{c_0^{c_0}}  =   \frac{c_0^{c_0}- (\mathrm{Stretch} \big)^{-1}  }{c_0^{c_0}}  + \frac{\mu^1[\mathcal{V}]}{c_0^{c_0}} \\ \\  = 1 - c_0^{-c_0} (\mathrm{Stretch} \big)^{-1}   +  c_0^{-c_0} \mu^1[\mathcal{V}]    \text{. } 
\end{align*} }

\noindent The homeomorphism,

\begin{align*}
    f \big[ x \big] =  1 - c_0^{-c_0} (\mathrm{Stretch} \big)^{-1}   +  c_0^{-c_0} x , 
\end{align*}

\noindent for $0 < x < 1$, can be read off from the inequality, hence establishing its existence. \boxed{}

\subsection{Crossing events for Lemma $9$}

For a fixed ordering of all $6$ edges that enclose any $H_j \subset \mathcal{S}_{T,L}$, $\{1_j , 2_j , 3_j , 4_j , 5_j , 6_j \}$, crossing events $\mathcal{C}$ to obtain hexagonal symmetric domains with rotational and reflection symmetry will be defined. To obtain generalized regions from their symmetric counterparts in the plane from {[14]}, we make use of comparison between boundary condition with the $n \geq 2$ extension measure. For $\mu$, we are capable of readily proving a generalization of the union bound with the following prescription.

First, we define $5$-armed crossing events across the box $H_j \in \mathcal{S}_{T,L}$, from which families of crossing probabilities across a countable number of domains are introduced.

\noindent \textbf{Definition} \textit{$4.2$} (\textit{crossings events across the hexagonal box}) Fix $x = \{ k \}$. From a partition of $x$ into equal $k$ subintervals, each of length $\frac{s}{k}$, we define a countable family of crossing events from the partition $\mathcal{S}_j$ of $1_j$ to the corresponding topmost edge $4_j$ of $H_j$, as well as crossing events from $\mathcal{S}_j$ to all remaining edges of $H_j$. We consider crossing events across finite volumes arranged as follows,

\begin{itemize}
    \item[$\bullet$] From our choice of $1_j$, we horizontally position the line $\mathcal{L} \equiv [0,\delta] \times \{0\} \subsetneq  \textbf{H}$ for arbitrary $\delta$. We denote the horizontal translate $H_{j+ \delta^{\prime}}$ of $H_j$ along $\mathcal{L}$ by $\delta^{\prime}$ where $\delta^{\prime} << \delta$.

    \item[$\bullet$] From crossing events across the series $\{H_{j-\delta^{\prime}}, H_j , H_{j + \delta^{\prime}} \}_{ \{ \delta^{\prime}\in \textbf{R} \} }$ of hexagons we additionally introduce crossing events across translates by stipulating that the crossing starting from the partition of $\mathcal{L}$ into $k$ subintervals to any of the remaining edges $\{2,3,4,5,6 \}$ of $H$ occur in other regions, namely $H_{j - \delta^{\prime}} \cap H_j$ and $H_j \cap H_{j+\delta^{\prime}}$.

\end{itemize}

With the properties of the crossings provided above, by (FKG) and the finite energy property we conclude by sending $L \rightarrow + \infty$, generalizing the crossing events on $\mathcal{S}_T$ in the weak limit along the infinite hexagonal strip. 

Differences emerge in the proofs for the dilute Potts model in comparison to those of the random cluster model, not only in the encoding of boundary conditions for $\mu$ but also in the construction of the family of crossing probabilities, and the cases that must be considered to prove the union bound. We gather these notions below; denote the quantities corresponding to the partition $\mathcal{S}_j \subsetneq  1_j$ with the following events,

{\small \begin{align*}
        \mathscr{C}_{2_j} =   \{  \mathcal{S}_j \overset{H_{j+\delta^{\prime}}}{\longleftrightarrow}  2_{j-\delta^{\prime}} \}      \text{, } 
    \mathscr{C}_{3_j} =  \{ \mathcal{S}_j \overset{H_{j+\delta^{\prime}}}{\longleftrightarrow} 3_{j-\delta^{\prime}} \} \text{   ,       }    \text{             } \mathscr{C}_{4_j} =  \{ \mathcal{S}_j  \overset{H_j}{\longleftrightarrow} 4_j \}      \text{, } 
    \mathscr{C}_{ 5_j}= \{ \mathcal{S}_j \overset{H_{j-\delta^{\prime}}}{\longleftrightarrow} 5_{j+\delta^{\prime}} \}  \text{, } \\ \text{      } \mathscr{C}_{6_j} =  \{ \mathcal{S}_j \overset{H_{j-\delta^{\prime}}}{\longleftrightarrow} 6_{j+\delta^{\prime}} \}     \text{, } 
\end{align*} }

\noindent as well as the following dual crossing events,

% \draw[thick] plot [smooth,tension=1.5] coordinates{(0.3,-2.5) (4.1,2.2) (2.5,1.2) (1.5,1.5) (0.9,0.2) (1.0,0.3) (0.9,-1.9) (0.6,-2.4) (0.5,-2.6)};

\begin{align*}
        \mathscr{C}^{\prime}_{2_j} =  \{ \mathcal{S}_j  \underset{*}{\overset{H_{j+\delta^{\prime}}}{\longleftrightarrow}} 2_{j + \delta^{\prime}} \} \text{ } \backslash \text{ } \mathscr{C}_{2j} \text{, } 
    \mathscr{C}^{\prime}_{ 3_j}= \{ \mathcal{S}_j \underset{*}{\overset{H_{j+\delta^{\prime}}}{\longleftrightarrow}} 3_{j + \delta^{\prime}} \} \text{ } \backslash \text{ }  \mathscr{C}_{3j} \text{,       }   \text{             } \\ \mathscr{C}^{\prime}_{5_j} =  \{ \mathcal{S}_j \underset{*}{\overset{H_{j-\delta^{\prime}}}{\longleftrightarrow}} 5_{j + \delta^{\prime}} \}   \text{ } \backslash \text{ }  \mathscr{C}_{5j}  \text{, }  \mathscr{C}^{\prime}_{6_j} =  \{ \mathcal{S}_j \underset{*}{\overset{H_{j-\delta^{\prime}}}{\longleftrightarrow}} 6_{j + \delta^{\prime}} \} \text{ } \backslash  \text{ } \mathscr{C}_{6j}  \text{. } 
\end{align*}

\noindent Straightforwardly, an original crossing event and its dual are obtained through the relation,

\begin{align*}
    \big\{ A \overset{\textbf{H}}{\longleftrightarrow} B \big\} =  \big\{ A \underset{*}{\overset{\textbf{T}}{\longleftrightarrow}} B \big\}   ,
\end{align*}

\noindent for $A,B \in F \big( \textbf{H} \big)$. Before proceeding to make use of the $6$-arm events to create symmetric domains for \textbf{Lemma} $9$ (presented below), we prove \textbf{Proposition} $8$ below. We take the finite weak volume limit as the underlying graph exhausts $\textbf{H}$, as depicted in \textbf{Figure} \textit{9}, with possible choices for connected components depicted in \textbf{Figure} \textit{10} and \textbf{Figure} \textit{11}.

\bigskip

\noindent \textbf{Lemma $9$} (\textit{$6$-arm events, existence of the constant $c$}): For every $\lambda > 0$ there exists a constant $c, \lambda$ such that for every $n \in \textbf{Z}$,

\begin{align*}
          \mu[C_0]  \geq \frac{c}{\lambda^5} \mu[ \mathcal{V}^{\prime} ]^{5}  \text{ . }
\end{align*}

\bigskip

\noindent \textit{Proof of Proposition $8$}. Let $C_j = \{\mathcal{S}_j  \overset{H_j \cup H_{j+ 2 \delta^{\prime}}}{\longleftrightarrow} \mathcal{S}_{j + \delta^{\prime}} \cup \mathcal{S}_{j+2 \delta^{\prime}}  \}$.

Uniformly in boundary conditions, for \textbf{Proposition} $8$ horizontal (vertical) crossings $\mathcal{H}$ ($\mathcal{V}$) across $H_j$ can be pushed forwards under $\mu$ to obtain a standard lower bound for the probability of obtaining a longer vertical (horizontal) crossing $\mathcal{V}^{\prime}$ ($\mathcal{H}^{\prime}$) through one application of (FKG) to the finite intersection of shorter vertical (horizontal) crossings $\mathcal{H}^{\prime}_j$ ($\mathcal{V}^{\prime}_j$), 

\begin{align*}
       \mu[\mathcal{H}^{\prime}] \geq  \mu  \big[ \bigcap_{j \in \mathcal{J}} C_j \text{ }  \big]  \geq \prod_{j \in \mathcal{J}} \mu[\mathcal{V}^{\prime}_j] \geq  \big[ \frac{c}{\lambda^3} \mu[\mathcal{V}^{\prime}]^{5}\big]^{|\mathcal{J}|}   \text{ , }  
\end{align*}

\noindent where the product is taken over admissible $j \in \mathcal{J} \equiv \{j \in \textbf{R}: \exists \textit{ a regular hexagon with} \textit{ side length } j $ \&$  H_j \cap \mathcal{S}_{T,L} \neq \emptyset \}$, with $c, \lambda > 0$. This is the condition provided above in \textbf{Lemma} \textit{9}. We denote the sequence of inequalities with $(\mathrm{FKG})$ because the same argument will be applied several times for collections of horizontal and vertical crossings. From a standard lower bound from vertical crossings, the claim follows by setting $\lambda$ equal to the aspect ratio of $H_j$. \boxed{}

\bigskip

\noindent The lower bound of $(\mathrm{FKG})$ above is raised to the cardinality of $\mathcal{J}$. We apply the same sequence of terms from this inequality to several arguments in \textbf{Corollary} $11$, \textbf{Lemma} $1$,  \textbf{Lemma} $2$, \textbf{Lemma} $13$, \& \textbf{Lemma} $14$. We turn to a statement of \textbf{Lemma} $9$.

\begin{figure}[H]
\begin{center}
\includegraphics[width=0.66\columnwidth]{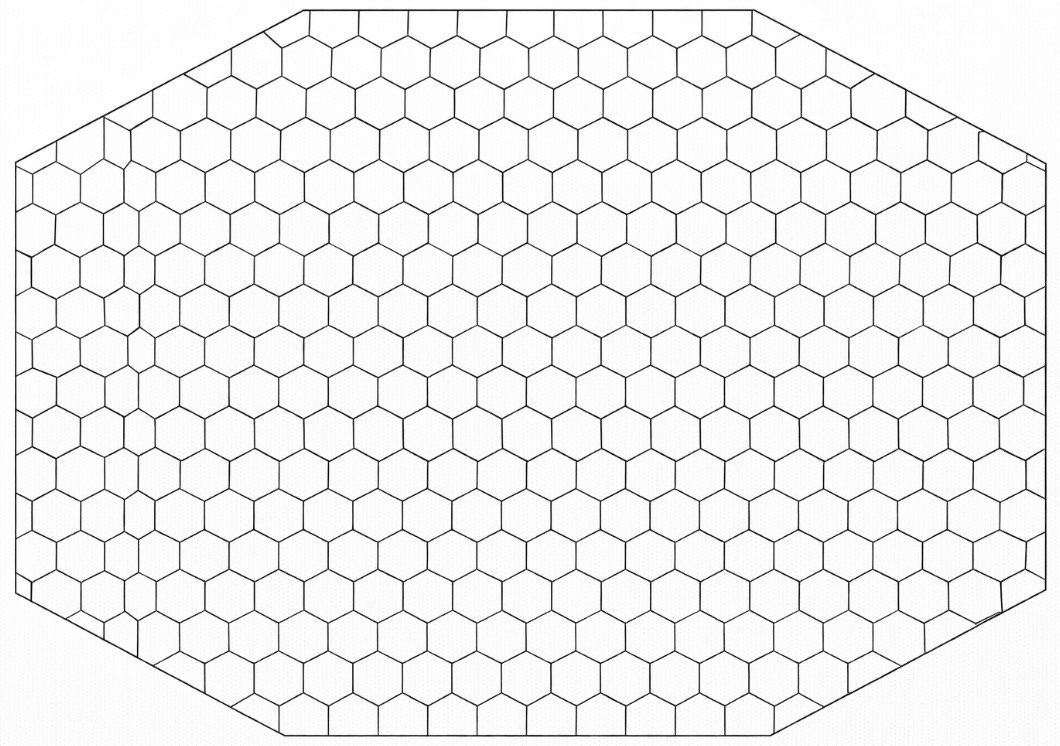}
\end{center}
\caption{A depiction of a hexagonal finite volume about which crossing events, such as those considered in Lemma $9$, will be considered below.}
\end{figure}

\section{Lemma $9$ arguments}

\textit{Proof of Lemma $9$}. For the \textit{$6$-arm lower bound}, the argument involves manipulation of symmetric domains. In particular, we must examine the crossing event that is the most probable from the union bound, in $3$ cases that are determined by the $\frac{2\pi}{3}$ rotational invariance of $\mu$. Under this symmetry, in the union bound it is necessary that we only examine the structure of the crossing events $\mathcal{C}$ in the following cases. We include the index $j$ associated with crossing events $\mathcal{C}_j$, executing the argument for arbitrary $j$ (in contrast to $j \equiv 0$ in {[14]}), readily holding for any triplet $j-\delta^{\prime},j, j+\delta^{\prime}$ which translates $H_j$ horizontally. Besides exhibiting the relevant symmetric domain in each case, the existence of $c$ will also be justified. Depending on the construction of $\mathrm{Sym}$, we either partition the outermost layer to $\mathrm{Sym}$, called the incident layer to $\partial \text{ } \mathrm{Sym}$, as well as sides of $\mathrm{Sym}$ with $L_{\mathrm{Sym}}$, $R_{\mathrm{Sym}}$, $T_{\mathrm{Sym}}$ and $B_{\mathrm{Sym}}$. Finally, we finish the section with bounds in \textit{5.4.3} to conclude the argument.

\subsection{$\mathcal{C}_j \equiv \mathscr{C}_{2_j}$}

In the first case, crossings across $2_j$ can be analyzed with the events $\mathscr{C}_{j}$ and $\mathscr{C}_{j + 2 \delta^{\prime}}$. To quantify the conditional probability of obtaining a $2_{j + \delta^{\prime}}$ crossing beginning from $\mathcal{S}_{j + \delta^{\prime}}$, denote $\Gamma_{2_j}$ and $\Gamma_{2_{j + 2 \delta^{\prime}}}$ as two independent and identically distributed random variables drawn from $\mathcal{S}_j$ and $\mathcal{S}_{j + 2 \delta^{\prime}}$ to $2_j$ and $2_{j + 2 \delta^{\prime}}$, the collections of random paths at $j$ and $j+2 \delta^{\prime}$, respectively. Denote $\gamma_1 \in \Gamma_{2_j}$, $\gamma_2 \in \Gamma_{2_{j + 2 \delta^{\prime}}}$ as two realizations of paths from $\mathcal{S}_j$ and $\mathcal{S}_{j + 2 \delta^{\prime}}$. To incorporate properties of the dilute Potts model, we also condition that the number of connected components $k_{\gamma_1}$ of $\gamma_1$ equal the number of connected components of $k_{\gamma_2}$ of $\gamma_2$ in the spin configuration sampled under $\mu$. For simplicity we denote $k_{\gamma_1} \equiv k_{\gamma_1}|_{\mathscr{C}_j  \cap   \mathscr{C}_{j + 2 \delta^{\prime}}}=1$ and $k_{\gamma_2} \equiv k_{\gamma_1}|_{\mathscr{C}_j  \cap  \mathscr{C}_{j + 2 \delta^{\prime}}}=1$. Finally, assign $\Omega \subsetneq \textbf{H}$ as the points to the left of $\gamma_1$ and to the right of $\gamma_2$, and the symmetric domain as $\mathrm{Sym} \equiv \mathrm{Sym}_{2_j} \equiv \mathrm{Sym}_{2_j}(\Omega)$. To obtain a crossing across $\mathrm{Sym}$, we conditionally pushforward the event

\begin{align*}
     \mu[C_0  |  \Gamma_{2_j} = \gamma_1 , \Gamma_{2_{j + 2 \delta^{\prime}}} = \gamma_2  ]  \text{, }
\end{align*}

\noindent which quantifies the probability of obtaining a connected component across $\mathcal{S}_{j + \delta^{\prime}} \cup \mathcal{S}_{j + 2 \delta^{\prime}}$. We condition $C_0$ through $\gamma_1$ and $\gamma_2$ because if there exits a spin configuration passing through $\mathrm{Sym}$ whose boundaries are determined by $\gamma_1$ and $\gamma_2$, then necessarily the configuration would have a connected component from $\mathcal{S}_j$ to $\mathcal{S}_{j+2} \cup \mathcal{S}_{j+4}$ hence confirming that $C_0$ occurs. To establish a comparison between this conditional probability and the conditional probability of obtaining a horizontal crossing across $\mathrm{Sym}$, consider,

\begin{align*}
         \mu[      \gamma_1 \overset{\Omega}{\longleftrightarrow} \gamma_2      | \Gamma_{2_j} = \gamma_1 , \Gamma_{2_{j + 2 \delta^{\prime}}} = \gamma_2  ]      \text{ , } 
\end{align*}

\noindent subject to wired boundary conditions on $R_{\mathrm{Sym}}$ and $L_{\mathrm{Sym}}$ and free boundary conditions elsewhere. Conditionally this probability is an upper bound for another probability supported over $\mathrm{Sym}$, as,

{\small \begin{align*}
    \mu[      \gamma_1 \overset{\Omega}{\longleftrightarrow} \gamma_2   |  \Gamma_{2_j} = \gamma_1 ,  \Gamma_{2_{j + 2 \delta^{\prime}}} = \gamma_2 , k_{\gamma_1} = k_{\gamma_2} = 1  ]    \geq \mu_{\Omega}^{\{\gamma_1 , \gamma_2\}}[    \gamma_1    \longleftrightarrow   \gamma_2    ]  \tag{$\mathrm{V}$-$\mathrm{SYM DOM}$}  \text{, }  
\end{align*} }

\noindent with the conditioning on the connected components applying to $+/-$ spin configurations, $\Omega$ is a region inside the symmetric domain, and the $\{ \gamma_1, \gamma_2\}$ superscript indicates boundary conditions wired along $\gamma_1$ and $\gamma_2$. We denote this inequality as $(\mathrm{V}$-$\mathrm{SYM DOM})$, which is short for vertical crossings across the symmetric domain that we introduce and further analyze in Section \textit{5.2.2}. Similarly, conditional on $\Gamma_{2_j} = \gamma_1 \text{ } \& \text{ } \Gamma_{2_{j + 2 \delta^{\prime}}} = \gamma_2$, $\{ \gamma_1 \overset{\Omega}{\longleftrightarrow} \gamma_2 \}$ occurs.

To quantify the probability of $\mathscr{C}_{2_{j + \delta^{\prime}}} \backslash (C_0 \cup C_2)$, conditionally that the connect components of the event not intersect those of $\mathscr{C}_{2_j} \cap \mathscr{C}_{2_{j + 2\delta^{\prime}}}$, we make use of the $(\mathcal{S}-\mathrm{SMP})$ property, which impact the boundary conditions of the symmetric domains that will be constructed, while making use of the $(\mathcal{S}-\mathrm{CBC})$ property impacts the number of paths that can be averaged over in $\Gamma_{2_j}$ and $\Gamma_{2_{j + 2 \delta^{\prime}}}$ given the occurrence of $C_0$.

\subsubsection{Incident layer of hexagons to the symmetric domain boundary}.

Under ($\mathcal{S} -\mathrm{SMP}$), we push boundary conditions away from nonempty boundary $\partial \text{ } \mathrm{Sym} \subsetneq \partial H_j$ with the edge of intersection towards $L_{\mathrm{Sym}}$, to construct $\mathrm{Sym}$ by reflecting one half of the region enclosed by the realizations $\{\gamma_1 ,  \gamma_2\} \subsetneq \big[ \mathscr{C}_{2_j} \cap \mathscr{C}_{2_{j + 2 \delta^{\prime}}} \big]$. Because the event $\mathscr{C}_{j + \delta^{\prime}}$ necessarily induces the existence of a loop configuration from $\mathcal{S}_j$ to $2_j$, under Dobrushin/mixed boundary conditions which stipulate the existence of a wired arc of length $\frac{\pi}{6}$ along $2_j$, the measure $\mu^{\tau}_{G,x,n} \big[ \cdot \big]$ over spin configurations satisfying $\mathscr{C}_{2_j}$ implies that the probability of a crossing across $\mathrm{Sym}$ supported on $\mu^{\mathrm{mix}}_{\mathrm{Sym}}$ \footnote{The $\mathrm{mix}$ boundary conditions are provided in two separate constructions of $\mathrm{Sym}$ below.}.

Formally, boundary conditions are pushed away from the boundary of $H_j$ onto boundaries of the symmetric domain as follows. 

\bigskip

\noindent \textbf{Definition} \textit{2} (\textit{pushing boundary conditions onto symmetric domains from boundary conditions on $H_j$}) From boundary conditions along $2_j$, before reflecting connected components induced by the crossings event $\mathscr{C}_{2_j}$ about $2_j$, boundary conditions along symmetric domains are obtained with the following procedure:

\begin{itemize}
    \item[$\bullet$] To partition vertices in $\mathrm{Sym}$ for constructing boundary conditions on vertices along the boundaries of symmetric domains, we assign $+$ boundary conditions to a partition of the first layer of hexagons outside of the crossing induced by $\mathscr{C}_{2_{j+\delta^{\prime}}} \backslash (C_0 \cup C_2)$, conditioned under realizations of paths $\gamma_1 \text{ } \& \text{ } \gamma_2$.
    
    \item[$\bullet$] To apply $(\mathcal{S}-\mathrm{SMP})$, given the crossing $\mathscr{C}_{2_{j+\delta^{\prime}}} \backslash (C_0 \cup C_2)$, the length of the boundary of the symmetric domain is determined by the number of connected components of the spin configuration, which corresponds to the the edges present in the configuration. From the total number of vertices on the boundary, we introduce boundary conditions with $(\mathcal{S}-\mathrm{CBC})$. Outside of $H_j$, the paths $\gamma_1$ and $\gamma_2$.

    \item[$\bullet$] After having identified the boundaries of the symmetric domain, reflection of one half of $\mathrm{Sym}$ is constructed by taking the union $\gamma^{x_{\gamma_1 , \gamma_2}}_1 \cup \gamma^{x_{\gamma_1, \gamma_2}}_2$, where the paths in the union denote the restriction of the connected components of $\gamma_1$ and $\gamma_2$ after $\mathscr{C}_j$ and $\mathscr{C}_{j + 2 \delta^{\prime}}$ have occurred. The remaining top half of $\mathrm{Sym}$ is obtained by reflection through $2_j$ that was crossed by $\gamma_1$ and $\gamma_2$, as with the remaining half of the lower part (the connected components of $\gamma_2$ constitute one half of the lower region of $\mathrm{Sym}$).
    
    \item[$\bullet$] The reflections $\tilde{\gamma}^{x_{\gamma_1 , \gamma_2}}_1$ and $\tilde{\gamma}^{x_{\gamma_1 , \gamma_2}}_2$ described in previous steps provide the remaining half of $\mathrm{Sym}$ after performing appropriate reflections with respect to the intersection with the side of the hexagonal box.
\end{itemize}

\subsection{$(\mathcal{S}-\mathrm{SMP})$ property}

We progress towards making use of another modification for the dilute Potts model through the symmetric domain construction previously described to ensure that such domains are conditionally bridged with good probability. 

\subsubsection{Upper bound for conditional crossing events across symmetric domains}. 

To proceed, we make use of $\mathrm{Sym}$, in addition to the modification of boundary conditions as follows. From an application of $(\mathcal{S}-\mathrm{CBC})$, the conditional probability introduced at the beginning of the proof, under spin configurations supported on $\mu_{\mathrm{Sym}}$ satisfies, under the conditional measure $\mu_{\Omega} \equiv \mu_{\Omega}[  \cdot |_{\Omega}  |      \gamma_1 \cap \gamma_2 = \emptyset , \gamma_1 \cap \gamma_3 = \emptyset  , k_{\gamma_1} = k_{\gamma_2}   = 1 ]$, for measurable events depending on finitely many edges in $\Omega$,

\begin{align*}
    \mu[  \mathscr{C}_{2_j} \backslash ( C_0   \cup   C_2 )    |  \Gamma_{2_j} = \gamma_1 ,  \Gamma_{2_{j + 2 \delta^{\prime}}} = \gamma_2  ] \leq \mu_{\Omega}^{\{\gamma_1 , \gamma_2\}^{c}} [    \mathscr{C}_{2_{j + \delta^{\prime}}  }   ]    \text{ , }
\end{align*}

\noindent after examining the pushforward of the conditional probability above under spin configurations supported in $\mathrm{Sym}$, where the superscript $\{\gamma_1 , \gamma_2 \}^{c}$ denotes free boundary conditions along $\gamma_1$ and $\gamma_2$ and wired elsewhere, the complement of $\{\gamma_1 , \gamma_2 \}$ given in the lower bound of $(\mathrm{V}$-$\mathrm{SYM DOM})$ (provided in Section \textit{5.1}). The stochastic domination above of the conditional probability under no boundary conditions on any side of $\mathrm{Sym}$ will be studied for paths $\gamma_3 \in \Gamma_{j + \delta^{\prime}}$. The event under $\mu_{\Omega}^{\{\gamma_1 , \gamma_2\}^{c}}$ demands that the connected components of $\gamma_3$ be disjoint for those of $\gamma_1$ and $\gamma_2$ for the entirety of the path.

Particularly, we remove the conditioning from the pushforward in the upper bound because the definition of $\Omega$ implies that connectivity holds in between $\gamma_1$ and $\gamma_2$. Pointwise, the connected components of $\gamma_3$ do not intersect those of $\gamma_1$ and $\gamma_2$. Recalling $(\mathrm{V}$-$\mathrm{SYM DOM})$ in Section \textit{5.1}, we present additional modifications to the renormalization argument through the lower bound of the inequality to exhaust the case for $\mathcal{C}_j \equiv \mathscr{C}_{2_j}$. Lower bounds for the pushforward under $\mu_{\Omega}$ can only be obtained for mixed boundary conditions along $\mathrm{Sym}$ precisely under partitions of the incident hexagonal layer given in Section \textit{5.1.1} $\&$ Section \textit{5.1.2}.

Under the conditions of $(\mathcal{S}-\mathrm{SMP})$, crossings in $\Omega$ with boundary conditions $\{\gamma_1 , \gamma_2 \}$, the lowermost bound for $(\mathrm{V}$-$\mathrm{SYM DOM})$ can only be established when boundary conditions are distributed under Section \textit{5.1.1} or Section \textit{5.1.2}. For completeness, we first establish the lower bound for Section \textit{5.1.1}, in which the boundary conditions for a crossing distributed under $\mu_{\mathrm{Sym}}^{\{\gamma_1 , \gamma_2 \}}$ can be compared to a closely related crossing distributed under $\mu_{\Omega}^{\{ \gamma_1 , \gamma_2 \}^c}$. 

To establish the comparison, the edges in $\mathrm{Sym} \cap \Omega$, we divide the proof into separate cases depending on whether the boundary conditions for vertcies along $\gamma_1$ or $\gamma_2$ are connected together under wired or free boundary conditions. One instance of pushing boundary conditions occurs for $\mathcal{C}_j \equiv \mathscr{C}_{2_j}$, while another instance of pushing boundary conditions occurs when $\mathcal{C}_j \equiv \mathscr{C}_{4_j}$ in Section \textit{5.4}. \footnote{In contrast to the planar case of {[14]}, considerations through the condition $k_{\gamma_1} = k_{\gamma_2} = 1$ impact the construction of $\mathrm{Sym}$ and the rotational symmetry the region enjoys.}

\subsubsection{Pushing wired boundary conditions away from $\Omega$ towards $\mathrm{Sym}$} .

One situation occurs as follows. It is possible that $+ \backslash -$ configurations distributed under $\mu_{\Omega}$ can be compared to configurations distributed under $\mu_{\mathrm{Sym}}$ by pushing boundary conditions away from the first partition of $\mathrm{Sym}$ towards $\Omega$; applying $(\mathcal{S} -\mathrm{CBC})$ between deterministic and random circuits yields

\begin{align*}
           \mu_{\Omega}^{\{\gamma_1 , \gamma_2 \}^c}[       \mathscr{C}_{2_{j + \delta^{\prime}}}     ]  \leq \mu_{\mathrm{Sym}}^{\{T,B\}}[     \mathscr{C}_{2_{j + \delta^{\prime}}}                ]            \text{, }   
\end{align*}

\noindent by virtue of monotonicity in the domain be cause $\Omega \subsetneq  \mathrm{Sym}$, where $\mu_{\Omega}$ is taken under boundary conditions $\{T,B\}$ wired along $T_{\mathrm{Sym}}$ and $B_{\mathrm{Sym}}$. Additionally, the comparison

\begin{align*}
        \mu_{\mathrm{Sym}}^{\{T,B\}}[     \mathscr{C}_{2_{j + \delta^{\prime}}}     ] \leq \mu_{\mathrm{Sym}}^{\{T,B\}}[   T_{\mathrm{Sym}} \longleftrightarrow B_{\mathrm{Sym}}   ]    \text{, }  
\end{align*}

\noindent holds by virtue of $(\mathrm{FKG})$ for the Spin measure, in which we suitably restricted our analysis of $\mu$ for $n \geq 1$ $\&$ $nx^2 \leq 1$, from which it follows that the event $\{T_{\mathrm{Sym}} \longleftrightarrow B_{\mathrm{Sym}} \}$ depends on more edges than the conditional event $\{ \mathscr{C}_{2_{j + \delta^{\prime}}} | \gamma_1 \cap \gamma_2 = \emptyset , \gamma_1 \cap \gamma_3 = \emptyset , k_{\gamma_1} = k_{\gamma_2} = 1 \}$ under $\mu_{\Omega}$ does and is an increasing event. Finally, the simplest comparison, namely the equality

\begin{align*}
        \mu_{\mathrm{Sym}}^{\{T,B\}}[  T_{\mathrm{Sym}} \longleftrightarrow B_{\mathrm{Sym}}    ]   = \mu_{\mathrm{Sym}}^{ \{L,R\}}[     L_{\mathrm{Sym}} \longleftrightarrow R_{\mathrm{Sym}}    ]        \text{, }
\end{align*}

\noindent holds by virtue of dual boundary conditions of $\mu_{\mathrm{Sym}}$, in which the pushforward of the event $\{T_{\mathrm{Sym}} \longleftrightarrow B_{\mathrm{Sym}} \}$ under boundary conditions ${\{T,B \}}$ is equal to the pushfoward of the event $\{ L_{\mathrm{Sym}} \longleftrightarrow R_{\mathrm{Sym}}\}$ under boundary conditions ${\{L,R\}}$. Hence by rotational invariance of boundary conditions of $\mathrm{Sym}$ one obtains the following upper bound,

{\small \begin{align*}
         \mu[  \mathscr{C}_{2_j} \backslash ( C_0   \cup   C_2 )   |  \Gamma_{2_j} = \gamma_1 , \Gamma_{2_{j + 2 \delta^{\prime}}} = \gamma_2 , k_{\gamma_1} = k_{\gamma_2} =1 ] \leq \mu_{\Omega}^{\{\gamma_1 , \gamma_2 \}^c}[   \mathscr{C}_{2_{j + \delta^{\prime}}}        ]        \text{, } 
\end{align*} }

\noindent which holds by $(\mathcal{S}-\mathrm{SMP})$, as wired boundary conditions for $\mathscr{C}_{2_j}$ in between $\gamma_1$ and $\gamma_2$ can be pushed away to obtain wired boundary conditions along $\gamma_1$ and $\gamma_2$ for $\mathscr{C}_{2_{j + \delta^{\prime}}}$, in turn transitively yielding,

{\small \begin{align*}
       \mu[  \mathscr{C}_{2_j} \backslash ( C_0   \cup   C_2 )   | \Gamma_{2_j} = \gamma_1, \Gamma_{2_{j + 2 \delta^{\prime}}} = \gamma_2 , k_{\gamma_1} = k_{\gamma_2} = 1  ] \leq \mu_{\mathrm{Sym}}^{ \{L,R\}}[    L_{\mathrm{Sym}} \longleftrightarrow R_{\mathrm{Sym}}    ]      \text{. }   
\end{align*} }

Under $\frac{2 \pi}{3}$ rotational invariance of $\mu$, the argument for this case can be directly applied with $\mathcal{C}_j \equiv \mathscr{C}_{5j}$. Examining the pushforward of this crossing event, in addition to $\mathscr{C}_{5_{j - \delta^{\prime}}}$ which guarantees the existence of a connected component that necessarily crosses $5_j$ through $5_{j - \delta^{\prime}}$, leads to the same conclusion with wired boundary conditions from  to  along $\mathrm{Sym}$. Under duality, the identification between measures under nonempty boundary conditions over $\mathrm{Sym}$ readily applies. Hence a combination of $(\mathcal{S}-\mathrm{SMP})$, followed by $(\mathcal{S}-\mathrm{CBC})$, implies that $\{L_{\mathrm{Sym}} \longleftrightarrow R_{\mathrm{Sym}} \}$ occurs with substantial probability for $\mathcal{C}_j \equiv \mathscr{C}_{2_j}$ and $\mathcal{C}_j \equiv \mathscr{C}_{5_j}$.

\subsection{$\mathcal{C}_j \equiv \mathscr{C}_{3_j}$}

In the second case, one can apply similar arguments with the following modifications. To identify other possible symmetric regions $\mathrm{Sym}$ corresponding to $\mathscr{C}_{3_j}$ and $\mathscr{C}_{3_{j + 2 \delta^{\prime}}}$, fix path realizations $\gamma_1 \in \Gamma_{3_j}$ and $\gamma_2 \in \Gamma_{3_{j + 2\delta^{\prime}}}$. From $\gamma_1$ and $\gamma_2$, we construct $\mathrm{Sym}$ by reflecting half of the domain across $3_j$ instead of $2_j$. Under $\frac{2 \pi}{3}$ rotational invariance of $\mu$, $\mathrm{Sym}$ constructed in this case correspond to symmetric domains induced by the paths in $\mathscr{C}_{5_j}$ and $\mathscr{C}_{5_{j + 2 \delta^{\prime}}}$. Explicitly, the conditional probability is of the familiar form,

\begin{align*}
           \mu[  \mathscr{C}_{3_j} \backslash ( C_0   \cup   C_2 ) |  \Gamma_{2_j} = \gamma_1 ,  \Gamma_{2_{j + 2 \delta^{\prime}}} = \gamma_2  ]        \text{, } 
\end{align*}

\noindent which by the same argument applied to $\mathscr{C}_{3_j}$ is bounded above by

\begin{align*}
         \mu_{\mathrm{Sym}}^{ \{L,R\}}[     L_{\mathrm{Sym}} \longleftrightarrow R_{\mathrm{Sym}}   ]              \text{, } 
\end{align*}

\noindent for $\mathrm{Sym}(\Omega) \equiv \mathrm{Sym}$. Applying the same argument to push bondary conditions away from wired boundary conditions on $3_j$ ($5_j$), to $L_{\mathrm{Sym}}$ ($R_{\mathrm{Sym}}$) establishes the same sequence of inequalities, through contributions of $\mu, \mu \text{ } \& \text{ } \mu_{\mathrm{Sym}}$. $\mathrm{Sym}$ for $\mathscr{C}_{3_j}$ corresponds to rotating the crossings of loop configurations, and hence the symmetric region to $5_j$ from the symmetric domain corresponding to $2_j$.

\subsection{$\mathcal{C}_j \equiv \mathscr{C}_{4_j}$}

In the third case, we denote the events $C_{0}$ and $C_{2}$ as bottom to top crossings, respectively across $H_j$ and $H_{j + 2 \delta^{\prime}}$, with respective path realizations $\Gamma_1$ and $\Gamma_1$ as in the previous two cases. However, the final case for top to bottom crossings stipulates that the construction of $\mathrm{Sym}$ independently of $\Omega$. We present modifications to the square symmetric region of {[14]}, and partition the region over which connectivity events are quantified through points to the left and right of $\gamma_1$ and $\gamma_2$, respectively. In particular, we denote $\Omega$ as the collection of all points in the hexagonal box $\mathrm{Sym}$, along with the partition $\Omega = \Omega_L \cup \Omega_{(L \cup R)^c} \cup \Omega_R$. In the partition, each set respectively denotes the points to the left of $\gamma_1$, the points in between the left of $\gamma_1$ and the right of $\gamma_2$, and the points to the right of $\gamma_2$. With some abuse of notation we restrict the paths in $\Omega_L$, $\Omega_R$ and $\Omega_{(L \cup R)^c}$ to coincide with crossings in between the top most edge of $H_j$ and $\mathrm{Sym}$, in which $\Omega_R \equiv (\mathrm{Sym} \cap H_j) \cap \Omega_R$, $\Omega_L \equiv (\mathrm{Sym} \cap H_j) \cap \Omega_L$, and $\Omega_{(L \cup R)^c} \equiv (\mathrm{Sym} \cap H_j) \cap \Omega_{(L \cup R)^c}$. We provide such an enumeration to apply $(\mathcal{S}-\mathrm{SMP})$ and then $(\mathcal{S}-\mathrm{CBC})$, when comparing the spin representation measures supported over $\Omega$ and $\mathrm{Sym}$.

Besides the $\Omega$ partition, to apply $(\mathcal{S}-\mathrm{SMP})$ we examine $\mathcal{R}_1 \equiv ( \mathrm{Sym} \backslash H_j^c) \text{ }  \cap \Omega_L$ and $\mathcal{R}_2 \equiv ( \mathrm{Sym} \backslash  H_j^c) \text{ } \cap \Omega_R$ which denote the collection of points to the left of $\gamma_1$ and to the right of $\gamma_2$ in the region above $H_j$ that is contained in $\mathrm{Sym}$. To apply $(\mathcal{S}-\mathrm{CBC})$, it is necessary that we isolate $\mathcal{R}_1$ and $\mathcal{R}_2$ so that $(\mathcal{S}-\mathrm{CBC})$ can be applied to the outermost layer of hexagons incident to $\partial\Omega$ through a partition of the incident layer. 

We provide an upper bound for the pushforward of the following conditional probability, where,

\begin{align*} \mathcal{R} \equiv \bigg\{ \big\{ \Gamma_{2_j} = \gamma_1 \big\} , \big\{ \Gamma_{2_{j + 2 \delta^{\prime}}} = \gamma_2 \big\} ,  \big\{ \gamma_1 \cap \gamma_3 = \emptyset \big\} ,  \big\{ \gamma_1 \cap \gamma_2 = \emptyset  \big\} \bigg\} , 
\end{align*}

\noindent one defines,

\begin{align*}
          \mu[  \mathscr{C}_{4_j} \backslash ( C_0   \cup   C_2 )   |  \mathcal{R}  ]                 \text{, } 
\end{align*}

\noindent for the class of hexagonal box symmetric domains $\mathrm{Sym}$, with $\gamma_3 \in \Gamma_{j + \delta^{\prime}}$.

\subsubsection{Pushing boundary conditions away from $H_j$ towards $\mathrm{Sym}$}.

Next, we push boundary conditions away from $H_j$. Under the assumption that the upper half of $\mathrm{Sym}$ is endowed with wired boundary conditions while the lower half is endowed with free boundary conditions. We denote these boundary conditions with $\mathrm{Top}$, and will consider the measures supported over $\mathrm{Sym}$, respectively. From observations in previous cases, to analyze the conditional probability of $C_0$ given $\Gamma_j = \gamma_1$ and $\Gamma_{j + 2 \delta^{\prime}} = \gamma_2$, we introduce the following lower bound for a connectivity event between $\gamma_1$ and $\gamma_2$ in $\Omega$, with,

\begin{align*}
        \mu[ C_0  |    \mathcal{R}  ] \geq    \mu[     \gamma_1 \overset{\Omega}{\longleftrightarrow} \gamma_2      |  \mathcal{R}  ]       \text{, }        
\end{align*}

\noindent holds from arguments applied when $\mathcal{C}_j \equiv \mathscr{C}_{2_j}$. By construction, $H_j \subsetneq  \mathrm{Sym}$ implies

\begin{align*}
       \mu[      \mathscr{C}_{4_j}         |    \mathcal{R}    ]      \leq \mu_{\Omega_{(L \cup R)^c}}^{\{\gamma_1 , \gamma_2 \}^c}[    \mathscr{C}_{4_j}         |     \mathcal{R}   ]      \text{, } 
\end{align*}

\noindent due to montonicity in the domain, as the occurrence of $\mathscr{C}_{4_j}$ conditionally on disjoint connected components of $\gamma_3 \in \Gamma_{j + \delta^{\prime}}$ with those of $\gamma_1$ and $\gamma_2$. In comparison to the conditioning applied through $k_{\gamma_1} = k_{\gamma_2} = 1$ for $\mathscr{C}_{2_j}$ and $\mathscr{C}_{3_j}$, the sides of $\mathrm{Sym}$ are formed independently of the connected components of $\gamma_1$ and $\gamma_2$; a combination of montonicity of $\mu$, in addition to $(\mathcal{S}-\mathrm{SMP})$ through an equal partition of the incident layer outside of $\mathrm{Sym}$ equally into two sets along which $+ \backslash -$ spin is constant.

After pushing boundary conditions towards $\mathrm{Sym}$, we make use of rotational symmetry of $\mathrm{Sym}$. In particular, the distribution of boundary conditions from the incident layer partition of Section \textit{5.4.1} satisfies the following inequality,

{\small \begin{align*}
            \mu_{\Omega_{(L \cup R)^c}}^{\{\gamma_1 , \gamma_2\}^c}[       C_0    |  \mathcal{R}    ]     \leq        \mu_{\mathrm{Sym}}^{(\mathrm{Top \text{ } Half})}[   C_0    |  \mathcal{R}    ]    \leq   \mu_{\mathrm{Sym}}^{(\mathrm{Top \text{ } Half})} [    T_{\mathrm{Sym}} \longleftrightarrow B_{\mathrm{Sym}}  ]  = \mu_{\mathrm{Sym}}^{(\mathrm{Top \text{ } Half})^{\frac{2 \pi}{3}}}[  L_{\mathrm{Sym}} \longleftrightarrow \\ R_{\mathrm{Sym}}     ]   \text{, }
\end{align*} }

\noindent where $(\mathrm{Top \text{ } Half})$ denotes wired boundary conditions along the top half of hexagonal $\mathrm{Sym}$. Within the sequence of inequalities, the leftmost lower bound for $\mu_{\mathrm{Sym}}^{\{ L, R\}} [ \text{ }  C_0 \text{ }  | \text{ } \mathcal{R}    \text{ }  ]  $ holds because $\Omega_{(L \cup R)^c} \subset \mathrm{Sym}$, with $\{L,R \}$ denoting wired boundary conditions along $L_{\mathrm{Sym}}$ and $R_{\mathrm{Sym}}$. \footnote{In contrast to square symmetric domains of {[14]} for the random cluster model, hexagonal $\mathrm{Sym}$ have two left sides and two right sides, and in turn require that boundary conditions along $\mathrm{Sym}$ be rotated by a different angle than $\frac{\pi}{2}$.} The next lower bound for $\mu_{\mathrm{Sym}}^{\{L,R \}}[    T_{\mathrm{Sym}} \longleftrightarrow B_{\mathrm{Sym}}    ]$ holds because the event $\{T_{\mathrm{Sym}} \overset{\mathrm{Sym}}\longleftrightarrow B_{\mathrm{Sym}} \}$ depends on finitely many more edges in $\textbf{H}$ than $\{C_0 | \mathcal{R} \}$ does. Finally, the last inequality holds due to complementarity as in the argument for $\mathcal{C}_j \equiv \mathscr{C}_{2_j}$. $\{L,R\}$ denotes a $\frac{2 \pi}{3}$ rotation of the boundary conditions supported over $\mathrm{Sym}$. 

More specifically, rotating the boundary conditions $\{L,R \}$ by $\frac{2\pi}{3}$ to obtain the boundary conditions $\{L,R \}^{\frac{2 \pi}{3}}$ amounts to four $\frac{\pi}{6}$ rotations of $\mathrm{Sym}$. With each rotation, the boundary conditions $\{L,R \}^{\frac{\pi}{6}}$ are obtained by rotating the partition of the incident layer along $\partial \text{ } \mathrm{Sym}$ to its leftmost neighboring edge, in addition to modifications of the connectivity in $\mathscr{C}_{4_j}$.

Finally, the arguments imply the same result as in other cases, in which

\begin{align*}
        \mu[  \mathscr{C}_{4_j} \backslash ( C_0   \cup   C_2 )    |\mathcal{R}   ]     \leq   \mu_{\mathrm{Sym}}^{(\mathrm{Top \text{ } Half)}^{\frac{2\pi}{3}}}[ L_{\mathrm{Sym}} \longleftrightarrow R_{\mathrm{Sym}}      ]     \text{. }
\end{align*}

\bigskip

We conclude the argument for \textbf{Lemma} $9^{*}$, not only having shown that the same inequality holds for a different classes of symmetric domains in the $\mathcal{C}_j \equiv \mathscr{C}_{4_j}$ case, but also that rotation of boundary conditions wired along the top half of $\mathrm{Sym}$ for top to bottom crossings can be used to obtain boundary conditions for left to right crossings. 

\subsubsection{$(\mathcal{S} -\mathrm{CBC})$ lower bound for the conditional crossing event $\{  C_0   |  \mathcal{C}_0 \cap \mathcal{C}_4 \}$}.

We complete the argument by providing the following inequalities for each case. We make use of the special case of $(\mathcal{S}- \mathrm{CBC})$ from Section \textit{3.2}, in which for $\mathcal{C}_j \equiv \mathscr{C}_{2_j}$, conditionally on top to bottom crossings $\mathcal{C}_0$ and $\mathcal{C}_2$ from Section \textit{5.2.2}, the pushforward below satisfies,

{\small \begin{align*}
      \mu [   C_0    |  \mathcal{C}_0  \cap  \mathcal{C}_4    ]     \geq         \frac{1}{n ^{ k^{\tau_{H_0 \cup H_{2 \delta} } } (\sigma) -     k^{\tau_{H_1 \cup H_{1+\delta}}  } (\sigma)  } x \text{ } \mathrm{exp}(h) \text{ }  }      \text{ }  \mu  [      \mathcal{C}_2 \backslash (  C_0 \cup C_2    )   |        \mathcal{C}_0 \cap \mathcal{C}_4      ] \tag{($\mathcal{S}$-SMP)- ($\mathcal{C}_2$)}  \text{, } 
\end{align*} }

\noindent where the normalization to the crossing probability in the lower bound is dependent on the edge weight $x$. One obtains the same bound for crossings $\mathcal{C}_j \equiv \mathscr{C}_{4_j}$. The bound above corresponds to the partition of boundary conditions. Finally, the existence of $c$ such that the inequality in the statement of $\textbf{Lemma}$ $9$ holds is of the form. For the following quantities obtained below, denote,

{\small \[
\left\{\!\begin{array}{ll@{}>{{}}l}
  k^{\tau_{H_0 \cup H_{2 \delta} } } (\sigma)   \equiv \textit{Number of connected components $k$ associated with boundary conditions $\tau$ sup-} \\ \textit{ported over $H_0 \cup H_{2 \delta}(\sigma) $ for the spin configuration $\sigma$}   , \\ \\        k^{\tau_{H_1 \cup H_{1+\delta}}  } (\sigma)  \equiv \textit{Number of connected components $k$ associated with boundary conditions $\tau$ sup-} \\ \textit{ported over $H_1 \cup H_{1+ \delta}(\sigma) $ for the spin configuration $\sigma$}   , \\  \\    e^{\tau_{H_0 \cup H_{2 \delta}}}(\sigma)   \equiv      \textit{Edges $e$ associated with boundary conditions $\tau$ supported over $H_0 \cup H_{2 \delta}$ for the }  \\ \textit{spin configuration $\sigma$}     , \\  \\  e^{\tau_{H_1 \cup H_{1+\delta}}}(\sigma)   \equiv   \textit{Edges $e$ associated with boundary conditions $\tau$ supported over $H_1 \cup H_{1+\delta}$ for }  \\ \textit{the spin configuration $\sigma$}          , \\ \\  r_{H_0 \cup H_{\delta} }(\sigma_1)  \equiv     \textit{Summation of the spins of the face configuration $\sigma_1$ over $H_0 \cup H_{\delta}$}    , \\ \\        r_{H_1 \cup H_{1+\delta}}(\sigma_2)   \equiv \textit{Summation of the spins of the face configuration $\sigma_2$ over $H_1 \cup H_{1+\delta}$}    , \\ \\ r^{\prime}_{H_0 \cup H_{\delta}}(\sigma_1) \equiv \textit{Summation over the number of monochromatically colored triangles of $\sigma_1$} \\ \textit{over $H_0 \cup H_{\delta}$} , \\ \\ r^{\prime}_{H_1 \cup H_{1+\delta}}(\sigma_2) \equiv \textit{Summation over the number of monochromatically colored triangles of $\sigma_2$} \\ \textit{over $H_1 \cup H_{1+\delta}$}  . 
\end{array}\right.
\] }

\bigskip

\noindent To obtain several factors in the constant for the lower bound, the first of which is proportional to the reciprocal of the number of connected components and positions at which the connected components are located on the lattice, write,

\begin{align*}
    n^{ \big(  k^{\tau_{H_0 \cup H_{2 \delta}  } } (\sigma) \cap  k^{\tau_{H_1 \cup H_{1+\delta}} }(\sigma)    \big) }  x^{  \big(   e^{\tau_{H_0 \cup H_{2 \delta}}}(\sigma) \cap e^{\tau_{H_1 \cup H_{1+\delta}}}(\sigma)  \big)  } \text{, } 
\end{align*}

\noindent in addition to the reciprocal of the following factor dependent on the magnitude of the two external fields $h_1$, $h_2$, for which one writes,

\begin{align*}
 e^{ h \big(   r_{H_0 \cup H_{\delta} }(\sigma_1) - r_{H_1 \cup H_{1+\delta}}(\sigma_2)      \big) + \frac{h^{\prime}}{2} \big(   r^{\prime}_{H_0 \cup H_{\delta}}(\sigma_1) - r^{\prime}_{H_1 \cup H_{1+\delta}}(\sigma_2)  \big)} \text{. } 
\end{align*}

\noindent Fix some $I$ taken to be sufficiently large. Denoting the product of all the factors above as $\mathcal{P}$, the desired lower bound is of the form,

\begin{align*}
          c =  (3 I)^{-3}  \mathcal{P} \text{. } 
\end{align*}

\noindent because the superposition of crossing probabilities,

\begin{align*}
c  \text{ }   \mu[C_0] \text{ } +  \text{ }  \mu[C_2]  \text{, } 
\end{align*}

%   n^{\# \bigg\{  k^{\tau_{H_0 \cup H_{2 \delta}  } }(\sigma) \cap k^{\tau_{H_1 \cup H_{1+\delta}}} (\sigma)  \bigg\}       }  x^{ \# \bigg\{   e_{\tau_{H_0 \cup H_{2 \delta}}}(\sigma) \cap e_{\tau_{H_1 \cup H_{1+\delta}}}(\sigma)  \bigg\}  } \times \cdots \\ e^{ h \big(   r_{H_0 \cup }(\sigma_1) - r_{H_1 \cup H_1+\delta}(\sigma_2)      \big) + \frac{h^{\prime}}{2} \big(   r^{\prime}_{H_0 \cup H_{\delta}}(\sigma_1) - r^{\prime}_{H_1 \cup H_{1+\delta}}(\sigma_2)          \big)  }   

\noindent where the first crossing probability is magnified with respect to the product of edges, loops, and the exponential of the difference between the external fields that are scaled with respect to the summation of spins, in addition to the number of monochromatically colored hexagons. We denote the spin configurations $\sigma_1$ and $\sigma_2$ supported over crossings over $H_0 \cup H_{\delta}$ and $H_1 \cup H_{1+\delta}$, respectively. To provide a lower bound for this superposition of crossing probabilities, we bound each term below with,

\begin{align*}
\mu[  C_0 \cap  \mathscr{C}_1     \cap \mathscr{C}_2 ] + \mu[ ( C_0 \cup  C_2 ) \cap        \mathscr{C}_1     \cap    \mathscr{C}_2      ] \text{, } 
\end{align*}

\noindent where the crossing events $\mathscr{C}_1 \neq \mathscr{C}_2$ are disjoint and can be chosen from the crossing events defined in Section \textit{4.2}. The lower bound holds because the crossing probability event $C_0$ depends on more edges than the event $C_0 \cap \mathscr{C}_1 \cap \mathscr{C}_2$ does, with the same observation holding between the crossing events $C_2$ and $(C_0 \cup C_2) \cap \mathscr{C}_1 \cap \mathscr{C}_2$. Moreover, the superposition provided in the lower bound itself can be bounded below by the pushforward of the single crossing event,

\begin{align*}
        \mu[  \mathcal{C}_0 \cap \mathcal{C}_2 \cap \mathcal{C}_4       ]  \geq  \mu[ \mathcal{C}_0  ]^3  \text{, }  
\end{align*}

\noindent from which the form of the constant $c$ defined above follows, due to a previous application of $(\mathcal{S}-\mathrm{SMP})$ in the lower bound for the conditional
crossing event provided in ($\mathcal{S} -\mathrm{SMP}$)-($\mathcal{C}_2$). The ultimate inequality follows from the fact that each crossing event in the intersection is independently bounded below by the product of three crossing probabilities. As a result, the union bound,

\begin{align*}
      \mathrm{max} \big\{  \mu[     \mathscr{C}_{2_j}    ]  , \mu[     \mathscr{C}_{3_j}    ]  , \mu[     \mathscr{C}_{4_j}   ]    : 0 \leq j < I      \big\}   \geq \frac{\mu[    \mathcal{V}_{\mathcal{H}}   ]}{3I}  \text{, } 
\end{align*}

\noindent holds, which is equivalent to the maximum of the crossing events taken over $\mathscr{C}_{2_0}$, $\mathscr{C}_{3_0}$ and $\mathscr{C}_{4_0}$ ,

\begin{align*}
        \mathrm{max} \big\{ \mu[ \mathscr{C}_{2_0} ] , \mu[   \mathscr{C}_{3_0} ] , \mu[  \mathscr{C}_{4_0}   ]          \big\} \geq     \frac{\mu [ \mathcal{V}_{\mathcal{H}} ]}{3I}      \text{, } 
\end{align*}

\noindent which holds given that the maximum of each one of the crossing probabilities across $2_j$ or $3_j$, across $3I$ events, in addition to the fact that the lower bound for the intersection of crossing events below takes the form,

\begin{align*}
    \mu[    C_0       ]     \geq   \mu [  \mathscr{C}^1          \cap \mathscr{C}^2     ]    \geq   \frac{\mu[ \mathcal{V}_{\mathcal{H}}  ]^2}{9I^2}            \text{, }   
\end{align*}

\noindent where in the upper bound $\mathscr{C}^1$ can be any one of the crossings to a rightmost edge of $\mathcal{H}$, $\mathscr{C}_{5_j}$ or $\mathscr{C}_{6_j}$ , while $\mathscr{C}^2$ can be any one of the crossings to a leftmost edge of $\mathcal{H}$, $\mathscr{C}_{2_j}$ or $\mathscr{C}_{3_j}$. This concludes the argument for deterministic explorations of symmetric domains. Below, we describe how to make use of the above arguments for stochastic explorations of symmetric domains.

\bigskip

\noindent Below, we demonstrate how arguments obtained previouslhy in this section can be applied to the crossing events $\mathscr{C}^{\prime}_{2_j}, \mathscr{C}^{\prime}_{3_j}$ and $\mathscr{C}^{\prime}_{4_j}$, respectively.

\begin{figure}[H]
\begin{center}
\includegraphics[width=0.82\columnwidth]{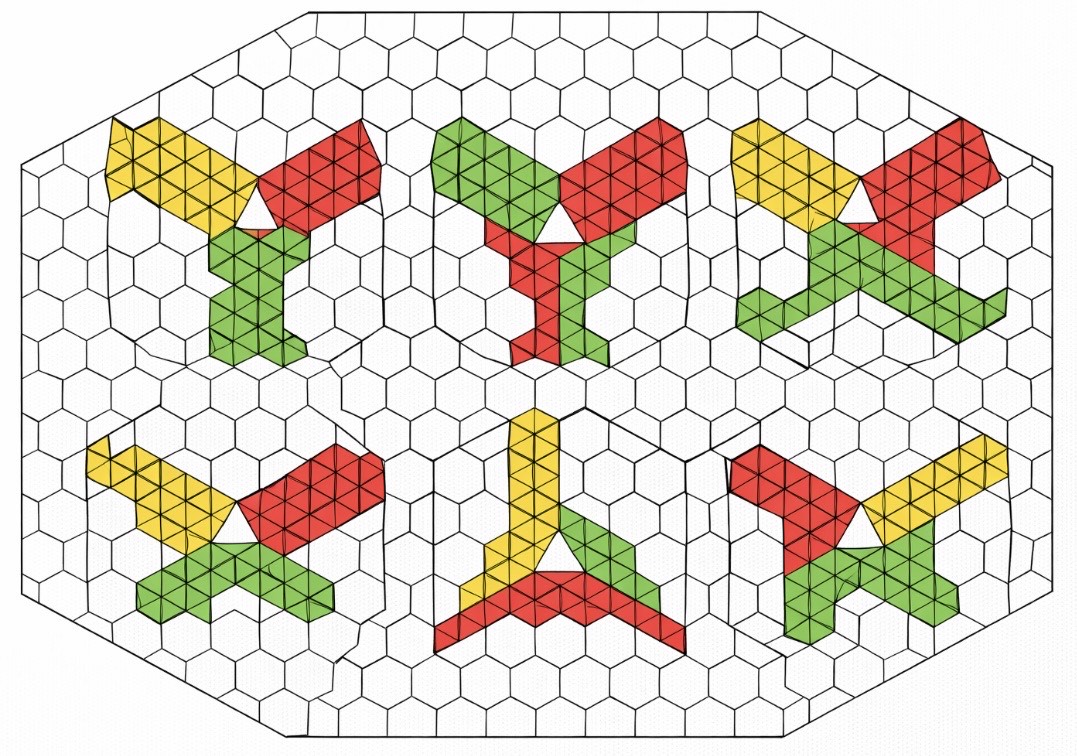}
\end{center}
\caption{A depiction of the interior of a finite volume over the hexagonal lattice for which connected components are colored in green, orange, and red.}
\end{figure}

\subsection{$\mathcal{C}_j \equiv \mathscr{C}^{\prime}_{2_j}$}

\noindent We make use of several observations previously mentioned for $\mathcal{C}_j \equiv \mathscr{C}_{2_j}$ in Section \textit{5.1}. In the fourth case, it suffices to make use of the following dual crossing events,

\begin{align*}
        \mathscr{C}^{\prime}_{2_j} =  \{ \mathcal{S}_j  \underset{*}{\overset{H_{j+\delta^{\prime}}}{\longleftrightarrow}} 2_{j + \delta^{\prime}} \} \text{ } \backslash \text{ } \mathscr{C}_{2j} \text{, } 
    \mathscr{C}^{\prime}_{ 3_j}= \{ \mathcal{S}_j \underset{*}{\overset{H_{j+\delta^{\prime}}}{\longleftrightarrow}} 3_{j + \delta^{\prime}} \} \text{ } \backslash \text{ }  \mathscr{C}_{3j} \text{,       }   \text{             } \\ \mathscr{C}^{\prime}_{5_j} =  \{ \mathcal{S}_j \underset{*}{\overset{H_{j-\delta^{\prime}}}{\longleftrightarrow}} 5_{j + \delta^{\prime}} \}   \text{ } \backslash \text{ }  \mathscr{C}_{5j}  \text{, }  \mathscr{C}^{\prime}_{6_j} =  \{ \mathcal{S}_j \underset{*}{\overset{H_{j-\delta^{\prime}}}{\longleftrightarrow}} 6_{j + \delta^{\prime}} \} \text{ } \backslash  \text{ } \mathscr{C}_{6j}  \text{, } 
\end{align*}

\noindent to $\mathscr{C}_{2_j}$, $\mathscr{C}_{3_j}$, $\mathscr{C}_{5_j}$ and $\mathscr{C}_{6_j}$, respectively. That is, it suffices to argue that there exists a symmetric domain,

\begin{align*}
    \mathrm{Sym} \big( \Omega^{\prime} \big) , 
\end{align*}

\noindent which admits the decomposition,

\begin{align*}
     \mathrm{Sym} \big( \Omega^{\prime} \big) \equiv \gamma^{\prime}_1 \cup \gamma^{\prime}_2 \cup \gamma^{\prime}_3 \cup \gamma^{\prime}_4       , 
\end{align*}

\noindent where,

\begin{align*}
  \gamma^{\prime}_1 \equiv  \textit{Left boundary of $  \mathrm{Sym} \big( \Omega^{\prime} \big)$}  , \\ \\ \gamma^{\prime}_2 \equiv   \textit{Top boundary of $  \mathrm{Sym} \big( \Omega^{\prime} \big)$}  , \\ \\ \gamma^{\prime}_3 \equiv  \textit{Bottom boundary of $  \mathrm{Sym} \big( \Omega^{\prime} \big)$}  , \\ \\ \gamma^{\prime}_4  \equiv   \textit{Right boundary of $  \mathrm{Sym} \big( \Omega^{\prime} \big)$}  , 
\end{align*}

\noindent for which,

\begin{align*}
\big\{   \mu \big[   C^{\prime}_0 \big| \Gamma_{2_j} \equiv \gamma^{\prime}_1 , \Gamma_{2_{j+2\delta^{\prime}}} = \gamma^{\prime}_2       \big] \text{ } \textit{whp} \big\} \Longrightarrow  \big\{ \mu \big[ \textit{the crossing }  C^{\prime}_0    \textit{ occurs across } \mathrm{Sym} \big( \Omega^{\prime} \big)     \big] \text{ } \textit{whp} \big\}  . 
\end{align*}

\noindent To demonstrate that the above correspondence holds for $C^{\prime}_0$ as it does for $C_0$, observe,

\begin{align*}
         \mu[      \gamma^{\prime}_1 \underset{*}{\overset{ \Omega^{\prime}  }{\longleftrightarrow}} \gamma^{\prime}_2      | \Gamma_{2_j} = \gamma^{\prime}_1 , \Gamma_{2_{j + 2 \delta^{\prime}}} = \gamma^{\prime}_2  ]      \text{ , } 
\end{align*}

\noindent can be lower bounded with the inequality,

{\small \begin{align*}
    \mu[      \gamma^{\prime}_1 \underset{*}{\overset{ \Omega^{\prime}  }{\longleftrightarrow}} \gamma_2   |  \Gamma_{2_j} = \gamma^{\prime}_1 ,  \Gamma_{2_{j + 2 \delta^{\prime}}} = \gamma^{\prime}_2 , k_{\gamma^{\prime}_1} = k_{\gamma^{\prime}_2} = 1  ]    \geq \mu_{\Omega}^{\{\gamma_1 , \gamma_2\}}[    \gamma^{\prime}_1    \underset{*}{\longleftrightarrow}   \gamma^{\prime}_2    ]  \tag{$\mathrm{V}^{*}$-$\mathrm{SYM DOM}$}  \text{. }  
\end{align*} } 

\noindent Furthermore, under the assumption that the above inequality is taken for wired boundary conditions of the high-temperature expansion of the loop $O(n)$ model $\mu \big[ \cdot \big]$, one can adapt the procedure for pushing boundary conditions onto symmetric domains for the dual crossings from that which is presented in $\textbf{Definition}$ \textit{2}.

\subsubsection{Incident layer of hexagons to the dual symmetric domain boundary}.

\noindent For quantifying the following dual crossing probability,

\begin{align*}
   \mu \big[    \mathscr{C}^{\prime}_{2_{j + \delta^{\prime}}}              \backslash            \big( C^{\prime}_0 \cup C^{\prime}_2 \big)       \big] ,
\end{align*}

\noindent we make use of the following procedure:

\bigskip

\noindent \textbf{Definition} \textit{3} (\textit{pushing boundary conditions onto symmetric domains from boundary conditions on the dual of $H_j$}) From boundary conditions along $2_j$, before reflecting connected components induced by the crossings event $\mathscr{C}_{2_j}$ about $2_j$, boundary conditions along symmetric domains are obtained with the following procedure:

\begin{itemize}
    \item[$\bullet$] To partition vertices in $\mathrm{Sym}$ for constructing boundary conditions on vertices along the boundaries of symmetric domains, we assign $+$ boundary conditions to a partition of the first layer of hexagons outside of the crossing induced by $\mathscr{C}^{\prime}_{2_{j+\delta^{\prime}}} \backslash (C^{\prime}_0 \cup C^{\prime}_2)$, conditioned under realizations of paths $\gamma^{\prime}_1 \text{ } \& \text{ } \gamma^{\prime}_2$.
    
    \item[$\bullet$] To apply $(\mathcal{S}-\mathrm{SMP})$, given the crossing $\mathscr{C}^{\prime}_{2_{j+\delta^{\prime}}} \backslash (C^{\prime}_0 \cup C^{\prime}_2)$, the length of the boundary of the symmetric domain is determined by the number of connected components of the spin configuration, which corresponds to the the edges present in the configuration. From the total number of vertices on the boundary, we introduce boundary conditions with $(\mathcal{S}-\mathrm{CBC})$. Outside of the dual of $H_j$ over $\textbf{T}$, the paths $\gamma^{\prime}_1$ and $\gamma^{\prime}_2$.

    \item[$\bullet$] After having identified the boundaries of the symmetric domain, reflection of one half of $\mathrm{Sym}$ is constructed by taking the union $\big( \gamma^{x_{\gamma_1 , \gamma_2}}_1 \big)^{\prime} \cup \big( \gamma^{x_{\gamma_1, \gamma_2}}_2 \big)^{\prime}$, where the paths in the union denote the restriction of the connected components of $\gamma^{\prime}_1$ and $\gamma^{\prime}_2$ after $\mathscr{C}^{\prime}_j$ and $\mathscr{C}^{\prime}_{j + 2 \delta^{\prime}}$ have occurred. The remaining top half of $\mathrm{Sym}$ is obtained by reflection through $2_j$ that was crossed by $\gamma^{\prime}_1$ and $\gamma^{\prime}_2$, as with the remaining half of the lower part.
    
    \item[$\bullet$] The reflections $\big( \tilde{\gamma}^{x_{\gamma_1 , \gamma_2}}_1 \big)^{\prime}$ and $\big( \tilde{\gamma}^{x_{\gamma_1 , \gamma_2}}_2 \big)^{\prime}$ described in previous steps provide the remaining half of $\mathrm{Sym}$ after performing appropriate reflections with respect to the intersection with the side of the hexagonal box.
\end{itemize}

\subsection{$(\mathcal{S}-\mathrm{SMP})$ property}

\noindent With the above procedure for pushing boundary conditions, it suffices to argue, under,

\begin{align*}
    \mu_{\Omega^{\prime}} \equiv \mu_{\Omega^{\prime}}[  \cdot |_{\Omega^{\prime}}  |      \gamma^{\prime}_1 \cap \gamma^{\prime}_2 = \emptyset , \gamma^{\prime}_1 \cap \gamma^{\prime}_3 = \emptyset  , k_{\gamma^{\prime}_1} = k_{\gamma^{\prime}_2}   = 1 ] ,
\end{align*}

\noindent that,

\begin{align*}
    \mu[  \mathscr{C}^{\prime}_{2_j} \backslash ( C^{\prime}_0   \cup   C^{\prime}_2 )    |  \Gamma_{2_j} = \gamma^{\prime}_1 ,  \Gamma_{2_{j + 2 \delta^{\prime}}} = \gamma^{\prime}_2  ] \leq \mu_{\Omega}^{\{\gamma_1 , \gamma_2\}^{c}} [    \mathscr{C}^{\prime}_{2_{j + \delta^{\prime}}  }   ]    \text{ , }
\end{align*}

\noindent as demonstrated in Section \textit{5.2.1}, where,

\begin{align*}
  C^{\prime}_0 \equiv \textit{Dual bottom to top crossing across $H_j$}  , \\ \\ C^{\prime}_2 \equiv \textit{Dual bottom to top crossing across $H_{j+2 \delta^{\prime}}$}   .
\end{align*}

\noindent To argue that the desired upper bound $\mu_{\Omega}^{\{\gamma_1 , \gamma_2\}^{c}} [    \mathscr{C}^{\prime}_{2_{j + \delta^{\prime}}  }   ]  $ holds, we make use of arguments in Section \textit{5.2.1}, and \textit{5.2.2}, for upper bounding conditional crossing events, and pushing wired boundary conditions away from $\Omega^{\prime}$ towards $\mathrm{Sym} \big( \Omega^{\prime} \big)$, respectively, below.

\subsubsection{Upper bound for conditional crossing events across dual symmetric domains and pushing wired boundary conditions away from $\Omega^{\prime}$ towards $\mathrm{Sym} \big( \Omega^{\prime} \big)$}.

\noindent From the lower bound provided in ($\mathrm{V}^{*}$-$\mathrm{SYM DOM}$), to obtain the desired upper bound,

\begin{align*}
   \mu_{\mathrm{Sym}( \Omega^{\prime} )}^{ \{L,R\}}[    L^{\prime}_{\mathrm{Sym}} \longleftrightarrow R^{\prime}_{\mathrm{Sym}}    ]  \equiv \mu_{\mathrm{Sym}( \Omega^{\prime} )}^{ \{L,R\}}[    L_{\mathrm{Sym}} \underset{*}{\longleftrightarrow }R_{\mathrm{Sym}}    ]      , 
\end{align*}

\noindent as presented in Section \textit{5.2.2} for crossing events before passing to the dual, observe,

{\small \begin{align*}
   \mu_{\Omega^{\prime}}^{\{\gamma_1 , \gamma_2 \}^c}[       \mathscr{C}^{\prime}_{2_{j + \delta^{\prime}}}     ]  \leq \mu_{\mathrm{Sym}}^{\{T,B\}}[     \mathscr{C}^{\prime}_{2_{j + \delta^{\prime}}}                ]                    \leq  \mu_{\mathrm{Sym}}^{\{T,B\}}[   T^{\prime}_{\mathrm{Sym}} \longleftrightarrow B^{\prime}_{\mathrm{Sym}}   ]  = \mu_{\mathrm{Sym} ( \Omega^{\prime}) }^{ \{L,R\}}[     L^{\prime}_{\mathrm{Sym}} \longleftrightarrow R^{\prime}_{\mathrm{Sym}}    ]   . 
\end{align*} }

\noindent By rotational invariance, the fact that,

{\small \begin{align*}
       \mu\big[  \mathscr{C}^{\prime}_{2_j} \backslash ( C^{\prime}_0   \cup   C^{\prime}_2 )   \mid \Gamma_{2_j} = \gamma^{\prime}_1, \Gamma_{2_{j + 2 \delta^{\prime}}} = \gamma^{\prime}_2 , k_{\gamma^{\prime}_1} = k_{\gamma^{\prime}_2} = 1  \big] \leq \mu_{\mathrm{Sym} ( \Omega^{\prime}) }^{ \{L,R\}}\left[    L^{\prime}_{\mathrm{Sym}} \longleftrightarrow R^{\prime}_{\mathrm{Sym}}    \right]      \text{. }   
\end{align*} }

\noindent implies,

\begin{align*}
 \big\{ \mu^{\{ \gamma^{\prime}_1 , \gamma^{\prime}_2 \}^c}_{\Omega^{\prime}} \big[ \mathscr{C}^{\prime}_{2_{j + \delta^{\prime}}  }          \big] \approx \mu^{\{ \gamma^{\prime}_1 , \gamma^{\prime}_2 \}^c}_{\mathrm{Sym} ( \Omega^{\prime} ) } \big[ \mathscr{C}^{\prime}_{2_{j + \delta^{\prime}}  }          \big]   \big\}  \Longleftrightarrow    \big\{ \mu_{\Omega}^{\{\gamma_1 , \gamma_2 \}^c}[   \mathscr{C}^{\prime}_{2_{j + \delta^{\prime}}}        ]         \approx   \mu_{\mathrm{Sym}( \Omega^{\prime} ) }^{ \{L,R\}}[    L^{\prime}_{\mathrm{Sym}} \longleftrightarrow R^{\prime}_{\mathrm{Sym}}    ]           \big\}  \\ \Longleftrightarrow     \big\{ \mu_{\Omega}^{\{\gamma_1 , \gamma_2 \}^c}[   \mathscr{C}^{\prime}_{2_{j + \delta^{\prime}}}        ]         \approx   \mu_{\mathrm{Sym} }^{ \{L,R\}}[    L_{\mathrm{Sym}} \underset{*}{\longleftrightarrow} R_{\mathrm{Sym}}    ]           \big\}     .
\end{align*}

\noindent Altogether, as previously described for $\mathcal{C}_j \equiv \mathscr{C}_{2_j}$ and $\mathcal{C}_j \equiv \mathscr{C}_{5_j}$ in Section \textit{5.2.2}, one has that,

{\small \begin{align*}
     \big\{   \mu\big[  \mathscr{C}^{\prime}_{2_j} \backslash ( C^{\prime}_0   \cup   C^{\prime}_2 )   \mid \Gamma_{2_j} = \gamma^{\prime}_1, \Gamma_{2_{j + 2 \delta^{\prime}}} = \gamma^{\prime}_2 , k_{\gamma^{\prime}_1} = k_{\gamma^{\prime}_2} = 1  \big]  \textit{ whp}                            \big\} \Longrightarrow \big\{                     \mu \big[ \big\{ \mathcal{C}_j \equiv \mathscr{C}^{\prime}_{2_j} \big\} , \big\{ L^{\prime}_{\mathrm{Sym}} \\ \longleftrightarrow R^{\prime}_{\mathrm{Sym}} \big\}  \big]     \textit{ whp}   \big\}   \Longrightarrow \big\{         \mu \big[ \big\{ \mathcal{C}_j \equiv \mathscr{C}^{\prime}_{5_j} \big\} , \big\{ L^{\prime}_{\mathrm{Sym}} \longleftrightarrow R^{\prime}_{\mathrm{Sym}} \big\}  \big]                     \textit{ whp}   \big\}   . 
       \end{align*}
       }

\subsection{$\mathcal{C}_j \equiv \mathscr{C}^{\prime}_{3_j}$}

In the fifth case, one can adapt similar arguments to dual crossings over symmetric domains as presented in the previous case for $\mathcal{C}_j \equiv \mathscr{C}^{\prime}_{2_j}$. Namely, for the crossings $\mathscr{C}^{\prime}_{3_j}$ and $\mathscr{C}^{\prime}_{3_{j+2\delta^{\prime}}}$, observe,

\subsection{$\mathcal{C}_j \equiv \mathscr{C}^{\prime}_{4_j}$}

In the sixth case, recall that for $\mathcal{C}_j \equiv \mathscr{C}_{4_j}$ one analyzed the conditionally defined crossing probability,

\begin{align*}
          \mu[  \mathscr{C}_{4_j} \backslash ( C_0   \cup   C_2 )   |  \mathcal{R}  ]                 \text{, } 
\end{align*}

\noindent where the conditioning is,

\begin{align*} \mathcal{R} \equiv \bigg\{ \big\{ \Gamma_{2_j} = \gamma_1 \big\} , \big\{ \Gamma_{2_{j + 2 \delta^{\prime}}} = \gamma_2 \big\} ,  \big\{ \gamma_1 \cap \gamma_3 = \emptyset \big\} ,  \big\{ \gamma_1 \cap \gamma_2 = \emptyset  \big\} \bigg\} . 
\end{align*}

\noindent Along these lines, under the same choice of dual top to bottom crossings $C^{\prime}_0$ and $C^{\prime}_2$ introduced in Section \textit{5.6}, under the same choices of paths $\Gamma_{2_j} = \gamma^{\prime}_1$ and $\Gamma_{2_{j+2 \delta^{\prime}}} \equiv \gamma^{\prime}_2$, observe, given the conditioning,

\begin{align*} \mathcal{R}^{\prime} \equiv \bigg\{ \big\{ \Gamma_{2_j} = \gamma^{\prime}_1 \big\} , \big\{ \Gamma_{2_{j + 2 \delta^{\prime}}} = \gamma^{\prime}_2 \big\} ,  \big\{ \gamma^{\prime}_1 \cap \gamma^{\prime}_3 = \emptyset \big\} ,  \big\{ \gamma^{\prime}_1 \cap \gamma^{\prime}_2 = \emptyset  \big\} \bigg\} . 
\end{align*}

\noindent that it suffices to express that the conditionally defined dual crossing event,

\begin{align*}
          \mu[  \mathscr{C}^{\prime}_{4_j} \backslash ( C^{\prime}_0   \cup   C^{\prime}_2 )   |  \mathcal{R}^{\prime}  ]                           \text{, } 
\end{align*}

\noindent implies that crossings across boundaries of $\mathrm{Sym} \big( \Omega^{\prime} \big)$ occur with high probability.

\subsubsection{Pushing boundary conditions way from the dual of $H_j$ towards $\mathrm{Sym} \big( \Omega^{\prime} \big) $}

\noindent To push boundary conditions away from the dual of $H_j$, we adopt a closely related set of assumptions for the symmetric domains considered under original crossing events defined in \textit{5.4.1}. That is, under the realization of paths,

\begin{align*} \big\{ \Gamma_{2_j} = \gamma^{\prime}_1 \big\} \\ \\  \big\{ \Gamma_{2_{j + 2 \delta^{\prime}}} = \gamma^{\prime}_2 \big\} , \end{align*}

\noindent observe, 

\begin{align*}
     \mu \big[ C^{\prime}_0 \big| \mathcal{R}^{\prime} \big] \geq  \mu \big[   \gamma^{\prime}_1  \overset{\Omega^{\prime}}{\longleftrightarrow} \gamma^{\prime}_2 \big|   \mathcal{R}^{\prime}    \big]                    , 
\end{align*}

\noindent holds from arguments applied when $\mathcal{C}_j \equiv \mathscr{C}^{\prime}_{2_j}$. Moreover, by construction, $\big( H_j \big)^{*} \subsetneq \mathrm{Sym} \big( \Omega^{\prime} \big)$ implies,

\begin{align*}
       \mu[      \mathscr{C}^{\prime}_{4_j}         |    \mathcal{R}^{\prime}    ]      \leq \mu_{\Omega^{\prime}_{(L^{\prime} \cup R^{\prime})^c}}^{\{\gamma^{\prime}_1 , \gamma^{\prime}_2 \}^c}[    \mathscr{C}^{\prime}_{4_j}         |     \mathcal{R}^{\prime}   ]      \text{. } 
\end{align*}

\noindent Moreover, given the realizations of the paths $\gamma^{\prime}_1$ and $\gamma^{\prime}_2$, as demonstrated in Section \textit{5.4.1}, observe, 

{\small \begin{align*}
            \mu_{\Omega^{\prime}_{(L^{\prime} \cup R^{\prime})^c}}^{\{\gamma^{\prime}_1 , \gamma^{\prime}_2\}^c}[       C^{\prime}_0    |  \mathcal{R}^{\prime}    ]     \leq        \mu_{\mathrm{Sym} ( \Omega^{\prime} ) }^{(\mathrm{Top \text{ } Half})}[   C^{\prime}_0    |  \mathcal{R}^{\prime}    ]    \leq   \mu_{\mathrm{Sym} ( \Omega^{\prime} ) }^{(\mathrm{Top \text{ } Half})} [    T^{\prime}_{\mathrm{Sym}} \overset{*}{\longleftrightarrow} B^{\prime}_{\mathrm{Sym}}  ]  = \mu_{\mathrm{Sym} ( \Omega^{\prime} ) }^{(\mathrm{Top \text{ } Half})^{\frac{2 \pi}{3}}}[  L^{\prime}_{\mathrm{Sym}} \overset{*}{\longleftrightarrow} \\ R^{\prime}_{\mathrm{Sym}}     ]   \text{, }
\end{align*} }

\noindent where $(\mathrm{Top \text{ } Half})$ denotes wired boundary conditions along the top half of hexagonal $\mathrm{Sym} \big( \Omega^{\prime} \big)$. Straightforwardly, by adapting the remarks for $\mu_{\mathrm{Sym}}^{\{L,R \}}[  T^{\prime}_{\mathrm{Sym}} \overset{*}{\longleftrightarrow} B^{\prime}_{\mathrm{Sym}}      ]$  in place of $\mu_{\mathrm{Sym}}^{\{L,R \}}[     T_{\mathrm{Sym}} \longleftrightarrow B_{\mathrm{Sym}}       ]$, one has that,

\begin{align*}
        \mu[  \mathscr{C}^{\prime}_{4_j} \backslash ( C^{\prime}_0   \cup   C^{\prime}_2 )    |\mathcal{R}^{\prime}   ]     \leq   \mu_{\mathrm{Sym} ( \Omega^{\prime} ) }^{(\mathrm{Top \text{ } Half)}^{\frac{2\pi}{3}}}[ L^{\prime}_{\mathrm{Sym}} \overset{*}{\longleftrightarrow} R^{\prime}_{\mathrm{Sym}}      ]     \text{. }
\end{align*}

\subsubsection{$(\mathcal{S} -\mathrm{CBC})$ lower bound for the conditional crossing event $\{  C^{\prime}_0   |  \mathcal{C}^{\prime}_0 \cap \mathcal{C}^{\prime}_4 \}$}.

\noindent We complete the arguments for the dual crossing events through the following inequalities. Similar to the inequality provided for the occurrence of $C_0$ conditionally upon the intersection of crossing events $C_0 \cap C_4$ in Section \textit{5.4.2}, it suffices to argue that an inequality of the form,

{\small \begin{align*}
      \mu [   C^{\prime}_0    |  \mathcal{C}^{\prime}_0  \cap  \mathcal{C}^{\prime}_4    ]     \geq         \frac{1}{n ^{\big(   k^{\tau_{H_0 \cup H_{2 \delta} } } (\sigma) \big)^{\prime}   - \big(     k^{\tau_{H_1 \cup H_{1+\delta}}  } (\sigma) \big)^{\prime} } x \text{ } \mathrm{exp}(h) \text{ }  }      \text{ }  \mu  [      \mathcal{C}^{\prime}_2 \backslash (  C^{\prime}_0 \cup C^{\prime}_2    )   |        \mathcal{C}^{\prime}_0 \cap \mathcal{C}^{\prime}_4      ] \tag{($\mathcal{S}$-SMP)- ($\mathcal{C}^{\prime}_2$)}  \text{, } 
\end{align*} }

\noindent holds.

\bigskip

\noindent Fix some $I^{\prime} \neq I$ taken to be sufficiently large. For the above inequality for dual crossing events, denote,

{\small \[
\left\{\!\begin{array}{ll@{}>{{}}l}
\big(   k^{\tau_{H_0 \cup H_{2 \delta} } } (\sigma) \big)^{\prime}  \equiv \textit{Number of connected components $k$ associated with boundary conditions $\tau$ sup-} \\ \textit{ported over the dual of $H_0 \cup H_{2 \delta}(\sigma) $ for the spin configuration $\sigma$}   , \\ \\    \big(     k^{\tau_{H_1 \cup H_{1+\delta}}  } (\sigma) \big)^{\prime} \equiv \textit{Number of connected components $k$ associated with boundary conditions $\tau$ sup-} \\ \textit{ported over the dual of $H_1 \cup H_{1+ \delta}(\sigma) $ for the spin configuration $\sigma$}   , \\  \\   \big(  e^{\tau_{H_0 \cup H_{2 \delta}}}(\sigma) \big)^{\prime}  \equiv      \textit{Edges $e$ associated with boundary conditions $\tau$ supported over the dual of $H_0 \cup H_{2 \delta}$ }  \\ \textit{for the  spin configuration $\sigma$}     , \\  \\ \big(  e^{\tau_{H_1 \cup H_{1+\delta}}}(\sigma) \big)^{\prime}  \equiv   \textit{Edges $e$ associated with boundary conditions $\tau$ supported over the dual of $H_1 \cup H_{1+\delta}$ }  \\ \textit{for the spin configuration $\sigma$}          , \\ \\ \big( r_{H_0 \cup H_{\delta} }(\sigma_1) \big)^{\prime} \equiv     \textit{Summation of the spins of the face configuration $\sigma_1$ over the dual of } \\ \textit{$H_0 \cup H_{\delta}$}    , \\ \\       \big(  r_{H_1 \cup H_{1+\delta}}(\sigma_2) \big)^{\prime}  \equiv \textit{Summation of the spins of the face configuration $\sigma_2$ over the dual} \\ \textit{ of $H_1 \cup H_{1+\delta}$}    , \\ \\ \big( r^{\prime}_{H_0 \cup H_{\delta}}(\sigma_1) \big)^{\prime} \equiv \textit{Summation over the number of monochromatically colored triangles of $\sigma_1$} \\ \textit{over the dual of $H_0 \cup H_{\delta}$} , \\ \\ \big( r^{\prime}_{H_1 \cup H_{1+\delta}}(\sigma_2) \big)^{\prime} \equiv \textit{Summation over the number of monochromatically colored triangles of $\sigma_2$} \\ \textit{over the dual of $H_1 \cup H_{1+\delta}$}  .  
\end{array}\right.
\] }

\bigskip

\noindent from which one can obtain the desired lower bound,

{\small \begin{align*}
   c^{\prime} \equiv \big( 3 I^{\prime} \big)^{-3} \bigg\{     n^{ \big(  \big( k^{\tau_{H_0 \cup H_{2 \delta}  } } (\sigma) \big)^{\prime}  \cap  \big( k^{\tau_{H_1 \cup H_{1+\delta}} }(\sigma)    \big) \big)^{\prime} }   x^{  \big(   \big( e^{\tau_{H_0 \cup H_{2 \delta}}}(\sigma) \big)^{\prime} \cap  \big( e^{\tau_{H_1 \cup H_{1+\delta}}}(\sigma) \big)^{\prime}  \big)  } \\ \times   e^{ h \big(  \big(  r_{H_0 \cup H_{\delta} }(\sigma_1) \big)^{\prime} - \big( r_{H_1 \cup H_1+\delta}(\sigma_2) \big)^{\prime}     \big) + \frac{h^{\prime}}{2} \big(   \big( r^{\prime}_{H_0 \cup H_{\delta}}(\sigma_1) \big)^{\prime} - \big( r^{\prime}_{H_1 \cup H_{1+\delta}}(\sigma_2) \big)^{\prime} \big)}   \bigg\} \\ \\ \equiv \big( 3 I^{\prime} \big)^{-3} \mathcal{P}^{\prime} , 
\end{align*}}

\noindent as appearing in the linear combination of dual crossing probabilities,

\begin{align*}
    c^{\prime} \mu \big[ C^{\prime}_0 \big] + \mu \big[  C^{\prime}_2 \big] . 
\end{align*}

\begin{figure}[H]
\begin{center}
\includegraphics[width=0.9\columnwidth]{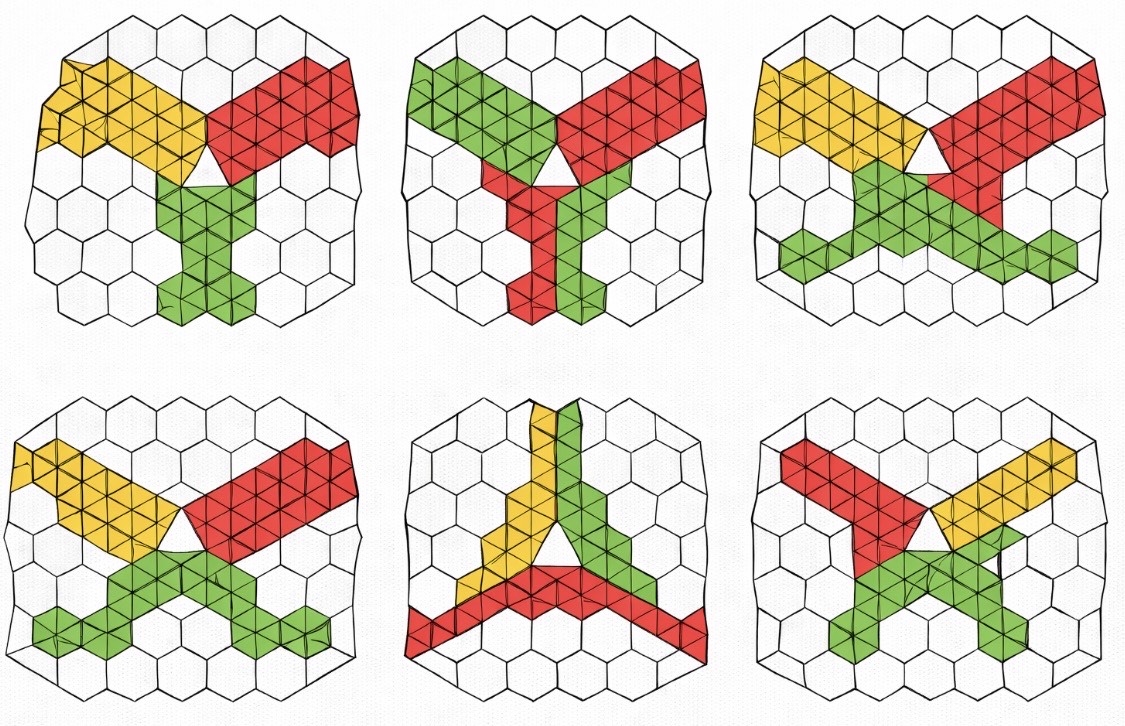}
\end{center}
\caption{A depiction of connected components, colored in orange, red and green, strictly contained within the interior of six finite volumes over the hexagonal lattice.}
\end{figure}

\noindent The desired inequality for dual crossing events is readily obtained from the following computations,

{\small \begin{align*}
   \mu \big[        C^{\prime}_0      \big|       \big\{ \Gamma_{2_j} = \gamma^{\prime}_1 \big\} ,       \big\{ \Gamma_{2_{j + 2 \delta^{\prime}}} = \gamma^{\prime}_2 \big\}                      \big]   \geq   f \big[        \mu_{\Omega}      \big[        \gamma^{\prime}_1 \overset{*}{\longleftrightarrow} \gamma^{\prime}_2 \big|    \big\{ \Gamma_{2_j} = \gamma^{\prime}_1 \big\} ,       \big\{ \Gamma_{2_{j + 2 \delta^{\prime}}} = \gamma^{\prime}_2 \big\}        \big]                \big]  \\ \geq           \mu^{\mathrm{Mix}}_{\Omega} \big[ \big\{ \textit{left boundary of $\Omega$} \big\} \overset{*}{\longleftrightarrow } \big\{ \textit{right boundary of $\Omega$} \big\}   \big]   \\ \geq   \mu^{\mathrm{Mix}}_{\Omega} \big[ \big\{ \textit{left boundary of $\mathrm{Sym} \big( \Omega \big) $} \big\} \overset{*}{\longleftrightarrow } \big\{ \textit{right boundary of $\mathrm{Sym} \big( \Omega \big) $} \big\}   \big]   .
\end{align*}}

\noindent Moreover, in light of the above lower bound between the left and right boundaries of $\mathrm{Sym} \big( \Omega \big)$, one can conclude that the dual crossing version of the conditionally defined probability,

\begin{align*}
 \mu \big[ C^{\prime}_2 \backslash \big( C^{\prime}_0 \cup C^{\prime}_4 \big) \big|   \big\{ \Gamma_{2_j} = \gamma^{\prime}_1 \big\} ,       \big\{ \Gamma_{2_{j + 2 \delta^{\prime}}} = \gamma^{\prime}_2 \big\}   ,   \pi_L , \pi_R     \big]   , 
\end{align*}

\noindent where,

\begin{align*}
  \pi_L \equiv     \textit{A random primal leftmost path so that $\gamma^{\prime}_1 \neq \gamma^{\prime}_2 $}               , \\ \\ \pi_R \equiv  \textit{A random primal rightmost path so that $\gamma^{\prime}_1 \neq \gamma^{\prime}_2 $}    , 
\end{align*}

\noindent can be upper bounded as follows,

{\small \begin{align*}
  \mu \big[ C^{\prime}_2 \backslash \big( C^{\prime}_0 \cup C^{\prime}_4 \big) \big|   \big\{ \Gamma_{2_j} = \gamma^{\prime}_1 \big\} ,       \big\{ \Gamma_{2_{j + 2 \delta^{\prime}}} = \gamma^{\prime}_2 \big\}   ,   \pi_L , \pi_R     \big]   \leq     f \big[         \mu \big[       \big\{ \textit{top boundary of $\Omega^{\prime}$}  \big\} \overset{*}{\longleftrightarrow }  \big\{ \textit{bot-} \\ \textit{tom boundary of $\Omega^{\prime}$}  \big\}         \big|   \big\{ \Gamma_{2_j} = \gamma^{\prime}_1 \big\}  ,       \big\{ \Gamma_{2_{j + 2 \delta^{\prime}}} = \gamma^{\prime}_2 \big\}   ,   \pi_L , \pi_R     \big]        \big]     \\     \leq    \mu^{(\mathrm{Mix} \big)^{\prime}}_{\Omega^{\prime}} \big[    \big\{ \textit{top boundary of $\Omega^{\prime}$}  \big\} \overset{*}{\longleftrightarrow }  \big\{ \textit{bottom boundary of $\Omega^{\prime}$}  \big\}         \big]    \\ \leq   \mu^{(\mathrm{Mix} \big)^{\prime}}_{\mathrm{Sym} ( \Omega^{\prime})} \big[    \big\{ \textit{top boundary of $\mathrm{Sym} \big( \Omega^{\prime} \big) $}  \big\} \overset{*}{\longleftrightarrow }  \big\{ \textit{bottom boundary of $\mathrm{Sym} \big( \Omega^{\prime} \big) $}  \big\}         \big]    \\ \leq   \mu^{(\mathrm{Mix} \big)}_{\mathrm{Sym} ( \Omega^{\prime})} \big[    \big\{ \textit{top boundary of $\mathrm{Sym} \big( \Omega^{\prime} \big) $}  \big\} \overset{*}{\longleftrightarrow }  \big\{ \textit{bottom boundary of $\mathrm{Sym} \big( \Omega^{\prime} \big) $}  \big\}         \big]           . 
\end{align*} }

\noindent Altogether, for,

\begin{align*}
  \mathcal{V}^{\prime}_{\mathcal{H}} \equiv \textit{Vertical dual crossing event across $H$}  ,
\end{align*}

\noindent as demonstrated in Section \textit{5.4.2}, for dual crossing events,

\begin{align*}
      \mathrm{max} \big\{  \mu[     \mathscr{C}^{\prime}_{2_j}    ]  , \mu[     \mathscr{C}^{\prime}_{3_j}    ]  , \mu[     \mathscr{C}^{\prime}_{4_j}   ]    : 0 \leq j < I^{\prime}      \big\}   \geq \frac{\mu[    \mathcal{V}^{\prime}_{\mathcal{H}}   ]}{3I}  \\ \Updownarrow \\  \mathrm{max} \big\{ \mu[ \mathscr{C}^{\prime}_{2_0} ] , \mu[   \mathscr{C}^{\prime}_{3_0} ] , \mu[  \mathscr{C}^{\prime}_{4_0}   ]          \big\} \geq     \frac{\mu [ \mathcal{V}_{\mathcal{H}} ]}{3I}  \\ \Updownarrow  \\   \mu[    C^{\prime}_0       ]     \geq   \mu [  \big( \mathscr{C}^1    \big)^{\prime}      \cap \big( \mathscr{C}^2   \big)^{\prime}  ]    \geq   \frac{\mu[ \mathcal{V}^{\prime}_{\mathcal{H}}  ]^2}{9I^2}             . 
\end{align*}

\noindent Hence, for both deterministic and stochastic explorations of symmetric domains, one has that,

{\small \begin{align*}
 \mu \big[ C_0 \big|    \big\{ \Gamma_{2_j} = \gamma_1 \big\} ,       \big\{ \Gamma_{2_{j + 2 \delta^{\prime}}} = \gamma_2 \big\}   ,   \pi_L , \pi_R        \big] \geq  \mu \big[    C_2 \backslash \big( C_0 \cup C_4 \big)        \big|    \big\{ \Gamma_{2_j} = \gamma_1 \big\} ,       \big\{ \Gamma_{2_{j + 2 \delta^{\prime}}} = \gamma_2 \big\}   ,   \pi_L , \pi_R         \big] , \\ \\ \mu \big[ C^{\prime}_0 \big|      \big\{ \Gamma_{2_j} = \gamma^{\prime}_1 \big\} ,       \big\{ \Gamma_{2_{j + 2 \delta^{\prime}}} = \gamma^{\prime}_2 \big\}   ,   \pi_L , \pi_R        \big]  \geq  \mu \big[ C^{\prime}_2 \backslash \big( C^{\prime}_0 \cup C^{\prime}_4 \big)    \big|    \big\{ \Gamma_{2_j} = \gamma_1 \big\} ,       \big\{ \Gamma_{2_{j + 2 \delta^{\prime}}} = \gamma_2 \big\}   ,   \pi_L , \pi_R         \big]  , 
\end{align*} }

\noindent for which $C_0$ and $C^{\prime}_0$, respectively, occur from which we conclude the argument. \boxed{}

\bigskip

\noindent \textbf{Remark}. The above constant $c$ for which \textbf{Lemma} $9$ is proved, unlike the accompanying constant for the constant provided for the Random-Cluster model, is dependent on the product of the number of loops, the number of edges, and the difference in the number of monochromatically colored hexagons, as defined with the quantities $r(\sigma)$ and $r^{\prime}(\sigma)$ in ($\mathcal{S}$-SMP), instead of solely on $q$.

\section{Volume of connected components from wired boundary conditions}

\subsection{Proof of Lemma $10$ with the $\mu$ homomorphism of Proposition $8$}

To study behavior of the dilute Potts model in the \textit{{Continuous Critical}} and \textit{{Discontinuous Critical}} cases, we turn to studying vertical crossings under $\mu$ under wired boundary conditions. To denote vertical translates of hexagons containing $H_j$, we introduce $H_{j + \delta}$ as the hexagonal box which has nonempty intersection with $H_j$, and is of side length $j + \delta$. We state the following Lemma and Corollary.

\bigskip

\noindent \textbf{Lemma $10$} (\textit{volume of connected components}): For $x \in H_j$ and $C \geq 2$, there exists $\epsilon > 0 $ such that, given $\mu^1_{H_{C j}}[   H_j \longleftrightarrow \partial \text{} H_{j +\delta}  \text{ } ] < \epsilon$ for some $k$, in $H_j \cap H_{ j + \delta}$ there exits a positive $c$ satisfying,

\begin{align*}
             \mu^1_{H_j}[   \mathrm{Vol}\big(    \text{connected components in the annulus } H_j \cap H_{j +\delta}   \big) = N  ]   \leq e^{-cN}  \text{, } 
\end{align*}

\noindent for every $j , N \geq 2$, taken under wired boundary conditions.

\bigskip

\noindent \textit{Proof of Lemma $10$}. The arguments require use of hexagonal annuli which for simplicity we denote with $\mathcal{H}_{\mathcal{A}} \equiv H^c_j \cap H_{j+\delta}$, in which one hexagonal box is embedded within another, and set $\mathcal{P} \equiv \{ \mathrm{Vol}\big(   \text{connected components in the annulus } H^c_j \cap H_{j +\delta}   \big) = N \}$. The existence of the quantity $\mu^{\mathscr{C}_l}$, where $\mu$ is a finite constant and $\mathscr{C}_l$ is the number of connected components of length $l$ is standard from {[29]}. To prove the statement, we measure the connected components of length $l$ from the center of $H_j$ in $\mathcal{H}_{\mathcal{A}}$. 

From the connected components of $x$ in $H_j$, we can restrict the connected components to the nonempty intersection given by $\mathcal{H}_{\mathcal{A}}$. The argument directly transfers from the planar case to the hexagonal one with little modification, as the restriction of the connected components $\mathscr{C}_l$ of length $l$ to the annulus implies the existence of a connected set of in $\textbf{H}$, denoted with $S \subsetneq  \mathcal{H}_{\mathcal{A}}$ of vertex cardinality $N \backslash |H_j|$ from which a subset of the connected components $\mathcal{S}_{\mathcal{C}} \subsetneq  S$ can be obtained. We conclude the proof by analyzing the pushforward of $\mathcal{P}$ under wired boundary conditions supported on $H_j$, in which the union bound below over $\mathcal{J}_S$ satisfies,

{\small \begin{align*}
    \mu^1_{H_j}[     \mathcal{P}       ]  \leq \bigcup_{i \text{ }  \in \text{ }  \mathcal{J}_S =  \{ \text{connected components of size } l \text{ of } \mathcal{H}_{\mathcal{A}} \text{ in } S \} } \text{ } \mu^1_{H_j}[   \mathcal{P}_{i}    ]   \leq \big[ \mu^1_{H_j}[\mathcal{P}] \big]^{|N \backslash |H_j|| } \\ \leq \big[ \mu \epsilon^{|\mathcal{S}_{\mathcal{C}}|} \big]^{| N \backslash |H_j| |}  \leq \text{           }     e^{-c N} \text{, }  
\end{align*}     }     

\noindent where the union is taken over the collection of connected components under the criteria that admissible vertices from $S$ are taken to be of distance $2j$ from one another in $\mathcal{J}_S$, and events $\mathcal{P}_{\mathcal{J}_S}$ denote measurable events under $\mu^1_{H_j}$ indexed by the number of admissible vertices from $\mathcal{S}_{\mathcal{C}}$. We also apply $(\mathcal{S}-\mathrm{SMP})$ and $(\mathcal{S}-\mathrm{CBC})$ in the inequality above to push boundary conditions away, with $\epsilon$ arbitrary and small enough. \boxed{ }

\bigskip

We turn to the statement of \textbf{Corollary} $11$ below which requires modification to vertical crossings across $H_j$, which can be accommodated with families of boxes $H_j$ with varying height dependent on the usual RSW aspect ratio factor $\rho$. We also make use of $\mathcal{S}_{T,L} \equiv \mathcal{S}$.

\bigskip

\noindent \textbf{Corollary $11$} (\textit{dilute Potts behavior outside of the supercritical and subcritical regimes}): For every $\rho > 0$, $L \geq 1$, there exists a positive constant $\mathcal{C}$ satisfying the following, in which 

\begin{itemize}

    \item[$\bullet$] for the \textbf{Non}(\textit{{Subcritical}}) regime, the crossing probability under wired boundary conditions of a horizontal crossing across $\mathcal{H}_j$ supported over the strip, $\mu^1_{\mathcal{S}}[ \mathcal{H}_{\mathcal{H}_j}  ] \geq \mathcal{C}$, 
    
    \item[$\bullet$] for the \textbf{Non}(\textit{{Supercritical}}) regime, the crossing probability under free boundary conditions of a vertical crossing across $H_j$, $\mu^0_{\mathcal{S}}[ \mathcal{V}_{\mathcal{H}_j} ] \leq 1 - \mathcal{C} $, also supported over the strip.
    
\end{itemize}

\bigskip

\noindent \textit{Proof of Corollary $11$}. We present the argument for the first statement in \textbf{Non}(\textit{{Subcritical}}) from which the second statement in \textbf{Non}(\textit{{Supercritical}}) follows. For $\mathcal{S}$, in the \textbf{Non}(\textit{{Subcritical}}) phase horizontal crossing probabilities across $\mathcal{S}_{T,L} \equiv \mathcal{S}$ are bound uniformly away from $0$, which for $\mu$ can be demonstrated through examination of crossing events $C_j$ first introduced in \textit{Proof} of \textbf{Proposition} $8$. For , the result under which the pushforward with wired boundary conditions takes the form, for any $j \geq 1$,

\begin{align*}
   \mu^1_{\mathcal{S}}[  C_j  ]       \geq    e^{-6c}   \text{, } 
\end{align*}

\noindent from an application of \textbf{Lemma} $10$ to a connected component with unit volume in $\mathcal{H}$ type annuli. 

Also, in the following arrangement, we introduce a factor $\rho$ for the aspect length of a regular hexagon in $\mathcal{S}_{T,L}$ which mirrors the role of $\rho$ in RSW theory for crossings across rectangles. About the origin, we pushforward vertical crossing events on each side of $\mathcal{H}_j = \cup_i H_{j + \delta_i}$, respectively given by $H_{j + \delta_k }$ and $H_{j + \delta_l }$ for $k$ such that $H_{j + \delta_k}$ and $H_{ j + \delta_l}$ are of equal distance to the left and right of the origin. By construction, in any $\mathcal{H}_j$ with the aspect length dependent on $\rho$, intermediate regular hexagons can be embedded within $\mathcal{H}_j$ corresponding to the partition of the aspect length $\rho$. Longer horizontal or vertical crossings can be constructed through applications of $(\mathrm{FKG})$ which are detailed below.

From the lower bound on the volume of a unit connected component, a vertical crossing across a hexagon of aspect height $\delta$, from reasoning as given in $(\mathrm{FKG})$ can be bound below by the product of crossing probabilities of $\delta_i$ translates of vertical crossings across hexagons of aspect height $\delta_i$.

The measure under wired boundary conditions, for a vertical crossing $\mathcal{V}$ across $H_{j + \delta_k}$, is

\begin{align*}
      \mu_{\mathcal{H}_j} [  \mathcal{V}_{H_{j+ \delta_k}}       ]               \text{, }
\end{align*}

\noindent where the measure for the vertical crossing event given above is supported over $\mathcal{H}_j$.

  From the upper bound of $(\mathrm{FKG})$, longer vertical horizontal crossings occur across $2^i$ vertical translates of shorter vertical crossings. The next ingredient includes making use of previous arrangements of horizontal translates of $H_j$, namely the left translate $H_{j - \delta^{\prime}}$ and the right translate $H_{j + \delta^{\prime}}$. Under the occurrence of vertical crossings across $H_{j + \delta_k}$ and $H_{j + \delta_l}$. From this event, to show that some box $H_j$ in between $H_{j + \delta_k}$ and $H_{j + \delta_l}$ is crossed vertically, under wired boundary conditions supported over $H_j$ we directly apply previous arguments from $(\mathrm{FKG})$, with the exception that $(\mathrm{FKG})$ is applied to a countable intersection of vertical, instead of horizontal, crossing events $\mathcal{V}$. 
  
  Conditionally, if vertical crossings in $H_{j + \delta_k}$ and $H_{j + \delta_l}$ occur about arbitrary $H_{j + \delta_i}$ with $k \leq i \leq l$, then the probability below satisfies, under wired boundary conditions,

\begin{align*}
    \mu^1[ \mathcal{V}_{H_{ j + \delta_k}}     \cap      \mathcal{V}_{H_{j + \delta_l}}      ]     \geq \mu^1_{\mathcal{H}_j}[      \mathcal{V}_{H_{ j + \delta_k}}         ]  \mu^1_{\mathcal{H}_j}[ \mathcal{V}_{H_{j + \delta_l}}     ] = \mu^1_{\mathcal{H}_j}[    \mathcal{V}_{H_{j + \delta_l}}   ]^2     \geq \overset{j}{\underset{i=1}{\prod}} \mu^1_{\mathcal{H}_j}[  \mathcal{V}_{H_{j + \delta_{l_i}}} ]  \\ = \mu^1_{\mathcal{H}_j}[                               \mathcal{V}_{H_{j + \delta_{l_i}}}   ]^{2^{1-i}}  \tag{$ \circ$}   \text{, } 
\end{align*}

\noindent where $\mathcal{V}_H$ denotes the vertical crossing across hexagons of aspect length which is the same as that of $H_{j + \delta_k}$, but with and aspect height $\delta_{l_i}$ where $\delta_l = \cup_i \delta_{l_i}$. The union over $i$ indicates a partition of the aspect height of $H_{j + \delta_l}$ into $2^{1-i}$ intervals. Finally,

\begin{align*}
  \mu^1_{\mathcal{H}_j}\text{ } \big[                               \mathcal{V}_{H_{j + \delta_{l_i}}} \big]^{2^{1-i}}          \geq e^{-c({2^{1-i})}} \tag{$\circ \circ$}  \text{. } 
\end{align*}

The lower bound for the inequality above is obtained from an application of \textbf{Lemma} $10^{*}$ to the volume of a connected component from vertical crossings in $H_{j + \delta_k}$ and $H_{j + \delta_l}$. Between the second and third terms in $\circ$, monotonicity in the domain allows for a comparison between the measure under wired boundary conditions respectively supported over $H_{j + \delta_i}$ and $\mathcal{H}_j$.

From the partition of $\mathcal{H}_j$, to apply $(\mathcal{S}-  \mathrm{CBC})$ we consider the region between vertical crossings across $\mathcal{H}_{j + \delta_l}$ and $\mathcal{H}_{j + \delta_k}$. From the previous upper bound, given some $u$ the vertical event$\{\mathcal{V}_{\mathcal{H}_{ j + \delta_k}} \cup \mathcal{V}_{H_{j + \delta_l}} \}$ about $H_{j + \delta_u}$ occurs for some $k , l < u$. Under wired boundary conditions, the conditional vertical crossing

\begin{align*}
        \mu_{\mathcal{H}_j}^1\big[ \big\{  \mathcal{V}_{H_{j + \delta_{k-1}}} \cup \mathcal{V}_{H_{j + \delta_{l-1}}}  \big\} \big|     \big\{ \mathcal{V}_{H_{j + \delta_k}} \cup \mathcal{V}_{H_{j + \delta_l}} \big\}  \big]      \text{, } 
\end{align*}

\noindent is bounded from below by the lower bound of $(\circ \text{ }  \circ)$. With conditioning on $\{ \mathcal{V}_{H_{j + \delta_{k}}} \cup \mathcal{V}_{H_{j + \delta_{l}}}    \}$, the probability of simultaneous vertical crossings in $H_{j + \delta_k}$ and $H_{j + \delta_l}$ and $j + \delta_k \equiv j + \delta_l$, the pushforward under wired boundary conditions of vertical crossings across two hexagons which entirely overlap with one another gives the lower bound,

{\small \begin{align*}
        \mu^1_{\mathcal{H}_j}  \big[  \mathcal{V}_{\textbf{1}_{\{j + \delta_k \equiv j + \delta_l\}}}     \big]      \geq \mu^1_{\mathcal{H}_j}\big[   \mathcal{V}_{H_{j + \delta_k}} \cup \mathcal{V}_{H_{j + \delta_l}}    \big]  \bigg\{ \prod_{i=1}^j   \mu^1_{\mathcal{H}_j}\big[ \big\{  \mathcal{V}_{H_{j + \delta_{k-1}}} \cup \mathcal{V}_{H_{j + \delta_{l-1}}} \big\}   \big|  \big\{ \mathcal{V}_{H_{j + \delta_{k}}} \\ \cup \mathcal{V}_{H_{j + \delta_{l}}} \big\}   \big]  \bigg\}   \geq e^{-c}    \text{, } 
\end{align*} }

\noindent where the vertical crossing $\mathcal{V}$ occurs when the indicator is satisfied. As $\rho \longrightarrow + \infty$, the finite volume measure over $\mathcal{H}_j$ under the weak limit of measures yields a similar inequality

\begin{align*}
    \mu_{\mathcal{S}}^1[    \mathcal{V}_{\textbf{1}_{\{j + \delta_k \equiv j + \delta_l\}}}  ]    \geq \mu^1_{\mathcal{S}}[  \mathscr{C}_0     ]   \geq e^{-c} \text{, }
\end{align*}

\noindent with the exception that $\mu$ under wired boundary conditions is supported along the strip $\mathcal{S}$, and $\mathscr{C}_0$ denotes the crossing event in which hexagons to the right and left of $\mathcal{H}_0$ are crossed vertically. The exponential bound itself can be bounded below with the desired constant,

\begin{align*}
    e^{-c} \geq \mathcal{C} \text{, }
\end{align*}

\noindent establishing the inequality for the Spin measure under wired boundary conditions. From the union of vertical crossings $\mathcal{V}_{H_{j + \delta_k}} \cup \mathcal{V}_{H_{j + \delta_l}}$, applying the $\mu$ homeomorphism under the conditions on $c_0$ in $\textbf{Theorem $1$}$, 

\begin{align*}
    f \big[ x \big]  = 1   -  (\mathrm{Stretch} \big)^{-1}   c_0^{-c_0} + c_0^{-c_0} x  \text{, } 
\end{align*}

\noindent for $x = \mu^1[\mathcal{V}]$ to the inequality for vertical crossings bounded below by $\mathcal{C}$ implies that the upper bound of $\mathcal{C}$ on can be translated into a corresponding upper bound dependent on $\mathcal{C}$ for horizontal crossings, obtaining a similar upper bound under free boundary conditions, 

\begin{align*}
        \mu^0_{\mathcal{S}}[ \mathcal{V}_{\mathcal{H}_j}        ] \leq 1 - \mathcal{C}     \text{, } 
\end{align*}

\noindent concluding the argument after having taken the infinite aspect length as $\rho \longrightarrow + \infty$ for a second time. From rotational symmetries in the proof of \textbf{Lemma} $9$, there are six possible rotations from which $C_j$ can occur, in which $\mathcal{C} \equiv \mathscr{C}_{2_j}$, $\mathcal{C} \equiv \mathscr{C}_{3_j}$ or $\mathcal{C} \equiv \mathscr{C}_{4_j}$. Each upper bound under wired and free boundary conditions has been shown. \boxed{ }

\section{Vertical and horizontal strip densities}

\subsection{Towards proving horizontal and vertical crossing densities in Definition $1$}

In this section, we make use of strip densities similar to those provided for the random cluster model in {[14]} (defined in Section \textit{3.3}) from which strip density and renormalization inequalities will be presented, in the infinite length aspect ratio limit. In the arguments below, we present boxes $\mathscr{H}, \mathscr{H}_i$ and $\mathscr{H}^{\prime}_i$ across which horizontal and vertical crossings are quantified. For the lower bound of the conditional probability of obtaining no vertical crossings across each $\mathscr{H}_i$, we introduce a slightly larger hexagonal box $\mathscr{H}^{\mathrm{Stretch}}$ which has an aspect height ratio larger than that of $\mathcal{H}_j$.

\bigskip

\noindent \textbf{Definition} \textit{4} (\textit{dilute Potts horizontal and vertical strip densities}): For $N \geq 1$, $x \leq \frac{1}{\sqrt{n}}$, $nx^2 \leq \mathrm{exp}(-|h^{\prime}|)$, and $(n,x,h,h^{\prime})$, with external fields $h, h^{\prime}$ taken as real parameters, the strip density for horizontal crossings across $\mathcal{H}_j$ under the Spin measure with free boundary conditions is,

\begin{align*}
        p^{\mu}_N =  \mathrm{lim \text{ } sup}_{\rho \rightarrow \infty} \big\{ \mu^0_{
        [ 0 ,\rho N] \times_H [0 ,\lambda \mathrm{Stretch} ]}  \big[  \mathcal{H}_{[0, \rho N ] \times_H [0 , \lambda \mathrm{Stretch}  ] }  \big]^{\frac{1}{\rho}}   \big\} \text{, } 
\end{align*}

\noindent while for vertical crossings across $\mathcal{H}+j$, under the Spin measure with wired boundary conditions, is,

\begin{align*}
      q^{\mu}_N =  \mathrm{lim \text{ } sup}_{\rho \rightarrow \infty} \big\{ \mu^1_{[ 0 , \rho N] \times_H [0 , \lambda \mathrm{Stretch} ]}\text{ } \big[  \mathcal{V}^c_{[0, \rho N ] \times_H [0 , \lambda \mathrm{Stretch}]} \big]^{\frac{1}{\rho}}   \big\} \text{. } 
\end{align*}

\bigskip

\noindent We denote $p_N \equiv p^{\mu}_N$ and $q_N \equiv q^{\mu}_N$. With these quantities, we prove the strip density formulas which describe how boundary conditions induced by vertical crossings under wired boundary conditions across $ H_{j + \delta_k} , H_{j + \delta_l}  \subsetneq  \mathcal{H}_j$ relate to horizontal crossings under free boundary conditions.

In the proof below, we make use of arguments from \textbf{Lemma} $11$ to study vertical crossings across hexagons, and through applications of $(\mathcal{S}-\mathrm{SMP})$ and $(\mathcal{S}-\mathrm{CBC})$. To prove \textbf{Lemma} $1$, we define additional crossing events as follows. First, the crossing event that three hexagons, with aspect width of $\mathcal{H}_j$ and aspect length $\mathrm{Stretch}$ placed on top of each other, is pushed forwards to apply FKG-type arguments, with $(\mathrm{FKG})$, over a countable intersection of horizontal crossings across hexagons with the same aspect height and smaller aspect length than that of $\mathcal{H}_j$. We denote this event with $\mathcal{E}$. Second, we also need the event of obtaining a horizontal crossing across $H_{j , j +\delta}$ and $H_{j , j - \delta}$, conditioned on $\mathcal{E}$ which we denote as $\{\mathcal{F} | \mathcal{E} \}$. We study the conditions under which wired boundary conditions distributed from a prescribed distance of $H_{j , j - \delta}$ and $H_{j , j + \delta}$ induce vertical crossings. Third, crossing events across a larger domain than those considered in $\{ \mathcal{F} | \mathcal{E} \}$ are formulated by making use of the monotonicity in the domain assumption, denoted as $\mathcal{G}$ which is independently of $\rho$.

\noindent Fourth, the intersection of the previous three events is pushed forwards, and by virtue of $(\mathcal{S}-\mathrm{SMP})$ and $(\mathcal{S}-\mathrm{CBC})$, yields a strip inequality relating $p_n$ to $q_n$, and $q_n$ to $p_n$. As the aspect length $\rho \longrightarrow + \infty$, inequalities corresponding to the horizontal and vertical strip densities are presented.

\noindent \textit{Proof of Lemma $1$}. The argument consists of six parts; we fix $\lambda \in \textbf{N}$, $N \in 3 \textbf{N}$. As a matter of notation, below we denote each of the three boxes below as the Cartesian product of the aspect length and height ratios, and let $\rho \longrightarrow \infty$ in the last step. In the boxes $\mathscr{H}$, $\mathscr{H}_i$ and $\mathscr{H}^{\prime}_i$ below, $\lambda$ is taken smaller relative to $\rho$. Under the definitions of $\mathcal{E}$, $\{ \mathcal{F} | \mathcal{E} \}$ and $\mathcal{G}$, we first define all hexagonal boxes across which horizontal crossings occur, which are defined as,

{\small \begin{align*}
           \mathscr{H}   =   [ 0 , \rho N]    \times_H  [  0   ,       (2 \lambda) \text{ }  \mathrm{Stretch}     +  \mathrm{Stretch}   ]   \text{, } \\  \\
            \mathscr{H}_i   =  [0 , \rho N]   \times_H  [ (2 i )  \mathrm{Stretch}  + \mathrm{Stretch}   ,  (2i) \mathrm{Stretch}  + 2 \mathrm{Stretch}  ]   \text{, }  \\   \\
           \mathscr{H}^{\prime}_i  =   [0 , \rho N] \times_H   [      (2i)  \mathrm{Stretch}     ,      (2i)  \mathrm{Stretch}  +  \mathrm{Stretch}    ]  \text{, } 
\end{align*} }

\noindent for every $0 \leq i \leq \lambda - 1$. As indicated above, the notation $\times_H$, for some nonempty subset $[ 0, a] \times_H [0 , b]$ of $\textbf{H}$, denotes that the finite volume over the hexagonal lattice of length $a$ and height $b$, for $a$ and $b$ strictly positive. In the construction, the aspect length is the same as that of $\mathscr{H}$, while the aspect height of each box is partitioned in $i$ relative to the scaling of the $\mathrm{Stretch}$ factor. Also, a final box with the $\mathrm{Stretch}$ scaling itself will be defined, 

{\small \begin{align*}
        \mathscr{H}_{\mathrm{Stretch}} =   [0 , \rho N ] \times_H [  0  ,  N   \lambda\mathrm{Stretch}   ]  \text{, } 
\end{align*} } 

\noindent which is supported over which the spin measure with wired bound conditions for a lower bound of $\mu^1_{\mathscr{H}}[\mathcal{F}  | \mathcal{E}]$, and $n$ is an integer parameter. Second, to apply $(\mathrm{FKG})$ previously used, if $\mathcal{H}_{\mathcal{D}}$ denotes a horizontal crossing across a finite domain $\mathcal{D}$ of $\mathcal{S}$, we make use of $ \mathscr{H}, \mathscr{H}_i, \mathscr{H}^{\prime}_i  \subsetneq  \mathcal{D}$ with smaller aspect lengths across which horizontal crossings occur. The lower bound for applying $(\mathrm{FKG})$ across a countable family of horizontal crossings $\mathcal{H}_{\mathscr{H}_i}$ is,

\begin{align*}
      \mu^1_{\mathscr{H}}[   \mathcal{E}   ]   \geq   \mu^1_{\mathscr{H}} \big[  \bigcap_{0 \leq i \leq \lambda-1} \mathcal{H}_{\mathscr{H}_i}  \big]    \geq  \prod_{0 \leq i \leq \lambda - 1} \mu^1_{\mathscr{H}}[      \mathcal{H}_{\mathscr{H}_i}   ]  \geq   \prod_{0 \leq i \leq \lambda-1} \frac{1}{\big( \lambda_i^C \big)^{\rho}} \geq \frac{1}{\big( \lambda^C \big)^{\lambda  \rho}  }  \text{, } 
\end{align*}

\noindent with the existence of the lower bound guaranteed by \textit{Corollary} $11$, and $\lambda$ is the minimum amongst all $\lambda_i$. Before letting $\rho \longrightarrow + \infty$, pushing forwards the horizontal crossing event across $\mathscr{H} \subset \mathcal{D}$ under wired boundary conditions for vertical crossings across $\mathscr{H}^{\prime}_i$ gives,

\begin{align*}
    \mu^1_{\mathscr{H}}[ \mathcal{F}   |   \mathcal{E}  ]   \geq  \mu^1_{\mathscr{H}^{\prime}_i} \big[   \bigcap_{0 \leq i \leq \lambda } \mathcal{V}^c_{\mathscr{H}^{\prime}_i} \big] \text{ }  \text{, }
    \end{align*}

    \noindent which is bounded below by the product of independent events through repeated applications of $(\mathrm{FKG})$ across $\lambda - 1$ crossing events,
    
    \begin{align*}
    \prod_{0 \leq i \leq \lambda - 1} \mu^1_{\mathscr{H}^{\prime}_i}[      \mathcal{V}^c_{\mathscr{H}^{\prime}_i}   ]  \geq \mu^1_{[ 0 , \rho n ] \times_H  [   0   , n_1 \lambda \mathrm{Stretch}  ]}   \big[ \mathcal{V}^c_{[ 0, \rho n ] \times_H [ 0 ,     n_1   \lambda \mathrm{Stretch}  ]   }           \big]^{\lambda+1} \text{, } 
\end{align*}

\noindent for $n_2 > n_1$, from applications of $(\mathcal{S} -\mathrm{SMP})$, monotonicity in the domain, and applying $(\mathrm{FKG})$ to vertical crossing events, instead of horizontal crossing events. 

\bigskip

\noindent By construction of $\mathcal{E}$, the following lower bound for the conditional event $\{ \mathcal{F}   |  \mathcal{E} \}$,

\begin{align*}
       \mu^1_{\mathscr{H}}[   \mathcal{E} \cap \mathcal{F} \cap \mathcal{G}    ] \\ \\  = \mu^1_{\mathscr{H}}[ (\mathcal{E} \cap \mathcal{F} ) \cap \mathcal{G}   ] \overset{(\mathrm{FKG})}{\geq} \mu^1_{\mathscr{H}} [\mathcal{E} \cap \mathcal{F}      ]  \mu^1_{\mathscr{H}}[    \mathcal{G}    ] \\ \\ = \bigg[ n^{k_{\mathcal{G}}(\sigma)} \text{ } \mathrm{exp} \bigg[ {   \#  \big\{ \textit{hexagonal boxes } \mathcal{H}  : 1 \leq   i \sim j   \leq 6  ,  \sigma_v \in \{ \pm 1 \}^{v (\textbf{H})}  ,  \textbf{1}_{ \{ \sigma_{v_i} = \sigma_{v_j}  = 1 \} }     }\big\} \bigg]  \\ \times    x^{2          \{ e \text{ } : \text{ }  e \in (2 \lambda) \mathrm{Stretch} + \mathrm{Stretch}  \}  }     \bigg] \mu^1_{\mathscr{H}}[ \mathcal{E} \cap \mathcal{F}      ]         \text{, }   
\end{align*}    

\noindent where the edge weight in the lower bound is representative of additional weights in the configuration supported under the wired Spin measure over $\mathscr{H}$, which we denote as $(\mathrm{FKG})- (\mathcal{S}-\mathrm{SMP})$. In the exponential of the first external field, $\mathcal{H} \subsetneq   (2 \lambda) \mathrm{Stretch} + \mathrm{Stretch}$. In the exponent of the edge weight $x$, by definition the number of connected components in the configuration $\sigma$ sampled under $\mu$ is,

\begin{align*}
         k_{\mathcal{G}}(\sigma) =  \sum_{k \text{ } \mathrm{cc\text{ } 's}}  \bigg[ \sum_{u \sim v} \textbf{1}_{  \{ \sigma_{u_k} \equiv \sigma_{v_k} \} } + 1 \bigg]  \text{, } 
\end{align*}

\noindent where the summation is taken over all connected components $k$ so that $\mathcal{G}$ occurs, and all neighboring vertices $u$ and $v$ with the same $\pm$ spin, clearly impacting the number of connected components counted under $\sigma$.

Before completing the next step, we combine the estimates on $\mu^1_{\mathscr{H}}[\mathcal{E}]$ and $\mu^1_{\mathscr{H}}[\mathcal{F} \text{ } | \text{ } \mathcal{E}]$ to obtain the strip inequality between horizontal and vertical crossings. The following comparison amounts to making use of $(\mathcal{S}-\mathrm{SMP})$ and (MON) to establish the following. First, we know that the measure $\mu^1_{\mathcal{S}}[ \cdot ]$ can be bound above with,

\begin{align*}
        \mu^1_{\mathcal{S}}[    \mathcal{E} \cap \mathcal{F} \cap \mathcal{G}      ]     \leq \mu^1_{\mathcal{S}}[  \mathcal{E}  |    \mathcal{F} \cap \mathcal{G}    ]    \leq    \mu^0_{[0 , \rho N] \times_H [ 0 ,  n_2 \lambda \mathrm{Stretch}]}\text{ }  \big[    \mathcal{H}_{[0, \rho N] \times_H [ 0 ,   n_3 \lambda  \mathrm{Stretch}   ]}         \big]^{\lambda}  \text{, } 
\end{align*}

\noindent because the event in the upper bound is more likely to occur than the event in the lower bound, in addition to an application of $(\mathrm{FKG})$ for $\lambda$ horizontal crossings, each event of which has equal probability, across $[0 , \rho N] \times_H [ 0 , n_3  \lambda   \mathrm{Stretch} ]$, and $n_3 > n_2$. Also, the upper bound to the conditional probability $\mu^1_{\mathcal{S}}[    \mathcal{E}  |  \mathcal{F} \cap \mathcal{G} ]$ above is established by making the comparison between measures with free boundary conditions. In comparison to the planar case, modifications to the argument with $(\mathcal{S}-\mathrm{SMP})$, while other properties of the random cluster measure $\phi \big[ \cdot \big]$ directly apply.

In light of the lower bound in $(\mathrm{FKG}) - (\mathcal{S}-\mathrm{SMP})$ dependent on the edge weight $x$ and $\mathrm{Stretch}$, we consider the horizontal pushforward from the previous upper bound with free boundary conditions,

\begin{align*}
        \mu^0_{[ 0 , \rho N ] \times_H [ 0 , n_2  \lambda  \mathrm{Stretch}] }  \big[  \mathcal{H}_{[0 , \rho N] \times_H [ 0 , n_1    \lambda \mathrm{Stretch}   ]}   \big]   \text{, } 
\end{align*}

\noindent which can be bound below by establishing comparisons between the measure under wired boundary conditions supported over a smaller hexagonal domain, 

\begin{align*}
           \mu^1_{[ 0 , \rho N ] \times_H [ 0 , n_1  \lambda  \mathrm{Stretch}] }  \big[      \mathcal{V}^{\text{ } c}_{[0 ,  \rho N  ] \times_H [ 0 ,  n_1  \lambda \mathrm{Stretch}]}      \big]   \text{, } 
\end{align*}

\noindent as a consequence yielding one estimate for the vertical strip density,

{\small \begin{align*}
     \mu^0_{[0 , \rho N] \times_H [0 \text{ } ,  n_2 \lambda \text{ } \mathrm{Stretch} ]}  \big[   \mathcal{H}_{[0 ,  \rho N  ] \times_H [ 0  ,   n_1 \lambda \text{ } \mathrm{Stretch}   ]}     \big]^{\lambda}   \overset{\rho  \longrightarrow +\infty }{\approx}  \big( p_{\mathrm{Stretch} \text{ }    N} \big)^{\lambda} \\ \\  \geq \mu^1_{[0 ,  \rho N  ] \times_H [  0  , n_1 \lambda \text{ } \mathrm{Stretch}  ] }\big[ \mathcal{V}^{\text{ }c}_{[ 0 , \rho N ] \times_H [ 0 ,   n_1 \lambda \text{ } \mathrm{Stretch}     ]}   \big]^{\lambda + 1}      \overset{\rho  \longrightarrow +\infty }{\approx}      \frac{ 1 }{\lambda^C} \big( q_{\mathrm{Stretch}\text{ } N}     \big)^{\lambda+1}         \text{, }  
\end{align*}     }

\noindent due to the fact that,

{\small \begin{align*}
 \bigg[  n^{k_{\mathcal{G}}(\sigma)}  x^{2   \# \{  e :  e  \in  (2 \lambda) \mathrm{Stretch}  + \mathrm{Stretch} \}  }  \mathrm{exp}  \bigg[ {   \#  \big\{ \textit{hexagonal boxes }  \mathcal{H}  :  1 \leq   i \sim j   \leq 6  ,  \sigma_v \in \{ \pm 1 \}^{v(\textbf{H})}     ,  \textbf{1}_{\{ \sigma_{v_i} = \sigma_{v_j}  = 1\} }   \big\}   }\bigg]     \bigg]^{{ \lambda \rho}^{-1}}    \text{, } 
\end{align*} }

\noindent approaches $1$ as $\rho \longrightarrow \infty$, and cancellations of $\lambda$ crossings in the inequality gives,

\begin{align*}
        p_{\mathrm{Stretch} \text{ } N} \geq \frac{1}{\lambda^C} \big( q_{\mathrm{Stretch} \text{ } N} \big)^{\mathrm{Stretch} + \frac{\mathrm{Stretch}}{\lambda}}  \text{, }  
\end{align*}

\noindent resulting from $(\mathcal{S} -\mathrm{CBC})$, as in the region below the connected component of the path associated with the crossing $\mathcal{H}_{\mathscr{H}_1}$ the induced boundary conditions dominate the measure supported over a smaller domain under wired boundary conditions.

The strip inequality for horizontal and vertical crossings is finally achieved by taking each side of the inequality to the power $\frac{1}{\rho \lambda}$, which preserves the direction of the inequality as a monotonic decreasing transformation. As $\rho \longrightarrow + \infty$, we recover the peculiar definition of the horizontal strip density, while the other inequality corresponding to the vertical crossing density can be easily achieved by following the same argument, with the exception of the inequalities leading to the final estimate for vertical crossing events instead of horizontal ones. \boxed{}

\subsection{Pushing lemma}

We turn to the following estimates. In \textbf{Lemma} $13$ and \textbf{Lemma} $14$ below, $\bar{\mathscr{H}}$ denotes the box with aspect length $\rho n$, and variable aspect height defined for each box in the proof. To prove $\textbf{Lemma} \text{ } 13$ (stated below), we make use of the following property for the Spin measure. With the Pushing Lemma, we provide arguments for the renormalization inequalities in the next section. Under this Lemma, we proceed to obtain results for the (PushPrimal) and (PushDual) conditions (listed below), from which are the combined with the accompanying (PushPrimal Strip), and (PushDual Strip), respectively, to probabilistically quantify crossing probabilities across hexagons  $\widetilde{\mathscr{H}}_i$, $ \widetilde{\mathscr{H}}^+_i$, and $ \widetilde{\mathscr{H}}^-_i$ (introduced in Section \textit{8} to obtain the vertical and horizontal strip density formulas). 

\bigskip

\noindent \textbf{Property} (\textit{Finite energy for the Spin measure}, {[8]}): For any $\tau \in \{-1,1 \}^{\textbf{T}}$ and $\sigma \in \Sigma(G,\tau)$, $\mu^{\tau}_{G , n ,x , h , h^{\prime}}[\sigma] \geq \epsilon^{|G|}$, for any $\epsilon > 0$ depending only on $(n,x,h,h^{\prime})$.

\subsubsection{Statement of the two Lemmas for strip and planar domains}.

\noindent \textbf{Lemma} $13$ (\textit{Pushing Lemma}): There exists positive $c_1 \equiv  c_1(n,k_i,x_i,\mathrm{Stretch})$, such that for every $n \geq 1$, with aspect length $\rho$, one of the following two inequalities is satisfied,

\begin{align*}
          \mu^{\mathrm{Mixed}}_{\bar{\mathscr{H}}}[     \mathcal{H}_{\mathscr{H}}      ]   \geq c_1^{\rho}   \tag{PushPrimal}   \text{ , } 
\end{align*}

\noindent or, 

\begin{align*}
             \mu^{(\mathrm{Mixed})^{\prime}}_{\bar{\mathscr{H}}}[      \mathcal{V}^{\text{ } c}_{\mathscr{H}}          ]    \geq  c_1^{ \rho}  \tag{PushDual}    \text{ , } \\
\end{align*}

\noindent for every $\rho \geq 1$, and the superscript $\mathrm{Mixed}$ denotes wired boundary conditions along the left, top and right sides of $\bar{\mathscr{H}}$, and free boundary conditions elsewhere. $\mathcal{H}$ is the same hexagonal box used in previous arguments for \textbf{Lemma} $1$. Under the PushDual condition, the analogous statement holds for the complement of vertical crossings across $\mathscr{H}$, under dual boundary conditions $(\mathrm{Mixed})^{\prime}$ to $\mathrm{Mixed}$.

\bigskip

\noindent \textbf{Lemma} $14$ (\textit{Pushforward of horizontal and vertical crossings under mixed boundary conditions}): There exists positive $c_2=c_2(n,k_i,x_i,\mathrm{Stretch})$ such that for every $N \geq 1$, with aspect length $\rho$, one of the following two inequalities is satisfied,

\begin{align*}
           \mu^{(\mathrm{Mixed})^{\prime}}_{\mathcal{S}}[  \mathcal{H}_{\mathscr{H}}   ]         \geq c_2^{\rho}   \tag{PushPrimal Strip}     \text{, } 
\end{align*}

\noindent or,

\begin{align*}
                \mu^{(\mathrm{Mixed})^{\prime}}_{\mathcal{S}}[  \mathcal{V}^{ c}_{\mathscr{H}}    ]         \geq c_2^{\rho}         \tag{PushDual Strip}     \text{, } 
\end{align*}

\noindent for every $\rho \geq 1$. $(\mathrm{Mixed})^{\prime}$ denote the same boundary conditions from \textbf{Lemma} $13^{*}$, which manifests in the following.

\subsubsection{Proof of Lemma $14$ for (PushPrimal Strip) $\&$ (PushDual Strip) conditions}.

\noindent With some abuse of notation, we denote the hexagonal boxes for this proof as, 

\begin{align*}
       \mathscr{H}_i = [0 , 2 N] \text{ }  \times_H \text{ } \big[ \frac{i}{3} N , \frac{i+1}{3} N  \big] \text{, } 
\end{align*}

\noindent for $i=0,1,2$. Furthermore, we introduce the vertical segments along the bottom of each $\mathscr{H}_i$, and hexagons with same aspect length as those of each $\mathscr{H}_i$, in addition to hexagon of the prescribed aspect height below, respectively,

\begin{align*}
         \mathcal{I}_i =   \big[  \frac{i}{3} N  ,  \frac{i+1}{3} N  \big]   \times \{0\}  \text{, } \\ \\ 
          \mathcal{K}_i =  \big[  \frac{i}{3} N  ,   \frac{i+1}{3} N  \big]   \times_H    [ - N  , N ]    \text{, }      
\end{align*}

\noindent each of which are also indexed by $i$, with the exception that $i$ also runs over $i= 4,5$. Before presenting more arguments for the connectivity between $\mathcal{I}_1$ and $\mathcal{I}_4$, suppose that either $\mu^{(\mathrm{Mixed})^{\prime}}_{\mathcal{S}}[    \mathcal{V}_{\mathscr{H}_i}    ] \geq \frac{1}{6}$, or $\mu^{(\mathrm{Mixed})^{\prime}}_{\mathcal{S}}[   \mathcal{H}^c_{\mathscr{H}_i}  ] \geq \frac{1}{6}$ for some $i$. In the first case for the pushforward of vertical crossings in $\mathscr{H}_i$, another application of the $\mu$ homeomorphism $f$ from arguments to prove \textbf{Corollary} $11$ implies that (PushPrimal Strip) holds, while in the second case for the pushforward of horizontal crossings in $\mathscr{H}_i$, an application of the same homeomorphism implies that (PushDual Strip) holds. By complementarity, under $\mu^{(\mathrm{Mixed})^{\prime}}_{\mathcal{S}}[\cdot ]$, the pushforward of the following events respectively satisfy the lower bounds, as $\mu^{(\mathrm{Mixed})^{\prime}}_{\mathcal{S}}[  \mathcal{V}^c_{\mathscr{H}_i}  ] \geq \frac{5}{6}$, and $\mu^{(\mathrm{Mixed})^{\prime}}_{\mathcal{S}}[   \mathcal{H}_{\mathscr{H}_i} ] \geq \frac{5}{6}$. The same argument that follows applies to lower bounds for crossing probabilities by other constants than $\frac{1}{6}$ or $\frac{5}{6}$ which are provided in {[14]}, modifications to obtaining identical lower bounds in place of different constants are provided with the following.

With such estimates, under the same boundary conditions listed in (PushPrimal Strip) $\&$ (PushDual Strip), the Spin measure satisfies

\begin{align*}
        \mu^{(\mathrm{Mixed})^{\prime}}_{\mathcal{S}}[ \mathcal{V}^c_{\mathscr{H}_0}    \cap  \mathcal{H}_{\mathscr{H}_1} \cap    \mathcal{V}^c_{\mathscr{H}_2}   ] \leq   \mu^{(\mathrm{Mixed})^{\prime}}_{\mathcal{S}}[    \mathcal{H}_{\mathscr{H}_1}          ]     \leq \mu_{\mathscr{H}}^{* - (\mathrm{Mixed})^{\prime}} [  \mathcal{H}_{\mathscr{H}_1}   ]     \text{, }   
\end{align*}

\noindent where the upper bound for the probability of the intersection of the three events above only holds under boundary conditions in which the incident layer to the configuration (as given in arguments for the proof of \textbf{Lemma} $9^{*}$), the boundary conditions for the measure dominating $(\mathrm{Mixed})^{\prime}$ boundary conditions holds because every vertex that is wired in the $(\mathrm{Mixed})^{\prime}$ boundary conditions is also wired in the boundary conditions for the pushforward in the upper bound. Moreover, the partition of boundary vertices in the boundary conditions for the upper bound is composed of the arc that is wired in the boundary conditions for , in addition to a singleton 

Under $(\mathrm{Mixed})^{\prime}$ boundary conditions introduced for the high-temperature spin measure above, the conditional probability

\begin{align*}
        \mu^{\mathrm{Mixed}}_{\mathcal{S}}[  \mathcal{I}_1   \overset{\mathscr{H}}{\longleftrightarrow}  \mathcal{I}_4   | \mathcal{H}_{\mathscr{H}_1}    ]       \text{, } 
\end{align*}

\noindent can be bound below by conditioning on a horizontal crossing $\mathcal{H}_{\mathscr{H}_1}$ across $\mathscr{H}_1$. In particular, conditionally on $\mathcal{H}_{\mathscr{H}_i}$, the connectivity event 

\begin{align*}
         \mu^{(\mathrm{Mixed})^{\prime}}_{\mathcal{S}}[    \mathcal{I}_1   \overset{\mathcal{K}_1}{\longleftrightarrow}  \mathcal{I}_4     |     \mathcal{H}_{\mathscr{H}_1}        ]              \text{ , } 
\end{align*}

\noindent can be bounded below as shown above through applications of $(\mathcal{S}-\mathrm{CBC})$ and $(\mathcal{S}-\mathrm{SMP})$. Each property is applied as follows; for $(\mathcal{S}-\mathrm{SMP})$, we make use of previous partitions of the incident layer of hexagons to a configuration, in which $(\mathcal{S}-\mathrm{SMP})$ can only be applied when the outermost layer of a configuration can be partitioned into two equal sets over which the $\pm$ spin is constant.

Concluding, we apply standard arguments for the crossing event below through a lower bound dependent on a conditional probability,

\begin{align*}
        \mu^{\mathrm{Mixed}}_{\mathcal{S}}[     \mathcal{I}_1 \overset{\mathscr{H}}{\longleftrightarrow} \mathcal{I}_4    ]   \geq \mu^{\mathrm{Mixed}}_{\mathcal{S}} [   \mathcal{I}_1 \overset{\mathscr{H}}{\longleftrightarrow} \mathcal{I}_4      |  \mathcal{H}_{\mathscr{H}_i}  ]  \mu^{\mathrm{Mixed}}_{\mathcal{S}}[  \mathcal{H}_{\mathscr{H}_i}   ]  \geq  \mathscr{C}   \mu^{\mathrm{Mixed}}_{\mathcal{S}} [     \mathcal{I}_1 \overset{\mathscr{H}}{\longleftrightarrow} \mathcal{I}_4    | \mathcal{H}_{\mathscr{H}_i}         ]    \\    \\
        \geq      \frac{\mathscr{C}}{        \prod_{i=1}^{\alpha}   n^{k_i}   x^{e_i}  \mathrm{exp}(  h_i   )  }     \text{ }           \text{, }  
\end{align*}

\noindent from which an application of $(\mathrm{FKG})$, given suitable $\mathscr{C}>0$, for the countable intersection, dependent on $i$, of horizontal crossings across hexagons of sufficiently small aspect length $\mathrm{Stretch}_i$. The inverse proportionality in the lower bound is dependent on the product $\mathcal{T}$, defined in the proof for \textbf{Lemma} $1$ with $i$ running over two configurations with respective number of connected components $k_1 + 1$ and $k_2 + 1$. The lower bound dependent on the edge weight $x$ arises from multiple applications of $(\mathcal{S} -\mathrm{SMP})$ and (MON), in which the modification to (SMP) from the random cluster model argument with $(\mathcal{S}-\mathrm{SMP})$ for the Spin Measure results in comparisons between $\pm$ configurations and partitions of the incident layer as described in Section \textit{5.2}. 

Instead, if we suppose that the lower bounds for $\mu^{(\mathrm{Mixed})^{\prime}}_{\mathcal{S}}[   \mathcal{V}^c_{\mathscr{H}_i}   ] \geq c^{\prime \prime}$ for real $c^{\prime \prime}$, the lower bound on the second line above takes the form,

\begin{align*}
           \mu^{\mathrm{Mixed}}_{\mathcal{S}}[   \mathcal{I}_1 \overset{\mathscr{H}}{\longleftrightarrow} \mathcal{I}_4     | \mathcal{H}_{\mathscr{H}_1}      ] \geq  \frac{\mathscr{C}  c^{\prime \prime}}{\prod_{i=1}^{\alpha}  n^{k_i} x^{e_i} \mathrm{exp}(h_i)   }      \text{. } 
\end{align*}

\noindent due to the fact that the boundary conditions from the special case of the inequality, where the power to which the product of the edge weight and difference in monochromatically colored triangles is raised to the aspect ratio $\mathrm{Stretch}$ of $\mathscr{H}_{\mathrm{Stretch}}$, and the number of connected components in the exponent of $n$ is the difference between the number of connected components of a $\pm$ configuration respectively sampled under $\mu_{\mathcal{S}}^{(\mathrm{Mixed})^{\prime}}$ and $\mu_{\mathscr{H}}^{* - (\mathrm{Mixed})^{\prime}}$. 

Furthermore, the lower bound dependence on the edge weight $x$,$n$ and $e$, emerges from an application of (FKG) to the pushforward below of the connectivity event between $\mathcal{I}_1$ and $\mathcal{I}_4$, bounded below above,

\begin{align*}
             \mu^{\mathrm{Mixed}}_{\mathcal{S}}[      \mathcal{I}_1 \overset{\mathscr{H}}{\longleftrightarrow}  \mathcal{I}_4    ]  \geq                  \mu^{\mathrm{Mixed}}_{\mathcal{S}} \big[    \mathcal{I}_1 \overset{\mathscr{H}}{\longleftrightarrow} \mathcal{I}_4 \big| \mathcal{H}_{\mathscr{H}_1}            \big]  \mu^{\mathrm{Mixed}}_{\mathcal{S}} [   \mathcal{H}_{\mathscr{H}_1}    ]   \text{, }    
\end{align*}

\noindent which can be further bounded below by the product of crossing probabilities,

{\small \begin{align*}
     \frac{\mathscr{C}    c^{\prime \prime}}  {\big[ {\prod_{i=1}^{\alpha}  \text{ } n^{k_i} \text{ }   x^{e_i} \text{ }  \mathrm{exp} (  h_i  )      }  \big]^{2} }    \text{, }   
\end{align*} }

\noindent Observe that the horizontal crossing pushed forwards in the inequality above yields the desired pushforwards in the PushDual condition, as the previously mentioned application of $(\mathrm{FKG})$ yields,

{\small \begin{align*}
\mu  [    \mathcal{H}_{\mathscr{H}}       ]   \geq  \mu  \big[    \bigcap_{i=1}^{\alpha}    \mathcal{H}_{\mathscr{H}_i}   \big] \geq \prod_{i=1}^{\alpha} \mu  [    \mathcal{H}_{\mathscr{H}_i}           ]  \geq      \frac{\mathscr{C} c^{\prime\prime}  }{\big[              {\prod_{i=1}^{\alpha}   n^{k_i}  x^{e_i}  \mathrm{exp}( h_i )}          \big]^{2}}       \times \overset{\alpha -2}{\cdots} \times    \frac{\mathscr{C} c^{\prime\prime}  }{\big[              {\prod_{i=1}^{\alpha}   n^{k_i}  x^{e_i}  \mathrm{exp}( h_i )}          \big]^{2}}   \\ \\  \equiv   \frac{ \big(\mathscr{C} c^{\prime\prime}       \big)^{\alpha}  }{ \big[             {\prod_{i=1}^{\alpha}   n^{k_i}  x^{e_i} \text{ } \mathrm{exp}( h_i  )}          \big]^{2 \alpha} }   \text{, }
\end{align*} }

\noindent for crossings across each of the hexagons $\mathscr{H} = \cup_{i=1}^{\alpha} \mathscr{H}_i$, and where the respective powers $k_i$, $e_i$ and $h_i$ appear in powers of the number of loops, edges, and exponential for the first external field. The form of the constant is provided in the lower bound above, and the same argument can be applied to obtain constant corresponding to (PushDual Strip) for vertical hexagonal crossings to obtain the desired constant in the lower bound, concluding the proof. \boxed{}

\subsection{Lemma $13$ arguments from Strip conditions}

\noindent \textit{Proof of Lemma $13$}. We show that either $(\mathrm{PushPrimal \text{ } Strip}) \Rightarrow (\mathrm{PushPrimal})$, or that $(\mathrm{PushDual \text{ } Strip}) \Rightarrow (\mathrm{PushDual})$. Without loss of generality, suppose that $(\mathrm{PushDual \text{ } Strip})$ holds; to show that $(\mathrm{PushDual})$ holds, we introduce the following collection of similarly defined boxes from arguments in \textbf{Lemma} $14$ on the previous page,

\begin{align*}
          \widetilde{\widetilde{\mathscr{H}_i}} = [ 0  ,  \rho N ]   \times_H  \big[   \frac{i}{3} N  ,    \frac{i+1}{3} N     \big]      \text{ , } 
\end{align*}

\noindent for $1 \leq i \leq N$, with $N$ sufficiently large. Under $(\mathrm{Mixed})^{\prime}$ boundary conditions,

\begin{align*}
          \mu^{(\mathrm{Mixed})^{\prime}}[\mathcal{V}^{c}_{\widetilde{\widetilde{\mathscr{H}_N}}}    ]   \geq c^{ \rho}    \text{, } 
\end{align*}

\noindent the probability of a complement of the vertical crossing across $\widetilde{\widetilde{\mathscr{H}_N}}$, and can be bounded below by $c^{\text{ } \rho}$ because by assumption $(\mathrm{PushPrimal \text{ } Strip})$ holds. Clearly, the probability of obtaining a vertical crossing across the last rectangle over all $i$ can be determined by applying the FKG inequality across each of the $N$ smaller hexagons, yielding an upper bound of $c^{N \rho}$ to the probability of obtaining a longer $N$-hexagon crossing. 

Next, with similar conditioning on horizontal crossings in previous arguments, the probability of a horizontal crossing across $\widetilde{\widetilde{\mathscr{H}_i}}$, given the occurrence of a horizontal crossing across $\widetilde{\widetilde{\mathscr{H}_{i+1}}}$, satisfies for every $i$,

\begin{align*}
    \mu_{\mathcal{S}}^{(\mathrm{Mixed})^{\prime}}[  \mathcal{V}^c_{\widetilde{\widetilde{\mathscr{H}_i}}}  | \mathcal{V}^c_{\widetilde{\widetilde{\mathscr{H}_{i+1}}}}    ]  \geq c^{\rho}     \text{, }
\end{align*}

\noindent with the exception that the pushforward  $\widetilde{\widetilde{\mathscr{H}_{i+1}}}$, taken under $(\mathrm{Mixed})^{\prime}$ boundary conditions, in comparison to previous arguments for the wired pushforward

\begin{align*}
          \mu^1_{\mathcal{H}_j}[ \mathcal{V}_{\textbf{1}_{\{j + \delta_k \equiv j + \delta_l\}}}   ]    \text{, } 
\end{align*}

\noindent below by $e^{-c}$ for $\textbf{Corollary}$ $11$, can also be applied to bound the intersection of conditional events, for the event $\{  \mathcal{V}^c_{\widetilde{\widetilde{\mathscr{H}_i}}}  | \text{ } \mathcal{V}^c_{\widetilde{\widetilde{\mathscr{H}_{i+1}}}}   \}$, for all $i$, 

\begin{align*}
          \prod_{0 \leq i \leq N} \mu^{(\mathrm{Mixed})^{\prime}}_{\bar{\mathscr{H}}} [  \mathcal{V}^c_{\widetilde{\widetilde{\mathscr{H}_i}}} |  \mathcal{V}^c_{\widetilde{\widetilde{\mathscr{H}_{i+1}}}}   ] \geq  c^{N  \rho }  \text{, } 
\end{align*}

\noindent implying that the identical lower bound from the $(\mathrm{PushPrimal \text{ } Strip})$ holds, across the countable intersection of horizontal crossings,

\begin{align*}
         \mu^{(\mathrm{Mixed})^{\prime}}_{\bar{\mathscr{H}}}[  \mathcal{V}^{\text{ } c }_{\widetilde{\widetilde{\mathscr{H}_1}}}   ]  \geq c^{ N \rho}   \text{. } 
\end{align*}

\noindent We conclude the argument, having made use of the previous application of (FKG) across $0 \leq i \leq \lambda - 1$, uniformly in boundary conditions $(\mathrm{Mixed})^{\prime}$. \boxed{}

\section{Renormalization inequality}

We now turn to arguments for the Renormalization inequality. We make use of notation already given in the proof for the vertical and horizontal strip inequalities of $\textbf{Lemma}$ $1$, namely that we make use of a similar partition of the hexagons to the left and right of some $\mathscr{H}$. To restrict the crossings to occur across hexagons of smaller aspect length, we change the assumptions on our choice of $n$, and follow the same steps in the argument of $\textbf{Lemma}$ $1$ to obtain a lower bound for the pushforward $\mu^1_{\mathscr{H}}[\widetilde{\mathcal{E}} \cap \mathcal{F} \cap \mathcal{G}]$, where $\widetilde{\mathcal{E}}$ denotes the event that each of the three boxes $\widetilde{\mathscr{H}}_i, \widetilde{\mathscr{H}}_i^+ , \widetilde{\mathscr{H}}_i^-$ which are defined in arguments below. The partition of the aspect length of $\widetilde{\mathscr{H}}_i, \widetilde{\mathscr{H}}_i^+ , \widetilde{\mathscr{H}}_i^-$ is dependent on $i$. Also, the smaller scale over which we force the horizontal crossings to occur in $\widetilde{\mathcal{E}}$ is reflected in the partition of the aspect length, which not surprisingly permits for applications of $(\mathrm{FKG})$ with domains that are indexed by an auxiliary parameter for $0 \leq i \leq \lambda - 1$. The partition of $\mathscr{H}_i$ into the three boxes $\widetilde{\mathscr{H}}_i, \widetilde{\mathscr{H}}_i^+ , \widetilde{\mathscr{H}}_i^-$ determines corresponding powers, dependent on $\lambda$ to which the horizontal or vertical strip densities are raised before taking $\rho \longrightarrow + \infty$. As previously mentioned, differences in $(\mathcal{S} -\mathrm{SMP})$ emerge in one step of the following argument. We discuss the arguments for the proof when $(\mathrm{PushPrimal})$ holds, and in the remaining case when $(\mathrm{PushDual})$ holds, a modification to the argument is provided.

\subsection{Arguments for obtaining renormalization inequalities in the thermodynamic limit}

\noindent \textit{Proof of Lemma $2$}. Suppose that $(\mathrm{PushDual})$ holds; the $(\mathrm{PushPrimal})$ case will be discussed at the end. In light of the brief remark of the argument at the beginning of the section, we introduce the three boxes to partition the middle of $\mathscr{H}_i$ from the proof of \textbf{Lemma} $1$,

{\small \begin{align*}
     \widetilde{\mathscr{H}^{1,-}_i} \equiv   \widetilde{\mathscr{H}}_i^- =  [0 , \rho N]  \times_H  [  (2 i ) \mathrm{Stretch}  + \beta  \mathrm{Stretch} + \widetilde{\alpha_1}   \mathrm{Stretch}    , (2i) \mathrm{Stretch}   +  \beta  \mathrm{Stretch} \\  + \widetilde{\alpha_2} \mathrm{Stretch}   ]        \text{, } \\ \\
   \widetilde{\mathscr{H}^{1}_i}   \equiv    \widetilde{\mathscr{H}}_i  =  [0 , \rho N ]    \times_H   [ \text{ } (2 i ) \mathrm{Stretch}  +  \beta  \mathrm{Stretch} +  \widetilde{\alpha_2} \mathrm{Stretch}    ,  (2i)  \mathrm{Stretch}  +   \beta  \mathrm{Stretch} \\  + \widetilde{\alpha_3}  \mathrm{Stretch}   ]             \text{  , }  \\   \\
    \widetilde{\mathscr{H}^{1,+}_i}   \equiv  \widetilde{\mathscr{H}}_i^+   =   [0 , \rho N ]    \times_H \text{ }  [  (2 i )  \mathrm{Stretch}  +  \beta \mathrm{Stretch} + \widetilde{\alpha_3}   \mathrm{Stretch}    ,   (2i)  \mathrm{Stretch}  +  \beta  \mathrm{Stretch} \\ + \widetilde{\alpha_4} \mathrm{Stretch}    ]     \text{, }     
\end{align*} }

\noindent for every $ 0 \leq i \leq \lambda - 1$, and will apply steps of the argument from the proof of \textbf{Lemma} $1$, in which we modify all pushforwards under the prescribed boundary conditions for $\widetilde{\mathcal{E}}$. By construction, the boxes $\widetilde{\mathscr{H}}_i$, $\widetilde{\mathscr{H}}_i^+$, and $\widetilde{\mathscr{H}}_i^-$, each have the same $\rho$ aspect ratio, yet differ in the increment of the factors $\widetilde{\alpha_i} \in \beta \textbf{N}$, given $\beta$ sufficiently large, is given by,

\begin{align*}
   \widetilde{\alpha_i} = 1 + \frac{1}{\beta + i - 1} \text{. } 
\end{align*}

\noindent Briefly, we recall the steps with the sequence of inequalities below. Under one simple modification through the lower bound, we analyze the intersection of crossing probabilities as given in $(\mathrm{FKG})$, implying,

\begin{align*}
    \mu^1_{\mathscr{H}}[    \widetilde{\mathcal{E}}    ] \geq \prod_{0 \leq i \leq \lambda - 1} \mu^1_{\mathscr{H}}[    \mathcal{H}_{\widetilde{\mathscr{H}_i}}     ] \geq \frac{1}{(\lambda^{\prime})^{C\lambda \rho}}  \text{, } 
\end{align*}

\noindent from which the conditional probability dependent on $\widetilde{\mathcal{E}}$ can be bound from below as follows,

\begin{align*}
      \mu^1_{\mathscr{H}}[     \mathcal{F} |   \widetilde{\mathcal{E}}      ]   \geq \prod_{0 \leq i \leq \lambda - 1} \mu^1_{\widetilde{\mathscr{H}_i}}[    \mathcal{V}^c_{\widetilde{\mathscr{H}_i}}    ]     \geq    \mu^1_{[0 , \rho N] \times_H [0 , n^{\prime}_1 \lambda \text{ }  \mathrm{Stretch} ]} \big[  \mathcal{V}^{\text{ } c}_{[  0 , \rho N  ] \times_H  [  0 , N \lambda \mathrm{Stretch}  ]}         \big]^{\lambda + 1}       \text{. }   
\end{align*}

\noindent Further arguments result in the following lower bound for the probability of $\{ \widetilde{\mathcal{E}} \cap \mathcal{F} \cap \mathcal{G} \}$,

{\small \begin{align*}
       \mu^1_{\mathscr{H}}[    \widetilde{\mathcal{E}} \cap \mathcal{F} \cap \mathcal{G}    ]  = \mu^1_{\mathscr{H}}[   (\widetilde{\mathcal{E}} \cap \mathcal{F} ) \cap \mathcal{G}    ] \geq \mu^1_{\mathscr{H}} [\widetilde{\mathcal{E}} \cap \mathcal{F}      ]   \mu^1_{\mathscr{H}}[     \mathcal{G}    ] \\ \\  = x^{k_{\mathcal{G}}(\sigma)} \bigg[         x^{ \# \big\{      e  :  e \in (2 \lambda) \mathrm{Stretch} + \mathrm{Stretch} \big\} }  \mathrm{exp}       \bigg[     \#  \big\{  \textit{hexagonal boxes } \mathcal{H} :  1 \leq   i \sim j   \leq 6  ,  \sigma_v \in \{ \pm 1 \}^{v (\textbf{H} )} \\ ,  \textbf{1}_{\{ \sigma_{v_i} = \sigma_{v_j} = 1 \} }  \big\}     \bigg]  \bigg]  \mu^1_{\mathscr{H}}[    \mathcal{E} \cap \mathcal{F}  ]    \text{, }       
\end{align*}   }

\noindent which is the same lower bound provided in $(\mathrm{FKG})- (\mathcal{S}-\mathrm{SMP})$. The exponents of the number of edges $x$ and the exponential of the first external field are respectively parametrized with respect to the number of edges $e$, and the number of hexagons $\mathcal{H}$ that are not monochromatically colored. On the other hand, under the $\mathrm{PushDual}$ condition, the conditional pushforward under wired boundary conditions supported over $\mathscr{H}$ satisfies,

\begin{align*}
        \mu^1_{\mathscr{H}} [    \widetilde{\mathcal{F}} | \widetilde{\mathcal{E}} \cap \mathcal{F} \cap \mathcal{G}       ]   \geq c^{\text{ }   \lambda  \alpha      \text{ }  \mathrm{Stretch}\text{ } }    \text{, }  
\end{align*}

\noindent which will be used to complete the remaining steps from the \textbf{Lemma} $1$ proof. In particular, the intersection $\{\widetilde{\mathcal{E}} \cap \mathcal{F} \cap \mathcal{G} \}$ can be bounded above by the product of $\lambda$ horizontal crossings below, from $(\mathrm{FKG})$, 

\begin{align*}
           \mu^1_{\mathscr{H}}[  \widetilde{\mathcal{E}}  \cap   \widetilde{\mathcal{F}} \cap   \widetilde{\mathcal{G}}             ]            = \mu^1_{\mathscr{H}} [  \widetilde{\mathcal{F}} |    \widetilde{\mathcal{E}} \cap \mathcal{F} \cap \mathcal{G}    ]        \mu^1_{\mathscr{H}} [     \widetilde{\mathcal{E}} \cap \mathcal{F} \cap \mathcal{G}       ]   \geq c^{\lambda \alpha \text{ }  \mathrm{Stretch}} \mu^1_{\mathscr{H}}[   \widetilde{\mathcal{E}} \cap \mathcal{F} \cap \mathcal{G}      ]   \text{, }     
\end{align*}

\noindent through the same application of $(\mathcal{S}-\mathrm{SMP})$ and (MON), and where $\mathcal{F}$ denotes the crossing event that neither of the three hexagonal boxes defined at the beginning of the proof are vertically crossed. As a result, the last application of $(\mathrm{FKG})$ yields, for $\lambda > 0$ horizontal crossings across thinner hexagons,

\begin{align*}
         \mu^1_{\mathscr{H}} [   \widetilde{\mathcal{E}} \cap \widetilde{\mathcal{F}} \cap \mathcal{G}               ]   \leq   \mu^1_{[0 , \rho N ] \times_H [ 0 , N \lambda \mathrm{Stretch}   ]} \big[     \mathcal{H}_1      \big]^{\lambda + 1}    \text{, } 
\end{align*}

\noindent under free boundary conditions.

Comparing the pushforward under free boundary conditions to the pushforward under wired boundary conditions yields, after taking the same infinite aspect length limit as obtained in Section \textit{7}. From previous applications of $(\mathcal{S} -\mathrm{SMP})$ and (MON) are used, in order to suitably compare boundary conditions, as a consequence imply a similar estimate as obtained in Section \textit{7}, 

\begin{align*}
       \mu^0_{[0 , \rho N] \times_H [0 ,  n_2 \lambda \mathrm{Stretch} ]}  \big[   \mathcal{H}_{[0 ,  \rho N  ] \times_H [ 0  ,   n_1 \lambda \mathrm{Stretch}   ]}     \big]^{\lambda}   \overset{\rho  \longrightarrow + \infty }{\approx}  \big( p_{\mathrm{Stretch}   \text{ } N} \big)^{\lambda}  \text{,                                       } 
    \end{align*}

    \noindent which, as in previous arguments for \textit{Lemma} $1$, is bounded below by the following infinite aspect ratio limit,
    
    \begin{align*}
   \mu^1_{[0 ,  \rho N  ] \times_H [  0 , n_1 \lambda  \mathrm{Stretch}  ] } \big[   \mathcal{V}^{\text{ }c}_{[ 0 ,  \rho N ] \times_H [  0 ,   n_1 \lambda \mathrm{Stretch}     ]}  \big]^{\lambda + 1}      \overset{\rho  \longrightarrow + \infty }{\approx}    \frac{ 1 }{\lambda^C} \big( q_{\mathrm{Stretch}\text{ } N}     \big)^{\lambda+1}         \text{, }   
\end{align*}

\noindent with the exception that the support of the measure with free (-) boundary conditions is over a hexagon with thinner aspect length. The result corresponds to the renormalization inequality for the horizontal crossing probability, concluding the argument under $(\mathrm{PushDual})$. Below, we briefly describe how the same sequence of inequalities applies for the remaining possibility.

Suppose that $(\mathrm{PushPrimal})$ holds. Under this assumption, denote $\widetilde{\mathcal{F}}$ as the crossing event that none of the boxes $\widetilde{\mathscr{H}_i}^{\pm}$ are vertically crossed. From this event, the assumption implies from the definition of the horizontal and vertical strip densities for the Spin Measure that the arguments to bound the conditional probability can be achieved by the same line of argument, possibly with larger $C$. \boxed{}

\section{Quadrichotomy proof}

In the final section we classify all possible behaviors of the model. Briefly, we remark that for the continuous critical case, the first part of the argument does not require use of (SMP) and (CBC) for original results in the Random-Cluster model, implying that the entirety of the argument immediately applies. Briefly, we summarize the steps of the argument. We consider horizontal crossing events across a regular hexagon, pushed forwards under free boundary conditions for the Spin measure supported over a slightly larger hexagon. From knowledge of the longest edge in the $+$ path of the horizontal crossing, removing the largest edge from the configuration easily yields a connectivity event along the common set of edges over a subgraph of the triangular lattice that excludes the length of the maximal edge along two points $x$ and $y$. These steps demonstrate the ingredients for the \textit{{Discontinuous Critical}} case, before obtaining the horizontal strip densities in the infinite aspect length limit as $\rho \longrightarrow + \infty$. For the discontinuous critical case, the second part of the argument requires use of $(\mathcal{S} -\mathrm{CBC})$ and (MON). Before proceeding, we cite the following theorem which classifies the probability of obtaining loop configurations of fixed length in the model. In the following application of the inequalities, we choose an aspect ratio of hexagons dependent on $\rho$, from which horizontal and vertical crossings will be studied.

\bigskip

\noindent \textbf{Theorem} $3$ ({\textit{Continuous \& \textit{Discontinuous critical} cases}, behaviors of the dilute Potts model quadrichotomy}, {[8]}): For configurations distributed under the $\mu^{\tau}_{G,x,n} \big[ \cdot \big]$, for $n \geq 1$ and $x \leq \frac{1}{\sqrt{n}}$, one of the two possible behaviors occurs,

\begin{align*}
      \mu^{\tau}_{\mathcal{S}} [   R \geq k      ]  \leq \mathrm{exp}(-ck)         \tag{Exponential decay of + paths}     \text{, } 
\end{align*}

\noindent where $R$ is the diameter of the largest loop surrounding the origin, demonstrating that $+$ paths are exponentially unlikely for any $k \geq 1$, or,

{\small \begin{align*}
     c \leq   \mu^{\tau}_{\mathcal{S}}[    \textit{There exits a + path which horizontally crosses a rectangle over the triangular lattice }          ]  \\  \leq 1 - c     \tag{RSW box-crossing property}  \text{, } 
\end{align*}   } 

\noindent demonstrating that the RSW box-crossing property is satisfied, under the admissibility condition,

{\small \begin{align*}
  \big\{  \tau \textit{ is admissible for the Russo-Seymour-Welsh property} \big\} \Longleftrightarrow \big\{ \tau \textit{ is placed sufficiently far away} \\ \textit{  from $\mathcal{H}$ such that }  \mu^{\tau}_{G,x,n}[\mathcal{H}] \neq 0 \big\} , 
\end{align*} }

\noindent on boundary conditions $\tau$. Each possibility holds for $\tau \in \{\pm 1 \}^{F (\textbf{H}^*)}$ and $c > 0$.

\bigskip

Observe that we have slightly rephrased the first condition provided in {[8]} which is stated instead for the loop measure $\textbf{P}^{\xi}_{G,x,n} \big[ \cdot \big]$ defined in Section \textit{3.1}. The equivalent condition of obtaining a loop configuration whose largest diameter about the origin is $k$ is equivalent to obtaining a path of $+$ spins about the origin. From the statement of \textbf{Theorem} $3$, we now study \textbf{Theorem} $2$.

\subsection{Subcritical $\&$ Supercritical behaviors}

\noindent \textit{Proof of Theorem $2$} (\textit{{Discontinuous Critical}} phase from non \textit{{Subpercritical}} phase). As mentioned at the beginning of the section, first suppose that the first possibility holds. To show that this condition implies that the phase transition is discontinuous, consider the following. Define a horizontal crossing across $\mathscr{H}$. From the existence of such an event, the longest edge in the crossing of arbitrary length $L$ then excluding the length of this longest edge from the crossing implies that another closely related crossing event occurs across a subgraph of the triangular lattice which excludes the maximal edge with length $L$. Hence there exists vertices in a subgraph of the triangular lattice, such that the vertices $x$ and $y$ are connected by a $+$ path in a hexagon of smaller aspect length that is not regular. Collecting these observations implies the following, where the upper and lower bounds of the inequality are taken under $-$ boundary conditions, by the union bound,

\begin{align*}
        \mu^0_{[0 , \rho N ] \times_H [   0  , 2N   ] }[ \mathcal{H}_{[0 , \rho N ] \times_H [0 ,  N   ]}    ]   \leq  c      N^2 \mu^0_{[0 , \rho N ] \times_H [   0 ,  2 N  ] }\text{ } \big[   x \overset{ [0  ,\rho N ] \times_H [ 0  , N ]  }{\longleftrightarrow} y    \big] \text{, } 
\end{align*}    

\noindent where $x$ and $y$ are the vertices, with $c$ an arbitrary positive constant. For the next step, we introduce horizontal translates of $x$ with $x_k = x + (4k N , 0)$. Across all horizontal translates of $x$, yields the following lower bound for the connectivity event between $x$ and each $x_k$, by (MON) and (FKG),

\begin{align*}
     \mu^0_{\mathcal{S}}[   x  {\longleftrightarrow} x_k     ]      \geq \mu^0_{[0 ,  \rho N ] \times_H [ 0 , 2 N]}  \big[  x   \overset{ [0 , \rho N ] \times_H [ 0  , N   ]}{\longleftrightarrow} y    \big] \text{ }   \text{. } 
\end{align*}

\noindent From previous remarks, the first upper bound given in the proof dependent on $c$ yields the inequality, as applied in $(\mathrm{FKG})$ several times previously in the argument,

\begin{align*}
         \mu^0_{\mathcal{S}}[  x {\longleftrightarrow} x_k    ] \geq \frac{1}{cN^2}  \mu^0_{[0 ,  \rho N] \times_H [0 , 2N]}\text{ }  \big[  \mathcal{H}_{[0 , \rho N] \times_H [ 0 ,  N]}      \big]^{2k} \text{, } 
\end{align*}

\noindent from which taking the infinite limit as in previous arguments implies, for $k \longrightarrow + \infty$,

\begin{align*}
     p_{2N}^2   \geq \frac{1}{cN^2} \mu^0_{[0 , \rho N ] \times_H [0 , 2 N]} \text{ }  \big[    \mathcal{H}_{[0 , \rho N ] \times_H [0 , N]}  \big] \text{ }    \text{, }
\end{align*}    

\noindent so that the pushforward of the spin measure under free boundary conditions satisfies the strip density estimate from the original definition provided in the beginning of Section \textit{7}, from the connected components of $+$ paths from the occurrence of $\{x \longleftrightarrow + \infty \}$. Finally, we observe that the upper bound for the horizontal strip density decays exponentially fast, implying that the pushforward in the lower bound taken under free boundary conditions does as well. As expected, to analyze the other possibility for infinitely long vertical crossings, repeating the same steps of the argument, with the exception that the horizontal crossing event is instead a vertical crossing event, simply yields a similar bound, from an application of \textbf{Lemma} $12$ for some integer $\lambda$ satisfying the conditions of the statement, that the probability of obtaining an infinitely long vertical crossings is an upper bound of the following inequality, 

\begin{align*}
       q_{2N}^2 \geq  \frac{1}{cN^2} \big(1 - \mu^1_{[0 , \rho N] \times_H [0 , 2N]} \big[ \mathcal{V}^{\text{ } c}_{[0 ,  \rho N ] \times_H [ 0 ,  N  ]}  \big]  \big)   \text{, } 
\end{align*}

\noindent which nevertheless still exponential decays for the same reason as $k \longrightarrow + \infty$. This conclude the argument for the model, demonstrating that the horizontal and vertical strip densities hold for infinite aspect ratios. \boxed{}

\bigskip

\noindent In the following, we analyze the {\textit{Continuous critical case}} to obtain RSW results. We refer the attention of the reader to the specific arrangement of the finite volumes $\mathscr{R}_1, \mathscr{R}_2, \mathscr{R}_3, \mathscr{R}_4$ in \textbf{Figure} \textit{12} above.

\bigskip

\noindent $\textit{Proof of Theorem}$ $2$ (\textit{{Continuous Critical}} phase from non \textit{{Subcritical}} phase). To prove RSW results, consider the following four finite volumes, with $\mathscr{R}_1 \subsetneq  \mathscr{R}_2$, and $\rho , N > 0$,

\begin{align*}
      \mathscr{R}_1 \equiv [  0 ,  \rho N ] \times_H [  0, N ] \text{, } \\ \\ 
  \mathscr{R}_2 \equiv [  - N ,  \rho N + N ] \times_H [  - N , 2 N  ] \text{, } \\   \\
\mathscr{R}_3 \equiv  \big[ - \frac{2N}{3},  - \frac{N}{3} \big] \times_H [  - N,  2 N]  \text{ , }  \\    \\
 \mathscr{R}_4 \equiv    \big[  \rho N + \frac{N}{3}  ,  \rho N + \frac{2N}{3} \big] \times_H [ -N  ,2 N   ]    \text{. } 
\end{align*}

\begin{figure}
\begin{center}
\begin{align*}
     \includegraphics[width=0.9\columnwidth]{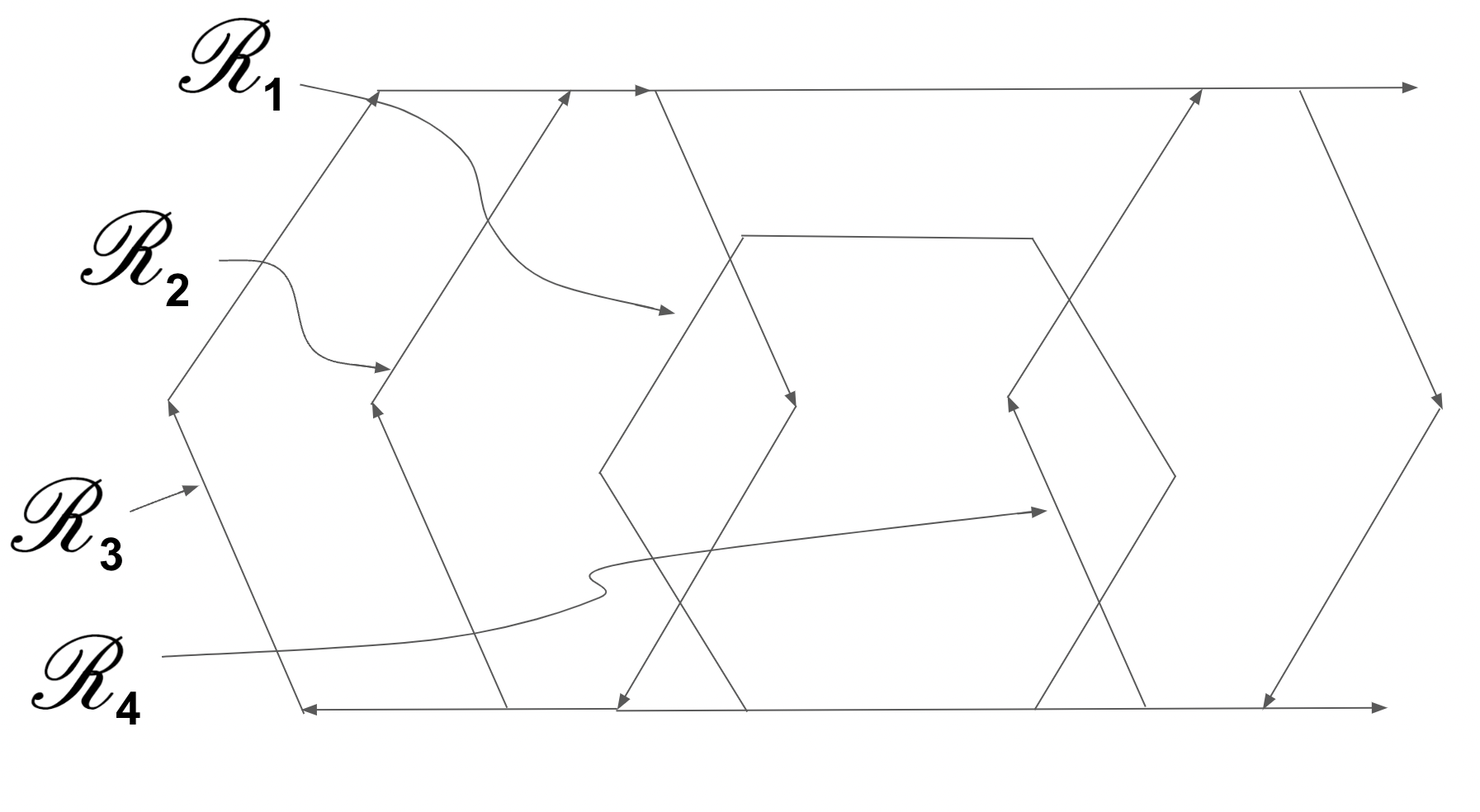}
\end{align*}
\caption{The configuration of finite volumes across which crossing events are quantified to obtain RSW results in the proof of the remaining case of \textbf{Theorem} $2$.}
\end{center}
\end{figure}

\noindent We bound the crossing probability of a the horizontal crossing event below with the product of vertical and horizontal strip densities (from \textbf{Definition} $1$). To ensure that the appropriate cancellation of crossing events occurs, we bound the probability of a horizontal crossing event with the horizontal strip density $p_n$, in which, 

\begin{align*}
        \mu^{\mathrm{Mixed}}_{\mathscr{R}_2} [ \mathcal{H}_{\mathscr{R}_1}  ]    \text{ }  \geq \text{ } p^{\alpha \rho}_N     \text{, }   
\end{align*}    

\noindent for every $\alpha \geq 1$ and $\mathrm{Mixed}$ boundary conditions (all vertices along the two edges of the hexagon on the left and right are wired), due to the fact that,

\begin{align*}
        \mu^{\mathrm{Mixed}}_{\mathscr{R}^{\alpha}_2} [   \mathcal{H}_{\mathscr{R}^{\alpha}_1}    ]      \geq   \mu^{\mathrm{Mixed}}_{\mathscr{R}_2} [ \mathcal{H}_{\mathscr{R}_1}    ]^{\alpha}       \geq  \mu^0_{[0, \rho N ] \times_H [ 0,   \lambda \mathrm{Stretch}]} \big[  \mathcal{H}_{[0   , \rho N ] \times_H [ 0 ,\lambda \mathrm{Stretch}  ]}   \big]   \text{, } 
\end{align*}

\noindent for $\mathscr{R}^{\alpha}_2 \equiv \big[- \alpha N , \alpha \rho N + N \big] \times_H \big[- N , 2 N \big] $, and $\mathscr{R}^{\alpha}_1 \equiv \big[ 0 ,  \alpha \rho N \big] \times_H \big[ 0, N \big]$, raising the inequality to $\frac{1}{\rho}$ as $\rho \longrightarrow + \infty$. On the other hand, below we bound the probability of the following intersection of crossing events, conditionally on $\mathcal{H}_{\mathscr{R}_1}$,

\begin{align*}
     \mu^{\mathrm{Mixed}}_{\mathscr{R}_2} \big[      \mathcal{H}_{\mathscr{R}_3} \cap \mathcal{H}_{\mathscr{R}_4}       \big|     \mathcal{H}_{\mathscr{R}_1}   \big]        \text{, } 
\end{align*}

\noindent which can similarly be lower bounded with the vertical crossing strip density, as, 

\begin{align*}
      \mu^{\mathrm{Mixed}}_{\mathscr{R}_2} \big[      \mathcal{H}_{\mathscr{R}_3} \cap \mathcal{H}_{\mathscr{R}_4}        \big|            \mathcal{H}_{\mathscr{R}_1}  \big]    \geq     q^{\gamma}_{\frac{N}{3}}      \text{ }     \text{, } 
\end{align*}        

\noindent for $\gamma$ sufficiently large, and as the inequality is raised to $\frac{1}{\rho}$ as $\rho \longrightarrow + \infty$. The two inequalities can be bounded below by the product of horizontal and vertical strip densities,

\begin{align*}
\mu^{\mathrm{Mixed}}_{\mathscr{R}_2} \big[   \mathcal{H}_{\mathscr{R}_1}   \big|  \mathcal{H}_{\mathscr{R}_3}  \cap \mathcal{H}_{\mathscr{R}_4}      \big]                       \geq \mu^{\mathrm{Mixed}}_{\mathscr{R}_2} \big[   \mathcal{H}_{\mathscr{R}_1}  \cap \mathcal{H}_{\mathscr{R}_3} \cap    \mathcal{H}_{\mathscr{R}_4}           \big]  \text{ }           \geq    p^{\alpha \rho}_N q^{\gamma}_{\frac{N}{3}} \text{ }     \text{. }    
\end{align*}

\noindent Finally, the ultimate term in the inequality above can be bounded above with the horizontal crossing event of interest, as the conditioning on the other two horizontal crossing events is removed,

\begin{align*}
     \text{ }  \mu^{\mathrm{Mixed}}_{\mathscr{R}_2}[     \mathcal{H}_{\mathscr{R}_1}   ]    \text{ }  \geq     \mu^{\mathrm{Mixed}}_{\mathscr{R}_2} \big[  \mathcal{H}_{\mathscr{R}_1}  \big| \mathcal{H}_{\mathscr{R}_3} \cap  \mathcal{H}_{\mathscr{R}_4}    \big]     \text{. }   
\end{align*}

\noindent Altogether, by duality and rotational invariance of $\mu$, one obtains,

\begin{align*}
 \text{ }   \mu^{\mathrm{Mixed}^{\prime}}_{\mathscr{R}_2} [    \mathcal{H}_{\mathscr{R}_1}        ]     \leq    1 - p^{\alpha \rho}_N q^{\gamma}_{\frac{N}{3}} \text{, } 
\end{align*}

\noindent where $\mathrm{Mixed}^{\prime}$ is a rotation of $\mathrm{Mixed}$, implying that the RSW inequality for horizontal crossing events is, with respective constants $c$ and $1-c$ provided in the lower and upper bounds,

\begin{align*}
  \text{ }    p^{\alpha \rho}_N q^{\gamma}_{\frac{N}{3}}    \leq  \mu^{\xi}_{\mathscr{R}_2}[    \mathcal{H}_{\mathscr{R}_1}    ]  \leq   1 - p^{\alpha \rho}_N q^{\gamma}_{\frac{N}{3}} \text{, } 
\end{align*}

\noindent independently of boundary conditions $\xi$, hence concluding the proof because the  \textit{{Continuous Critical}} phase occurs. \boxed{}

\bigskip

\subsection{Applications of different phases of the quadrichotomy}

\subsubsection{Subcritical regime}.

\noindent \textit{Proposition A (coexistence of wired and free high-temperature measures)}. The probability of connectivity to distance $n$ under the wired high-temperature Loop $O(n)$ measure is exponentially upper bounded, as

\begin{align*}
 \mu^1_{\mathscr{H}_N} [ 0 \longleftrightarrow \partial \mathscr{H}_N  ] \leq \mathrm{exp}(-c N) \text{, }    
\end{align*}

\noindent also implying that $\mu^0 = \mu^1$.

\bigskip

\noindent \textit{Proof of Proposition A}. To show that $\mu^0 = \mu^1$ in \textit{{Subcritical}}, we lower bound the horizontal crossing probability in terms of crossings of loops to the left and right sides, respectively $\mathcal{L}_N$ and $\mathcal{R}_N$ of $\mathscr{H}_N$, in which,

\begin{align*}
  \mu^{1}_{\mathscr{H}_{2N}}[     \mathcal{H}_{\mathscr{H}_N}  ]            \geq \mu^1_{\mathscr{H}_{2N}}[   0 \overset{\mathscr{H}_N}{\longleftrightarrow}  \text{ } \mathcal{L}_N   ] mu^1_{\mathscr{H}_{2N}}[   0 \overset{\mathscr{H}_N}{\longleftrightarrow}  \mathcal{R}_n ]   \geq \frac{1}{32}  \mu^1_{\mathscr{H}_{2N}}[   0 \longleftrightarrow  \partial \mathscr{H}_N       ]   \text{. } 
\end{align*}

\noindent As a consequence, by $(\mathcal{S}- \mathrm{SMP})$ the probability of connectivity to distance $n$ decays exponentially fast. The fact that the high-temperature measure under wired and free boundary conditions coincides follows from classical arguments. \boxed{}

\subsubsection{Supercritical regime}.

\noindent \textit{Proposition B (exponential unlikelilihood of obtaining finite connected components)}. In \textit{{Supercritical}}, the probability of an infinite connected component under free boundary conditions of being absent is exponentially unlikely, as

\begin{align*}
  \mu^0_{\mathscr{H}_N} [    \mathscr{H}_N \not\longleftrightarrow \infty       ] \leq \mathrm{exp}(-c N ) \text{. }  
\end{align*}

\bigskip

\noindent \textit{Proof of Proposition B}. By duality, if an infinite connected component does not exist in the primal configuration, then an infinite connected component exists in the dual configuration that is measurable over $\textbf{T}$. With the dual configuration, there does exist a loop whose maximum diameter is $n$, implying that the standard connectivity event at distance $n$ does occur. Hence the desired inequality follows. \boxed{}





\section{References}

\noindent[1] Beffara, V. \& Duminil-Copin, H. The self-dual point of the two-dimensional random-cluster model is critical for $q \geq 1$. \textit{Probability Theory and Related Fields} \textbf{153}: 511-542 (2012). https://link.springer.com/article/10.1007/s00440-011-0353-8.

\bigskip

\noindent[2] Beliaev, D., Muirhead, S. \& Wigman, I. Russo-Seymour-Welsh estimates for the Kostlan ensemble of random polynomials. \textit{Annales de l'Institut Henri Poincare} \textbf{57}(4): 2189-2218 (2017). https://doi.org/10.1214/20-AIHP1142.

\bigskip

\noindent[3] Crawford, N., Glazman, A., Harel, M. \& Peled, R. Macroscopic loops in the loop $O(n)$ model via the XOR trick. Ann. Probab. 53(2): 478-508 (March 2025). DOI: 10.1214/24-AOP1712.

\bigskip

\noindent[4] Clarence, J., Damasco, G., Frettloh, D.,. Loquias, M. Highly Symmetric Fundamental Domains for Lattices in $\textbf{R}^2$ and $\textbf{R}^3$. \textit{Arxiv} (2018).

\bigskip

\noindent[5] Duminil-Copin, H. Parafermionic observables and their applications to planar statistical physics models. \textit{Ensaois Matematicos} \textbf{25}: 1-371 (2013). https://www.unige.ch/~duminil/publi/parafermion.pdf.

\bigskip

\noindent[6] Duminil-Copin, H. Parafermionic observables and their applications. \textit{AMP Bulletin} (2015). https://www.unige.ch/~duminil/publi/2015IAMPbulletin.pdf.

\bigskip

\noindent[7] Duminil-Copin, H. Sharp threshold phenomena in statistical physics. \textit{Japanese Journal of Mathematics} \textbf{14}: 1-25 (2019). https://www.unige.ch/~duminil/publi/2018Takagi.pdf.

\bigskip

\noindent[8] Duminil-Copin, H., Glazman, A., Peled, R. \& Spinka, Y. Macroscopic Loops in the Loop $O(n)$ model at Nienhuis' Critical Point. \textit{J. Eur. Math. Soc.} \textbf{23}: 315-347 (2021).

\bigskip

\noindent[9] Duminil-Copin, H., Hongler, C. \& Nolin, P. Connection probabilities and RSW-type bounds for the FK Ising Model. \textit{Communications on Pure and Applied Mathematics} \textbf{64}(9): 1165-1198 (2011). https://www.unige.ch/~duminil/publi/RussoSeymourWelsh.pdf.

\bigskip

\noindent[10] Duminil-Copin, H., Manolescu, I. \& Tassion, V. Planar random-cluster model: fractal properties of the critical phase. \textit{Probability Theory and Related Fields} \textbf{181}: 401-449 (2020).  https://www.unige.ch/~duminil/publi/2020criticalFK.pdf.

\bigskip

\noindent[11] Duminil-Copin, H. \& Smirnov, S. Conformal invariance of lattice models. \textit{Lecture Notes} (2012). https://www.unige.ch/~smirnov/papers/clay-j.pdf.

\bigskip

\noindent[12] Duminil-Copin, H. \& Smirnov, S. The convective constant of the honeycomb lattice equals $\sqrt{2+\sqrt{2}}$. \textit{Annals of Mathematics} \textbf{175}(3): 1653-1665 (2012). https://www.unige.ch/~smirnov/papers/saw-j.pdf.

\bigskip

\noindent [13] Duminil-Copin, H., Sidoravicius, V. \& Tassion, V. Continuity of the phase transition for planar random-cluster and Potts models for $1 \leq q \leq 4$. \textit{Communications in Mathematical Physics} \textbf{349}: 47-107 (2017). https://www.unige.ch/~duminil/publi/pottsPhaseTransition.pdf.

\bigskip

\noindent[14] Duminil-Copin, H. \& Tassion, V. Renormalization of crossing probabilities in the planar random-cluster model. \textit{Moscow Mathematical Journal} \textit{20} (4): 711-740. https://www.unige.ch/~duminil/publi/2019renormalizationCrossings.pdf.

\bigskip

\noindent[15] Feher, G. \& Nienhuis, B. Currents in the dilute $O(n=1)$ model. \textit{Arxiv} 1510.02721v2 (2018).

\bigskip

\noindent[16] Fradkin, E. Disorder Operators and their Descendants. \textit{Lecture Notes}. https://doi.org/10.1007/s10955-017-1737-7.

\bigskip

\noindent[17] Glazman, A. \& Manolescu, I. Uniform Lipschitz functions on the triangular lattice have logarithmic variations. \textit{Communications in Mathematical Physics} \textbf{381}: 1153-1221 (2021). https://doi.org/10.1007/s00220-020-03920-z.

\bigskip

\noindent[18] Gheissari, R. \& Lubetzky, E. Quasi-polynomial Mixing of Critical two-dimensional Random Cluster Models.  \textit{Random Structures \& Algorithms} \textbf{56}(2) (2019). $\mathrm{https://cims.nyu.edu/~eyal/papers/FK\text{ }  mixing.pdf}$.

\bigskip

\noindent[19] Guo, W., Blote, H. \& Nienhuis, B. First and Second Order Transitions in the Dilute $O(n)$ models \textit{International Journal of Modern Physics C} \textbf{10}(1): 291-300 (1999). https://www.worldscientific.com/doi/abs/10.1142/S012918319900019X.

\bigskip

\noindent[20] Guo, W., Blote, H. \& Nienhuis, B. Phase Diagram of a Loop on the Square Lattice. \textit{International Journal of Modern Physics C} \textbf{10}(1): 301-308 (1999). https://www.worldscientific.com/doi/abs/10.1142/S0129183199000206.

\bigskip

\noindent[21] Hongler, C. Percolation on the triangular lattice (2007). https://hongler.org/lattice-models/supp-mat/triangular-lattice.pdf.

\bigskip

\noindent[22] Nienhuis, B. \& Guo, W. Tricritical $O(n)$ models in two dimensions. \textit{American Physical Society} \textbf{78}: 061104 (2008). https://doi.org/10.1103/PhysRevE.78.061104.

\bigskip

\noindent[23] Russo, L. A note on percolation. \textit{Zeitschrift fur Wahrscheinlichkeitstheorie und Verwandte Gebiete} \textbf{43} 39-48 (1978). https://scispace.com/pdf/a-note-on-percolation-1eh2pg6jgv.pdf.

\bigskip

\noindent[24] Seymour, P. \& Welsh, D. Percolation probabilities on the square lattice. \textit{Annals
Discrete Math}, \textbf{3}: 227–245 (1978). https://cse.buffalo.edu/~regan/DJAW/SeymourWelsh1978.pdf.

\bigskip

\noindent[25] Smirnov, S. Discrete Complex Analysis and Probability. \textit{Proceedings of the International Congress of Mathematicians} (2010). https://www.unige.ch/~smirnov/papers/icm-2010.pdf.

\bigskip

\noindent[26] Smirnov, S. Conformal invariance in random cluster models. I. Holomorphic fermions in the Ising model. \textit{Annals of Mathematics} \textbf{172}: 1435-1467 (2010). https://doi.org/http://doi.org/10.4007/annals.2010.172.1435.

\bigskip

\noindent[27] Tassion, V. Crossing probabilities for Voronoi percolation. \textit{The Annals of Probability} \textbf{44}(5): 3385-3398 (2016). DOI: 10.1214/15-AOP1052.

\bigskip

\noindent[28] Zeng, X. A Russo Seymour Welsh Theorem for critical site percolation on $\textbf{Z}^2$. \textit{Arxiv} (2013).

\bigskip

\noindent [29] Grimmett, G. \textit{Percolation}, Volume 321 of \textit{Fundamental Principles of Mathematical Sciences}. Springer-Verlag, second edition (1999).
 

\end{document}